\documentclass{article}
\usepackage[utf8]{inputenc}

\title{Projection onto quadratic hypersurfaces}
\author{Loïc Van Hoorebeeck$^a$, P.-A. Absil$^a$, Anthony Papavasiliou$^b$}

\date{%
\small
    $^a$ICTEAM, UCLouvain, Louvain-la-Neuve, Belgium.\\%
    $^b$CORE, UCLouvain, Louvain-la-Neuve, Belgium.\\[2ex]%
    \today
}

\usepackage[dvipsnames]{xcolor}
\usepackage[utf8]{inputenc}

\usepackage{mathtools}

\newtagform{parenthesest}{(}{$)_t$}
\usepackage{amssymb}
\usepackage{algorithm}
\usepackage{algorithmic}
\usepackage{bbold}
\usepackage{pgfplots}
\usepackage{siunitx}
\usepackage{graphicx}
\usepackage{lmodern}
\usepackage{amsthm} 
\usepackage{bm}
\usepackage{bbm}
\usepackage[fleqn]{nccmath}
\usepackage{hyperref}
\usepackage[capitalise]{cleveref}
\usepackage{footnote}
\usepackage{mathrsfs}
\makesavenoteenv{tabular}
\makesavenoteenv{table}
\usepackage{mathabx}
\usepackage{caption}
\usepackage{subcaption}

\usepackage{tikz}

\usetikzlibrary{fit,backgrounds}

\pgfdeclarelayer{background}
\pgfdeclarelayer{foreground}

\pgfplotsset{compat=newest}
\usepgfplotslibrary{groupplots}
\usepgfplotslibrary{dateplot}
\colorlet{myred}{red!60}
\colorlet{mygreen}{ForestGreen}
\colorlet{myblue}{RoyalBlue}
\definecolor{colorSol}{RGB}{206, 68, 216}

\usepackage{arrayjobx}
\usepackage[comparestr,comparenum,randompart]{arraysort}
\usepackage{booktabs}

\usepackage{lipsum}

\theoremstyle{plain}
\newtheorem{theorem}{Theorem}[section]

\newtheorem{proposition}[theorem]{Proposition}

\theoremstyle{definition}
\newtheorem{definition}{Definition}[section]

\theoremstyle{remark}

\newcommand{\alt}{\textrm{APC}}
\newcommand{\A}{\bm{B}}
\newcommand{\Ab}{\bm{b}}
\newcommand{\Ac}{c}
\newcommand{\argmin}{\operatornamewithlimits{\arg\,\min}}
\newcommand{\B}{\bm{B}}

\newcommand{\D}{\bm D}
\renewcommand{\d}{\bm d}

\newcommand{\dalt}{d_{\exact}}
\newcommand{\dexact}{d_{\alt}}

\newcommand{\eg}{\emph{e.g.}}
\newcommand{\exact}{\textrm{Gur}}

\newcommand{\I}{I}
\newcommand{\Ip}{I^+}
\newcommand{\In}{I^-}
\newcommand{\ie}{\emph{i.e.}}

\newcommand{\muI}{\mu^{\mathrm{I}}}

\newcommand{\norm}[1]{\left | \left | #1 \right | \right |}

\newcommand{\niter}{n_{\text{iter}}}

\renewcommand{\P}{(\ref{eq:P})}

\newcommand{\bp}[1]{\bm{p}}

\newcommand{\Q}{\mathcal{Q}}
\newcommand{\Qtot}{\mathcal{Q}^{\text{tot}}}
\newcommand\numberthis{\addtocounter{equation}{1}\tag{\theequation}}

\newcommand{\Rn}{\mathbb{R}^{n}}

\newcommand{\size}[1]{\left | #1 \right |}

\newcommand{\spec}[1]{\textrm{spec}(#1)}
\newcommand{\st}{\left | \right .}

\newcommand{\talt}{t_{\alt}}
\newcommand{\texact}{t_{\exact}}

\newcommand{\Tr}[1]{{#1}^{\intercal}}

\newcommand{\x}{\bm x}
\newcommand{\xmin}{\underline{x}}
\newcommand{\xmax}{\overline{x}}

\newcommand{\abs}[1]{\left | #1 \right |}

\renewcommand{\Pr}{\textrm{Pr}}
\DeclareMathOperator{\prox}{prox}
\DeclareMathOperator{\argmax}{argmax}

\crefformat{equation}{(#2#1#3)}
\crefformat{figure}{Figure~#2#1#3}

\begin{document}

\maketitle

\abstract{We address the problem of projecting a point onto a quadratic hypersurface, more specifically a central quadric. We show how this problem reduces to finding a given root of a scalar-valued nonlinear function. We completely characterize one of the optimal solutions of the projection as either the unique root of this nonlinear function on a given interval, or as a point that belongs to a finite set of computable solutions. We then leverage this projection and the recent advancements in splitting methods to compute the projection onto the intersection of a box and a quadratic hypersurface with alternating projections and Douglas-Rachford splitting methods. We test these methods on a practical problem from the power systems literature, and show that they outperform IPOPT and Gurobi in terms of objective, execution time and feasibility of the solution.


}



\section{Introduction}%
\label{sec:Introduction}
This paper discusses the projection of a given point onto a nonsingular quadratic hypersurface, or nonsingular \emph{quadric}. Quadrics are a natural generalization of hyperplanes. The projection onto a quadric appears, e.g., in splitting algorithms for the projection between a quadric---or the Cartesian product of quadrics---and a polytope. This problem has direct applications, \eg, in power systems~\cite{vh22}, where the power losses can be approximated as a quadratic hypersurface. Numerical experiments on such problems are developed in this paper. Other applications of quadratic projections emerge in the context of the security region of gas networks~\cite{song_security_2021}, or in local learning methods~\cite{scott_thesis_2020}.

It is therefore intriguing that few studies of this problem can be found in the literature. Indeed, projections onto quadratic surfaces have been studied for the 2D and 3D cases, see \eg,~\cite{morera_distance_2013,lott_iii_direct_2014,huang_shape_2020}. However, to the best of our knowledge, the extension to an arbitrary dimension has not been pursued, with the exception of the short discussion at the end of~\cite{lott_iii_direct_2014} and in~\cite{sosa_algorithm_2020}.
Although the method proposed in~\cite{sosa_algorithm_2020} can handle the singular case, \ie, the case where the matrix that defines the quadratic surface is singular, it does not always return the exact projection. Moreover, the two-level iterative scheme that is proposed in~\cite{sosa_algorithm_2020} can be computationally expensive.

The structure of this paper is twofold. Firstly, we tackle the problem of projecting onto an $n$-dimensional nonsingular quadric. Secondly, we leverage this projection in the context of splitting methods.

The projection considered in the first part of this paper (\cref{sec:Projection_onto_a_quadric}) is not unique in general, due to the nonconvexity of the feasible set. This implies that we cannot rely on first or second-order methods, since such methods may converge to a local minimum. This projection can also be handled by black-box (commercial) solvers, \eg, Gurobi or IPOPT~\cite{gurobi,wachter06}, however these methods suffer from two main problems: i) the execution time rockets when the dimension of the problem increases to mid or large-scale sizes---this phenomenon is present in our numerical results---and ii) Gurobi is not a local method, and does not exploit the local structure of the problem, even if in certain applications, the starting point is close to the feasible set. The first problem is highlighted by the power system application that is considered here, where the projection step is only a small part of the overall procedure, which renders execution time an important factor in our analysis. 

 Using the Lagrange multiplier technique, we reduce this quadratically constrained quadratic program (QCQP) to the problem of finding the roots of a nonlinear real function. Then, we completely characterize the solutions of the nonconvex projection, and compute one of the solutions of this projection as either the (unique) root of a scalar function on a computable interval, or among a finite set of closed-form solutions. We also show how to find this root using Newton's method, which guarantees a quadratic convergence. Thus, the proposed method provides an efficient way to obtain the exact projection onto a nonsingular quadric. Finally, to further reduce execution time, we also introduce a heuristic based on a geometric construction. This allows us to quickly map a point to the quadric. We detail two variants of this heuristic.


We note that our proposed approach for projecting onto a quadric is not unusual. For example, \cite[\S 6.2.1]{golub_matrix_2013} uses a similar construction for the problem of least squares minimization over a sphere. However, this problem is easier to tackle than ours, since the (unique) solution of this convex problem is the (unique) root of the \emph{secular equation} defined by the KKT conditions. In~\cite[\S 7.3]{conn_trust_2000}, the authors also use a similar procedure, and taxonomy of secular equations, for finding the $\ell_2$-norm model minimizer. However, while the discussion is analogous to what is proposed in this paper, \ie, searching for a specific root of a given scalar-valued nonlinear function on a specific domain, the domain and the function are different in our work. Moreover, our discussion on degenerate cases is not present in~\cite{conn_trust_2000}, since such cases do not appear in the problem that the authors tackle, which is linked to the trust-region subproblem.

In the second part of this paper (\cref{sec:Splitting_methods}) we test our method of projecting onto a quadric in order to solve the problem of projecting onto the intersection between a polytope and a quadric. Both projections can be easily computed: the first one is even trivial if we consider a box, and the second one is efficiently obtained with the method proposed in the first part of the paper. We then leverage the rich literature on splitting methods for nonconvex programming. We consider, in this work, two splitting methods: the alternating projections and Douglas-Rachford splitting. References for these schemes can be found in~\cite{attouch_proximal_2010,bauschke_phase_2002,drusvyatskiy_transversality_2015,lewis_alternating_2008,lewis_local_2009}, and in~\cite{bauschke_phase_2002,li_douglasrachford_2016}, respectively.

Depending on the splitting method considered, and whether or not we use exact projection onto the quadric or a heuristic, we detail five different methods for projecting a point onto the intersection of a box and a quadric. We analyse these methods in~\cref{sec:Numerical_experiments}. The five methods are benchmarked against IPOPT in the ellipsoidal and hyperboloidal case, both for small and large-scale problems. In these experiments, we observe that one of the proposed methods, namely the alternating projections with exact projections, attains the best objective. We also observe that the alternating projections used with one of the heuristics (the gradient-based heuristic) reaches competitive objectives in a reduced amount of run time. All the methods considered outperform IPOPT in terms of execution time, with a difference of several orders of magnitude. Finally, we benchmark one of the proposed methods against Gurobi. We use Gurobi in order to find the optimal solution for a problem inspired from the power systems literature. Since, in this context, the starting point is close to the feasible set, our proposed method clearly outperforms Gurobi, both in terms of execution time and objective. Using the lower bound computed by Gurobi, we can also conclude that, in the context of this specific problem, the proposed method finds the optimal solution, even if there is no guarantee for finding it in general.

\section{Projection onto a quadric}%
\label{sec:Projection_onto_a_quadric}


\subsection{Problem formulation}%
\label{sub:Problem_formulation}

In this section, the problem of interest is introduced. This problem consists in the projection of a given point $\tilde{\x}^0 \in \Rn$ onto a feasible set $\Q$:
\begin{align*}
	\min_{\bm x \in \Rn} 	& \norm{\bm x - \tilde{\x}^0}_2^2 \numberthis \label{eq:proj_quadric} \\ 
	\text{s.t. } 	& \bm x \in \Q ,
\end{align*}
where $\Q$ is a nonempty and non-cylindrical central quadric~\cite[Theorem 3.1.1]{odehnal_universe_2020}. In other words, $\Q$ is nonempty and there exists a quadratic function
 \[\Psi: \mathbb{R}^n \to \mathbb{R}: \x \mapsto \Psi(\x) = \Tr{\x} \B \x + \Tr{\Ab} \x + \Ac ,\]
 with $\A$ nonsingular and $c \neq \frac{\Tr{\Ab} \A^{-1} \Ab}{4}$, such that
\begin{equation}\label{eq:quadric}
	\Q = \left \{ \bm x \in \Rn \st{\Psi(\x) = 0}  \right \}= \Psi^{-1}(0).
\end{equation}
See~\cite[\S 3.1]{odehnal_universe_2020} and ~\cite[Chapter 21]{2010mathematik} for a complete classification of quadrics.

This quadratic surface, or \emph{quadric}, is denoted as quadric with middle point, in the sense of~\cite{2010mathematik}. The middle point or \emph{centre}, $\bm d$, corresponding to the centre of symmetry, is computed as $- \frac{\A^{-1} \Ab}{2}$, and the condition  $c \neq \frac{\Tr{\Ab} \A^{-1} \Ab}{4}$ is equivalent to $\bm d \notin \Q$.  

Note that, under these assumptions, we can prove that the feasible set defined by~\cref{eq:quadric} is a manifold, see~\cite[\S 3.2.1]{vh22} for more details. This centre, and the characterization of the surface as a manifold, will be used in~\cref{sub:Quasi-projection_on_the_quadric} to build a fast but inexact projection mapping, referred to as a \emph{quasi-projection}.

The problem defined by~\cref{eq:proj_quadric} is invariant to translations and rotations. Hence, without loss of generality we can consider the following problem in \emph{normal form}~\cite[Theorem 3.1.1]{odehnal_universe_2020}:

\begin{align*}
	\min_{\x \in \Rn} 	& \norm{\x - \x^0}_2^2  \numberthis \label{eq:proj_quadric_std} \\
	\text{s.t. }		& \sum_{i=1}^n \lambda_i x_i^2 -1 = 0 ,
\end{align*}
where $\bm \lambda = \spec{\A}$, contains the eigenvalues of $\A$ sorted in descending order, and $\x^0$ is the appropriate transformation of $\tilde{\x}^0$. Note that, since the feasible set is nonempty, we have $\lambda_1 >0$. We will refer to the solution(s) of~\cref{eq:proj_quadric_std} as the \emph{(true) projection(s)} of $\bm x^0$ onto the quadric.

Note that the centre $\d$ is now the origin $\bm 0$. Since this problem is symmetric with respect to the axes, we consider that $\x^0 \geq 0$, \ie, $\x^0$ is located inside the first orthant ($\mathbb{R}^{n}_+:=\left \{ \x \in \Rn \st{\x \geq0} \right \} = \mathbb{R}^{n,*}_+ \bigcup \left \{ \bm 0 \right \}$). 

Remark that if $\lambda_n > 0$, \ie, the quadric is an ellipsoid, and if we have $\sum_{i=1}^n \lambda_i (x^0_i)^2 - 1 > 0$, then the solution of~\cref{eq:proj_quadric_std} is identical to the solution of 

\begin{align*}
	\min_{\x \in \Rn} 	& \norm{\x^0 - \x}_2^2  \\
	\text{s.t. }		& \sum_{i=1}^n \lambda_i x_i^2 -1 \leq 0 ,
\end{align*}
which is a \emph{convex} optimization problem that is easy to solve, \eg, using interior points method (IPM), see~\cite{boyd_convex_2004,nesterov_lectures_2018} for more details, or a black-box commercial solver such as Gurobi~\cite{gurobi}. On the other hand, if $\A$ is indefinite or if $\sum_{i=1}^n \lambda_i (x^0_i)^2 - 1 < 0$, then we are confronted with a nonconvex optimization problem.

First, let us show that the problem is well-posed, \ie, that there exists a global optimum of~\cref{eq:proj_quadric_std}.
\begin{proposition}
	\label{prop:existence_projection_quadric}
	There exists a global optimum $\x^*$ of~\cref{eq:proj_quadric_std}.
\end{proposition}
\begin{proof}
	The objective of~\cref{eq:proj_quadric_std} is a real-valued, continuous and coercive function defined on a nonempty closed set, therefore there exists a global optimum~\cite[Theorem 2.32]{beck_2014}.
\end{proof}

\subsection{KKT conditions}%
\label{sub:KKT_conditions}

Since the feasible set is nonconvex, the projection operator does not always return a singleton, see \cite[Theorem 3.8]{fletcher_chebyshev_2015}. The set of solutions may be a singleton (\cref{fig:x_trajectories_mu}), a finite set (\cref{fig:x_trajectories_mu_degenerate_2}), or an infinite set (suppose that $\bm{x}^0$ is the centre of a sphere, \ie, $\bm \lambda = \mathbb{1}$ and $\x^0 = \bm{0}$).
Using the KKT conditions, we can characterize the solutions of~\cref{eq:proj_quadric_std}. The Lagrangian of~\cref{eq:proj_quadric_std}, with Lagrange multiplier $\mu$ and with $\D = \textrm{diag}(\bm \lambda) \in \mathbb{R}^{n\times n}$, reads

\begin{equation}
	\label{eq:lagrangian_not_sphere}
	\mathcal{L}(\x, \mu) = \Tr{(\x - \x^0)} (\x - \x^0) + \mu(\Tr{\x} \D \x - 1),
\end{equation}
and the gradient,
\begin{align*}
	\label{eq:grad_lagrangian_not_sphere_real}
	\bm \nabla \mathcal{L} (\x, \mu)	&= \begin{pmatrix} 
							2 (\x - \x^0) + 2 \mu \D \x \\ \numberthis
							 \Tr{\x} \D \x - 1 			
						\end{pmatrix}.
\end{align*}
Note that we can write the $i$-th equation of~\cref{eq:grad_lagrangian_not_sphere_real} as
\begin{equation}
	\label{eq:x_i}
	x_i (1+\mu \lambda_i) = x^0_i .
\end{equation}
Any point $(\x, \mu)$ that satisfies
\begin{equation}
	\label{eq:grad_lagrangian_not_sphere}
	\bm \nabla \mathcal{L} (\x, \mu)	= \bm 0 
\end{equation}
is referred to as a \emph{KKT point}. An optimal solution, $\x^*$, of~\cref{eq:proj_quadric_std} must either meet the KKT conditions~\cref{eq:grad_lagrangian_not_sphere}, or fail to satisfy the linear independence constraint qualification (LICQ) criterion. In~\cref{eq:proj_quadric_std}, the latter occurs if $\nabla \Psi (\x^*) = \bm 0$. This corresponds to the case where the centre belongs to the quadric, and is ruled out by the condition $c \neq \frac{\Tr{\Ab} \A^{-1} \Ab}{4}$.

Isolating $\x$ in the first $n$ equations of~\cref{eq:grad_lagrangian_not_sphere} yields, for $\mu \notin \pi(\A) := \left \{ -\frac{1}{\lambda} \st{\lambda \text{ is an eigenvalue of }\A} \right \}$,
\begin{equation}
	\label{eq:x_not_sphere}
	\x (\mu) =(\bm I + \mu \D)^{-1}  \x^0  ,
\end{equation}
and the $i-$th component can be rewritten as
\begin{equation}
	\label{eq:x_not_sphere_j}
	x_i (\mu) = \frac{x_i^0}{1+\mu \lambda_i}.
\end{equation}
Note that the set $\pi(\A)$ corresponds to the \emph{poles} of the rational function~\cref{eq:f_mu}.

We distinguish two cases:
\begin{itemize}
	\item Case 1: $\mu \notin \pi( \A )$. The matrix $(\bm I + \lambda \bm D)$ is nonsingular and we have~\cref{eq:x_not_sphere}.
	\item Case 2: $\mu = -\frac{1}{\lambda_i}$ for some $i =1, \hdots, n$. The $i$-th equation of~\cref{eq:grad_lagrangian_not_sphere} reads $-2 x_i^0 = 0$, therefore $\mu$ is a solution only if $x_i^0 = 0$. We first treat the case $\bm x^0 >0$, denoted as nondegenerate case, in~\cref{sub:Ellipsoid_case,sub:Hyperboloid_case}. Then, the degenerate case $\bm x^0 \geq 0$ is tackled in~\cref{sub:degenerate_case}.
\end{itemize}

In the first case, if we insert~\cref{eq:x_not_sphere} in the quadric equation, $\Psi(\x) = 0$, we obtain a \emph{univariate}, \emph{extended-real-valued} function
\begin{align*}
	f : \mathbb{R} \to \overline{\mathbb{R}} : \mu \mapsto f(\mu) &= \Psi(\x(\mu)) = \Tr{\x(\mu)} \D \x(\mu) -1 \\
							   &= \sum_{i=1, x^0_i \neq 0}^n \lambda_i \left(\frac{x_i^0}{1+\mu \lambda_i}\right )^2 -1  \numberthis \label{eq:f_mu}
\end{align*}
of which we want to obtain the roots. Notice that, as the roots correspond to the values $\mu^*$ for which $\Psi(\x(\mu^*)) = 0$, they can be geometrically understood as the intersections between $\left \{\x(\mu): \mu \in \mathbb{R} \right \}$ and the quadric $\mathcal{Q} = \Psi^{-1}(0)$. This is illustrated in~\cref{ssub:2D_example_of_nondegenerate_projection_onto_an_ellipse,ssub:2D_example_of_nondegenerate_projection_onto_an_hyperbola}.

In the following, we show how to efficiently solve~\cref{eq:proj_quadric_std} by computing a specific root of~\cref{eq:f_mu}. We first consider the case where $\x^0 >0$ for the ellipsoid case in~\cref{sub:Ellipsoid_case} and the hyperboloid case in~\cref{sub:Hyperboloid_case}. Then we discuss the case $\x^0 \geq \bm 0$ in~\cref{sub:degenerate_case}. Finally, we bring everything together into a single algorithm, \cref{alg:exact_projection_quadric}, in~\cref{sub:Bringing_everything_together}. We also propose in~\cref{sub:Quasi-projection_on_the_quadric} a simpler procedure that allows us to map a point to the quadric without having to diagonalize the matrix $\A$. As this mapping does not return the true projection, we refer to it as \emph{quasi-projection}.

\subsection{Ellipsoid case, \texorpdfstring{$\mathbf{x}^0 >\mathbf{0}$}{x0>0}}%
\label{sub:Ellipsoid_case}

Here, we assume that the quadric is an ellipsoid, \ie, $\bm \lambda >0$ and that the initial point lies (strictly) in the first orthant, \ie, $\x^0 > \bm{0}$.

The goal of this section is twofold. First, we derive several successive results~(\cref{prop:same_orthant_ellipsoid,prop:f_decreasing,prop:uniqueness_root,prop:uniqueness_of_solution_orthant}) that characterize the roots of $f$ and the solutions of~\cref{eq:proj_quadric_std}. The combination of these results yields~\cref{prop:optimal_is_root} which states that \emph{\cref{eq:proj_quadric_std} can be solved by finding the root of $f$ on a given interval $\mathcal{I}$}. Second, we provide a starting point for the Newton root-finding algorithm for efficiently computing this root.

\begin{proposition}\label{prop:same_orthant_ellipsoid}
	Under the standing assumptions, every solution $\x^*$ of~\cref{eq:proj_quadric_std} satisfies $\x^*>0$.
\end{proposition}
\begin{proof}
		Recall that any solution of~\cref{eq:proj_quadric_std} is a KKT point. Using~\cref{eq:x_i} we see that if $(\x^*, \mu^*)$ is a KKT point then the positivity of $x_i^0$ for all $i = 1, \hdots, n$ implies that $x^*_i \neq 0$ for all $i$.

Let us suppose, for the sake of contradiction, that $\x^*$ is a minimizer of~\cref{eq:proj_quadric_std} and that there exists a nonempty set of indices $J \subseteq \left \{ 1, 2, \hdots, n \right \}$ with $x^*_j < 0$ for all $j \in J$.
		By symmetry, we can construct $\x^{**}$ defined as
		\[
			\x^{**} :=
				\begin{cases}
					x^*_i & \text{if } i \notin J, \\
					-x^*_i & \text{if } i \in J,
				\end{cases}
		\]
		and we have $\Psi(\x^{**}) = 0$, \ie, the point belongs to the quadric. The (squared) objective can be computed:
		\begin{align*}
			\norm{\x^{**} - \x^0}_2^2 	&= \sum_{i=1}^n (x^0_i - x^{**}_i)^2 \\
							&= \sum_{i=1, i\notin J}^n (x^0_i - x^{*}_i)^2 + \sum_{j \in J} \underbrace{(x_j^0 + x_j^*)^2}_{<(x^0_j - x^*_j)^2} \\
							&< \norm{\x^{*} - \x^0}_2^2
		\end{align*}
		This contradicts the optimality of $\x^*$.
	\end{proof}

	\begin{proposition}\label{prop:f_decreasing}
		$f$, defined as in~\cref{eq:f_mu}, is strictly decreasing on $\mathcal{I}:= \left ] -\frac{1}{\lambda_1}, +\infty \right [$.
	\end{proposition}
	\begin{proof}
		Since $f \in \mathcal{C}^1$ on $\mathcal{I}$, we compute
		\begin{equation} \label{eq:df_mu}
			f'(\mu) = -2 \sum_{i=1}^n \frac{(\lambda_i  x^0_i)^2}{{\underbrace{(1+\mu \lambda_i)}_{>0\text{ for }\mu \in \mathcal{I}}}^3},
		\end{equation}
		and this function is negative on $\mathcal{I}$. 

	\end{proof}

	\begin{proposition}
		\label{prop:uniqueness_root}
		Function $f$ restricted to $\mathcal{I}$ has one and only one zero.
	\end{proposition}
	\begin{proof}
		By~\cref{prop:f_decreasing}, $f$ is strictly decreasing, and hence has \emph{at most} one zero.
		Moreover, $\lim\limits_{\mu \to +\infty}f(\mu) = \Psi(\d) = -1 < 0$ and $\lim\limits_{\mu \to -\frac{1}{\lambda_1}} f(\mu) = + \infty$; the continuity of $f$ on $\mathcal{I}$ implies the existence of the zero on $]-\frac{1}{\lambda_1}, +\infty[$.
	\end{proof}
	\begin{proposition}\label{prop:uniqueness_of_solution_orthant}
		If $\mu^*$ is a root of $f$, and $\mu^* \notin \mathcal{I}$, then $\x(\mu^*) \notin \Rn_{+}$.
	\end{proposition}
	\begin{proof}
		If $\mu^* \notin \mathcal{I}$, then $\mu^* < -\frac{1}{\lambda_1}$ and therefore the first component of $\x(\mu^*)$ reads
		\begin{equation}
			x_1(\mu^*) = \frac{x^0_1}{1+\mu^* \lambda_1} 
		\end{equation}
		As the denominator is negative, $\x(\mu^*)$ belongs to a different orthant than $\x^0$.
	\end{proof}

	\begin{proposition}\label{prop:optimal_is_root}
		If $x^0_i \neq 0$ for all $i = 1, \hdots, n$ and $\bm \lambda >0$, then the optimal solution of~\cref{eq:proj_quadric_std} is given by the unique root $\mu^*$ of $f$ restricted to $\mathcal{I}$. 
	\end{proposition}
	\begin{proof}
		As shown in~\cref{sub:KKT_conditions}, the optimal solution $\x^*$ is a KKT point, meaning that it satisfies~\cref{eq:grad_lagrangian_not_sphere}.
Using~\cref{prop:same_orthant_ellipsoid}, $\x^*$ belongs to the same orthant as $\x^0$, we are therefore interested in the best KKT solution in the first orthant. However,~\cref{prop:uniqueness_of_solution_orthant} shows that the corresponding $\mu^*$ of the KKT solutions belonging to the same orthant of $\x^0$ are located in $\mathcal{I}$ and~\cref{prop:uniqueness_root} proves the existence and uniqueness of a root on $\mathcal{I}$, which corresponds therefore to the optimal solution of~\cref{eq:proj_quadric_std}.
	\end{proof}

	\begin{proposition}\label{prop:f_convex}
		$f$ is strictly convex on $\mathcal{I}$.
	\end{proposition}
	\begin{proof}
		$f \in \mathcal{C}^2(\mathcal{I})$ and we compute
		\[ f''(\mu) = 6 \sum_{i=1}^n \frac{(\lambda_i x^0_i)^2 \lambda_i}{(1+\mu \lambda_i)^4},  \]
		which is positive on $\mathcal{I}$.
		
	\end{proof}

	\begin{proposition}\label{prop:convergence_newton}
		Let $\mu^0 \in \mathcal{I} =] -\frac{1}{\lambda_1}, +\infty[$ with $f(\mu^0)>0$. The Newton-Raphson algorithm with starting point $\mu^0$ converges to $\mu^*$, the unique root of $f$ on $\mathcal{I}$ (as in~\cref{prop:optimal_is_root}).
	\end{proposition}
	\begin{proof}
		Let us now prove by induction on $k$ that the sequence  $(\mu^k)_{k \in \mathbb{N}}$  provided by Newton's method is an increasing sequence upper bounded by $\mu^*$.
		The Newton-Raphson iterate for $k = 0,1,\hdots$ is given by
		\begin{equation}\label{eq:newton_iterate}
			\mu^{k+1} = \mu^{k} - \frac{f(\mu^k)}{f'(\mu^k)} \, .
		\end{equation}
		Using the induction hypothesis, which implies that $f(\mu^k) > 0$ for $\mu^k \neq \mu^*$, and~\cref{prop:f_decreasing}, we have \[ \mu^k < \mu^{k+1}.\] 
		Since $f$ is strictly convex on $\mathcal{I}$ (\cref{prop:f_convex}), the tangent of $f$ at a given point is below any chord starting from this point. In particular we have
		\[ f'(\mu^k) < \frac{f(\mu^k) - f(\mu^*)}{\mu^k - \mu^*},  \]
		Using the definition of $\mu^*$ and rearranging, we obtain
		\begin{align*}
			\mu^k - \frac{f(\mu^k)}{f'(\mu^k)}  	&< \mu^*, \\
			\mu^{k+1}				&< \mu^*.
		\end{align*}
		Since the sequence $(\mu^k)_{k\in\mathbb{N}}$ is strictly increasing (for $\mu^k \neq \mu^*$) and bounded, it must converge to a fixed point of~\cref{eq:newton_iterate} which corresponds to a root of $f$. This concludes the proof as, by~\cref{prop:optimal_is_root}, there is a unique root of $f$ on $\mathcal{I}$ corresponding to the optimal solution of~\cref{eq:proj_quadric_std}.
	\end{proof}

	\subsubsection{2D example of a nondegenerate projection onto an ellipse}%
	\label{ssub:2D_example_of_nondegenerate_projection_onto_an_ellipse}

\begin{figure}
	\begin{center}
		\includegraphics[width=\textwidth]{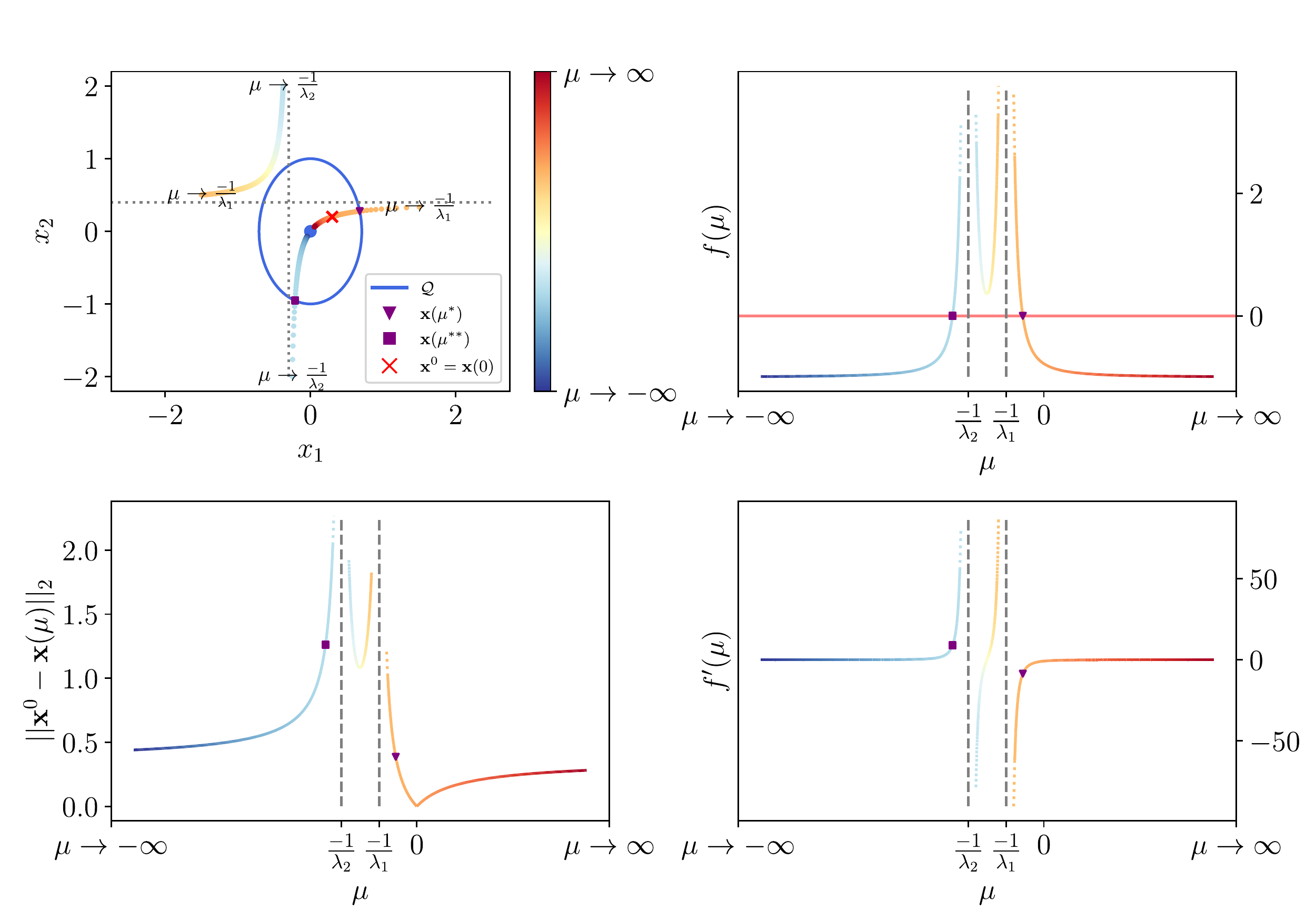}
	\end{center}
	\caption{Image of $\x(\mu)$ and graphs of $f(\mu)$, $f^\prime(\mu)$, $\norm{\x(\mu) - \x^0}_2$ for $\mu$ ranging on $]-\infty, +\infty[$ in the nondegenerate elliptic case, here $\Q (\Psi) := \Q$.}
	\label{fig:x_trajectories_mu}
	\end{figure}

	Figure~\ref{fig:x_trajectories_mu} presents an example of a nondegenerate projection, that is with $\bm x^0 >0$, onto an ellipse. We plot $\x(\mu)$, $f(\mu)$, $f'(\mu)$ and $\norm{\x(\mu)- \x^0}_2$ for $\mu$ ranging on $]-\infty, \infty[$. Let us describe how $\x(\mu)$, in the top left subfigure, varies when $\mu$ decreases from $+\infty$ to $-\infty$. For $\mu \to +\infty$, we have $\x(\mu)= \d$, where $\d = \bm{0}$ is the quadric centre depicted as a blue dot. Then, while decreasing $\mu$ to $0$, we reach $\x(0) = \x^0$. For $\mu \to -\frac{1}{\lambda_1}$, $\x(\mu)$ follows an asymptote \emph{and crosses the quadric on $\x(\mu^*)$}, the optimal solution of~\cref{eq:proj_quadric_std}, depicted as a purple triangle. Further decreasing $\mu$, $\x(\mu)$ reappears on the left part of the asymptote ($x_1 \to -\infty$) and tends to the asymptote ($x_2 \to +\infty$) defined by the other eigenvalue. Finally, $\x(\mu)$ converges to the quadric centre when $\mu \to -\infty$, passing again through the quadric in $\x(\mu^{**})$, the max point of~\cref{eq:proj_quadric_std}, depicted as a purple square.

	The function $f$ is also depicted with its two roots. Note that, depending on the parameters of the problem, it may have one or two additional roots corresponding to local minima or maxima. 
	We also observe in the bottom right figure that $f'(\mu)$ is negative on $\mathcal{I}$, and show the distance to $\x^0$ for the different values of $\x(\mu)$ in the bottom left figure.

	\subsection{Hyperboloid case, \texorpdfstring{$\mathbf{x}^0 > \mathbf{0}$}{x0 > 0}}%
	\label{sub:Hyperboloid_case}

	In the hyperboloid case, there is at least one positive and one negative eigenvalue of $\A$. Let $1 \leq p \leq n-1$ be the number of positive eigenvalues.
	Let us consider $e_1:= -\frac{1}{\lambda_1}$ and $e_2 := -\frac{1}{\lambda_{n}}$. We have $0 \in ]e_1, e_2[$. We will work analogously as in the ellipsoidal case, with $\mathcal{I} := ]e_1, e_2[$.
	\cref{prop:same_orthant_ellipsoid} has no assumption with respect to the positivity of $\lambda$, and therefore remains valid. The other propositions can be successively adapted: \cref{prop:f_decreasing_hyperboloid} adapts~\cref{prop:f_decreasing}, \cref{prop:uniqueness_root_hyperboloid} adapts~\cref{prop:uniqueness_root}, \cref{prop:uniqueness_of_solution_orthant_hyperboloid} adapts~\cref{prop:uniqueness_of_solution_orthant}, and finally the main result remains valid, \ie, \cref{prop:optimal_is_root_hyperboloid} adapts~\cref{prop:optimal_is_root}.

	We then propose \cref{alg:double_newton}, also based on Newton-Raphson, to efficiently compute the root of $f$ in $\mathcal{I}$, and hence one of the optimal solutions of~\cref{eq:proj_quadric_std}.
	An example of the hyperbolic case is provided in~\cref{fig:x_trajectories_mu_hyperbola}.

	\begin{proposition}\label{prop:f_decreasing_hyperboloid}
		$f$, defined as in~\cref{eq:f_mu}, is strictly decreasing on $\mathcal{I}:= \left ] e_1, e_2 \right [$.
	\end{proposition}
	\begin{proof}
		Since $f \in \mathcal{C}^1(\mathcal{I})$, we compute
		\begin{equation} \label{eq:df_mu_hyperboloid}
			f'(\mu) = -2 \sum_{i=1}^p \frac{(\lambda_i  x^0_i)^2}{{\underbrace{(1+\mu \lambda_i)}_{>0 \text{ if }\mu > e_1}}^3} -2 \sum_{i=p+1}^n \frac{(\lambda_i  x^0_i)^2}{{\underbrace{(1+\mu \lambda_i)}_{>0 \text{ if }\mu < e_2}}^3} ,
		\end{equation}
		and this function is negative on $\mathcal{I}$. 

	\end{proof}
	\begin{proposition}\label{prop:uniqueness_root_hyperboloid}
		Function $f$ restricted to $\mathcal{I}$ with $\x^0 >0$ has one and only one zero.
	\end{proposition}
	\begin{proof}
		By~\cref{prop:f_decreasing_hyperboloid}, $f$ is strictly decreasing, and hence has \emph{at most} one zero.
		Moreover, $\lim\limits_{\mu \to e_1}f(\mu) = +\infty$ and $\lim\limits_{\mu \to e_2} f(\mu) = - \infty$; the continuity of $f$ on $\mathcal{I}$ implies the existence of the zero on $\mathcal{I}$.
	\end{proof}
	\begin{proposition}\label{prop:uniqueness_of_solution_orthant_hyperboloid}
		If $\mu^*$ is a root of $f$, and $\mu^* \notin \mathcal{I}$, then $\x(\mu^*) \notin \Rn_+$.
	\end{proposition}
	\begin{proof}
		If $\mu^* \notin \mathcal{I}$, then either $\mu^* < -\frac{1}{\lambda_1}$ or $\mu^* > -\frac{1}{\lambda_{n}}$.
		The first case is already treated in the proof of~\cref{prop:uniqueness_of_solution_orthant}. For the second case, we note that
		\begin{equation}
			x_{n}(\mu^*) = \frac{x^0_{n}}{1+\mu^* \lambda_{n}}. 
		\end{equation}
		As $\lambda_{n}<0$, the denominator is negative, and thus $\x(\mu^*)$ belongs to a different orthant than $\x^0$.
	\end{proof}

	\begin{proposition}\label{prop:optimal_is_root_hyperboloid}
		If $x^0_i \neq 0$ for all $i = 1, \hdots, n$, then the optimal solution of~\cref{eq:proj_quadric_std} is given by the unique root $\mu^*$ of $f$ restricted to $\mathcal{I}:=]e_1, e_2[$. 
	\end{proposition}
	\begin{proof}
		Since~\cref{prop:same_orthant_ellipsoid,prop:uniqueness_of_solution_orthant,prop:uniqueness_root} are also valid in the hyperboloid case with $\mathcal{I}:=]e_1, e_2[$, the proof is identical to~\cref{prop:optimal_is_root}.
	\end{proof}

	\begin{proposition}\label{prop:uniqueness_inflexion}
		There exists a unique inflexion point, $\mu^\mathrm{I}$, of $f$ on $\mathcal{I}$.
	\end{proposition}
	\begin{proof}
		This follows from the monotonicity of $f^{\prime\prime} \in \mathcal{C}^{\infty}(\mathcal{I})$, \ie, $f^{\prime\prime\prime}(\mu) < 0$ for all $\mu \in \mathcal{I}$, and because $\lim\limits_{\mu \to e_1^+} f^{\prime\prime}(\mu) = - \lim\limits_{\mu \to e_2^-} f^{\prime\prime}(\mu) = +\infty $.
	\end{proof}

	Since there is a single inflexion point $\muI$, we can launch in parallel two Newton's algorithms and guarantee that at least one will converge.

	\begin{algorithm}
		\caption{Double Newton}
		\label{alg:double_newton}
		\begin{algorithmic}
			\IF{$f(0) <0$} 
			\STATE Use bisection method (see~\cite[Chapter 2.1]{Burden01}) to find $\mu_\mathrm{s} \in ]e_1, 0[$ s.t. $f(\mu_\mathrm{s}) >0$
			\ELSIF{$f(0) > 0$}
			\STATE Use bisection method to find $\mu_\mathrm{s} \in ]0, e_2[$ s.t. $f(\mu_\mathrm{s}) <0$
			\ELSE
			\RETURN 0
			\ENDIF
			\STATE $\mu_0 \gets \textrm{Newton}(0)$ \COMMENT{This and the next line are run in parallel}
			\STATE $\mu_1 \gets \textrm{Newton}(\mu_{\textrm{s}})$
			\RETURN $\mu_0, \mu_1$ \COMMENT{Returns the output solution of the first of the two parallel Newtons that is finished}
		\end{algorithmic}
	\end{algorithm}

	\begin{proposition}
		One of the two Newton's methods of~\cref{alg:double_newton} converges to $\mu^*$, with $\x(\mu^*)$ (defined in~\cref{eq:x_not_sphere}) the optimal solution of~\cref{eq:proj_quadric_std}.
	\end{proposition}
	\begin{proof}
		This proof relies on the double initiation of Newton's method in~\cref{alg:double_newton}: one starting from a positive value, and the other from a negative value.
		We comment on the sign of $f(\mu^\mathrm{I})$:
		\begin{itemize}
			\item If $f(\muI) < 0$, then $\mu^* \in ]e_1, \muI[$ and the function is convex on this interval. The situation is similar to~\cref{prop:convergence_newton}, and any starting point $\mu_\mathrm{s}$ with $f(\mu_\mathrm{s}) >0$ is a valid starting point, in the sense that the sequence of iterates converges to $\mu^*$.
			\item If $f(\muI) > 0$, then $\mu^* \in ]\muI, e_2[$ and the function is concave on this interval. Using a similar argument as~\cref{prop:convergence_newton}, any starting point $\mu_\mathrm{s}$ with $f(\mu_\mathrm{s}) <0$ is a valid starting point.
			\item If $f(\muI) = 0$, any starting point in $\mathcal{I}$ is a valid starting point.
		\end{itemize}
		
	\end{proof}
	Remark that with the knowledge of the value of $\muI$, we could launch a single Newton scheme with the appropriate starting point. Unfortunately, computing $\muI$ amounts to computing the root of $f''$ which is at least as costly as finding the root of $f$.

	\subsubsection{2D example of a nondegenerate projection onto a hyperbola}%
\begin{figure}
	\begin{center}
		\includegraphics[width=\textwidth]{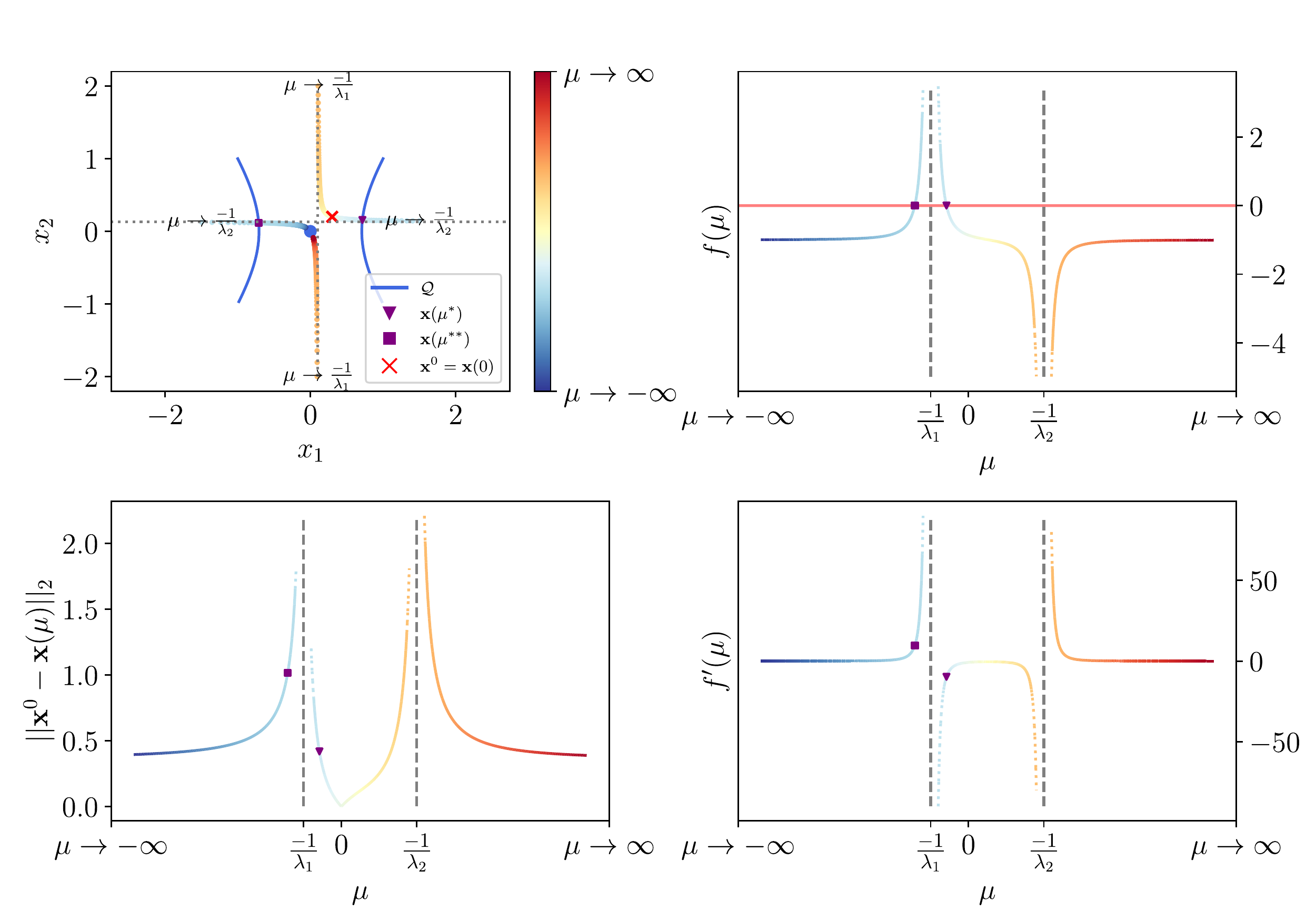}
	\end{center}
	\caption{Image of $\x(\mu)$ and graphs of $f(\mu)$, $f^\prime(\mu)$, $\norm{\x(\mu) - \x^0}_2$ for $\mu$ ranging on $]-\infty, +\infty[$ in the nondegenerate hyperbolic case, here $\Q (\Psi) := \Q$.}
	\label{fig:x_trajectories_mu_hyperbola}
	\end{figure}

	\label{ssub:2D_example_of_nondegenerate_projection_onto_an_hyperbola}
	Figure~\ref{fig:x_trajectories_mu_hyperbola} shows an example of a nondegenerate projection onto a hyperbola. We observe a similar image of $\x(\mu)$, with two asymptotes. We see that $f(\mu)$ has a unique inflexion point on $\mathcal{I}$. In this example, the inflexion point is on the right of the root, and thus we know that starting a Newton-Raphson scheme in some $\mu_{\mathrm{s}}$ with $f(\mu_\mathrm{s}) >0$ yields a sequence that converges to $\mu^*$.

	\subsection{Degenerate case, \texorpdfstring{$\mathbf{x}^0 \geq \mathbf{0}$}{x0 >= 0}}%
	\label{sub:degenerate_case}

Let us first assume that all eigenvalues of $\A$ are distinct, the case with repeated eigenvalues is treated at the end of the current section.

\subsubsection{All eigenvalues are distinct}%
\label{ssub:All_eigenvalues_are_distinct}

	The discussion is largely similar to~\cref{sub:Hyperboloid_case}, with the following differences: i) $f$ is continuous at $-\frac{1}{\lambda_i}$ if the associated component of $\x^0$ is equal to zero and ii) at most two additional KKT points can be obtained for each component of $\x^0$ that is equal to zero.

	We tackle these issues as follows. First, we arbitrarily decide to single out a solution in the first orthant. Second, we change the definition of $e_1$, $e_2$, to account for the continuity of $f$ in $-\frac{1}{\lambda_i}$: $\lim\limits_{\mu \to -\frac{1}{\lambda_i}}f(\mu) \neq \infty$ if $x^0_i = 0$. Finally, we show how to analytically compute these additional solutions.

	Let $I := \left \{1, 2, \hdots, n \right \}$ and $K \subseteq I$, it is clear that~\cref{prop:same_orthant_ellipsoid} is not valid any more if $\size{K}$ entries of $\x^0$ are equal to zero. Indeed, if $\x^*$ is an optimal solution that belongs to the same orthant as $\x^0$, then $\x^{**}$ defined as follows
	\begin{equation}
		x^{**}_i = 
		\begin{cases}
			x^*_i & \forall i\notin K,\\
			-x^*_i & \forall i\in K,\\
		\end{cases}
	\end{equation}
	is also an optimal solution. In fact, up to ${2^{\size{K}}} -1$ solutions outside the first orthant can be obtained by mirroring $\x^*$ along a selected set of components in $K$.

As we are interested in finding \emph{one} of the optimal solutions, we note that we can also restrict our search to the first orthant. 
	\begin{proposition}\label{prop:same_orthant_degenerate}
		Given $\x^0 \geq \bm{0}$,  there exists an optimal solution $\x^*$ of~\cref{eq:proj_quadric_std} such that $\x^* \geq \bm{0}$.
	\end{proposition}
	\begin{proof}
		Let $\x^*$ be an optimal solution. The existence of $\x^*$ follows from~\cref{prop:existence_projection_quadric}. Using a similar argument as in the proof of~\cref{prop:same_orthant_ellipsoid}, we have $\mathrm{sign}(x^0_j) = \mathrm{sign}(x^*_j) \quad \forall j \in I \setminus K$. Let
		\[x^{**}_i = 
		\begin{cases}
			-x^*_i & \forall i\in K \text{ with } \textrm{sign}(x^*_i) \neq \textrm{sign}(x^0_i),\\
			x^*_i 	& \text{elsewhere.}
		\end{cases}
	\]
		This feasible point has the same objective as $\x^*$ and is located in the same orthant as $\x^0$.
	\end{proof}

	For $\mu \notin -\frac{1}{\spec{\A}}$, we change the definition of $e_1, e_2$ as
	\begin{align*}\label{eq:e_1_e_2_degenerate}
		e_1 &= \max_{\{i\in I \st{\lambda_i>0, x^0_i \neq 0}\}} -\frac{1}{\lambda_i} \numberthis  \\
		e_2 &= \min_{\{i\in I \st{\lambda_i < 0, x^0_i \neq 0}\}} -\frac{1}{\lambda_i}  \\
	\end{align*}
	and $e_1 = - e_2 := -\infty$ if the $\max$ or $\min$ is empty. This takes into account the continuity of $f$ at $\mu = -\frac{1}{\lambda_i}$ if $x_i^0 = 0$ for some index $i$.

	Let us adapt~\cref{prop:f_decreasing_hyperboloid,prop:uniqueness_root_hyperboloid,prop:uniqueness_of_solution_orthant_hyperboloid} to the degenerate case.

	\begin{proposition}\label{prop:f_decreasing_hyperboloid_degenerate}
		$f$, defined as in~\cref{eq:f_mu}, is strictly decreasing on $\mathcal{I}:= \left ] e_1, e_2 \right [$.
	\end{proposition}
	\begin{proof}
		Since $f \in \mathcal{C}^1(\mathcal{I})$, we compute
		\begin{equation} \label{eq:df_mu_hyperboloid_degenerate}
			f'(\mu) = -2 \sum_{i=1,i \notin K}^p \frac{(\lambda_i  x^0_i)^2}{{\underbrace{(1+\mu \lambda_i)}_{>0 \text{ if }\mu > e_1}}^3} -2 \sum_{i=p+1, i\notin K}^n \frac{(\lambda_i  x^0_i)^2}{{\underbrace{(1+\mu \lambda_i)}_{>0 \text{ if }\mu < e_2}}^3} ,
		\end{equation}
		and this function is negative on $\mathcal{I}$. 
	\end{proof}
	\begin{proposition}\label{prop:uniqueness_root_hyperboloid_degenerate}
		Function $f$ restricted to $\mathcal{I}$ has one and only one zero if $\exists i \in \Ip:= \left \{ i\in I \st{\lambda_i >0} \right \}$ with $x^0_i \neq 0$.
	\end{proposition}
	\begin{proof}
		We note that the technical assumption on $x^0_i$ ensures that $\lim\limits_{\mu \to e_1} f(\mu) >0$. For $\lim\limits_{\mu \to e_2} f(\mu)$ we distinguish two cases:
	\begin{itemize}
		\item either $e_2 = +\infty$ and $\lim\limits_{\mu \to e_2} = -1$;
		\item or $e_2 = \min\limits_{i \in I \st{\lambda_i <0, x^0_i \neq 0}} -\frac{1}{\lambda_i}$ and $\lim\limits_{\mu \to e_2} f(\mu) = -\infty$.
	\end{itemize}
	Since in both cases the limit is negative, and $f$ is continuous and strictly decreasing on $\mathcal{I}$, there exists a unique zero on this interval.
	\end{proof}
	
	\begin{proposition}\label{prop:uniqueness_of_solution_orthant_hyperboloid_degenerate}
		If $\mu^*$ is a root of $f$, and $\mu^* \notin \mathcal{I}$, then $\x(\mu^*) \notin \Rn_+$.
	\end{proposition}
	\begin{proof}
		If $\mu^* \notin \mathcal{I}$, then either $e_1 \neq -\infty$ and $\mu < e_1$ or $e_2 \neq +\infty$ and $\mu > e_2$.
		
		The proof follows from the definition of $e_i$, \eg, in the first case we note that 

		\[
			x_{i_1}(\mu^*) = \frac{x^0_{i_1}}{1+\mu \lambda_{i_1}},
		\]
		where $i_1 =   \argmax_{\{i\in I \st{\lambda_i>0, x^0_i \neq 0}\}} -\frac{1}{\lambda_i} $. This implies that $\x(\mu^*)$ belongs to a different orthant that $\x^0$ since the numerator is nonzero and the denominator is negative.
	\end{proof}

	Remark that if $x^0_i = 0 \quad \forall i \in \Ip$, meaning that the assumption on $\x^0$ of~\cref{prop:uniqueness_root_hyperboloid_degenerate} does not hold, then $f(\mu)$ reads 
	\[ f(\mu) = \sum_{i \in \In} \lambda_i \left( \frac{x^0_i}{1+\mu \lambda_i} \right)^2 -1, \]
	where $\In := \I \setminus \Ip = \left \{p+1, p+2, \hdots, n \right \}$.
	This function is negative on $\mathbb{R}$. In this specific case, \emph{$f$ does not provide any KKT point}, such a situation is depicted in~\cref{fig:x_trajectories_mu_hyperbola_degenerate_no_root}. However, the problem is solvable, due to additional KKT points that appear when $\x^0$ is located on the axes.

	Indeed, if $\mu = -\frac{1}{\lambda_k}$ for $k \in K$ then the $k-$th entry of~\cref{eq:grad_lagrangian_not_sphere} reads

		\begin{equation}
			2(x_k - x_k^0) + 2 \mu \lambda_k x_k = 0 \label{eq:lagrangian_degenerate}
		\end{equation}
which is true no matter $x_k$. Therefore, we obtain at most \emph{two additional} solutions of the Lagrangian system~\cref{eq:grad_lagrangian_not_sphere}. Geometrically, this corresponds to looking at the intersection between i) a line perpendicular to the axis corresponding to the component $k$ where $x^0_k = 0$ and ii) the quadric. These solutions, if they exist, can be computed as 
\begin{equation}
	\label{eq:x_d_degenerate}
	(\bm x^{\mathrm{d}}_k)_i =
	\begin{cases}
		\frac{x^0_i}{1- \frac{\lambda_i}{\lambda_k}} & \text{if } i \neq k ,\\
		\pm  \sqrt{\frac{1}{\lambda_k} \left (  1-\sum\limits_{j\in I,j \neq k} \lambda_j  \left( \frac{x^0_j}{1-\frac{\lambda_j}{\lambda_k}} \right )^2 \right ) }			      & \text{if } i=k ,
	\end{cases} 
\end{equation}
	where $(\cdot)_i$ selects the $i-$th component, and we choose the ``$+$'' solution that lies in the first orthant. 

	Such a situation is depicted in~\cref{fig:x_trajectories_mu_degenerate}. We observe that $\x(\mu)$ moves around the axis corresponding to the component of $\x^0$ which is equal to zero. Moreover, the additional solution is depicted in green in~\cref{fig:x_trajectories_mu_degenerate,fig:x_trajectories_mu_degenerate_2}. Note that in~\cref{fig:x_trajectories_mu_degenerate}, the optimal solution is a root of $f$ and in~\cref{fig:x_trajectories_mu_degenerate_2}, it is the additional solution.

	Remark that there is no intersection, and therefore no additional solution to~\cref{eq:grad_lagrangian_not_sphere}, if $\frac{1}{\lambda_k} \left ( 1-\sum\limits_{j \in I, j \neq k} \lambda_j \left( \frac{x^0_j}{1-\frac{\lambda_j}{\lambda_k}}\right)^2 \right ) < 0$, see, \eg,~\cref{fig:x_trajectories_mu_hyperbola_degenerate_no_green}.
	
	\subsubsection{2D examples of degenerate projections}%
\begin{figure}
	\begin{center}
		\includegraphics[width=\textwidth]{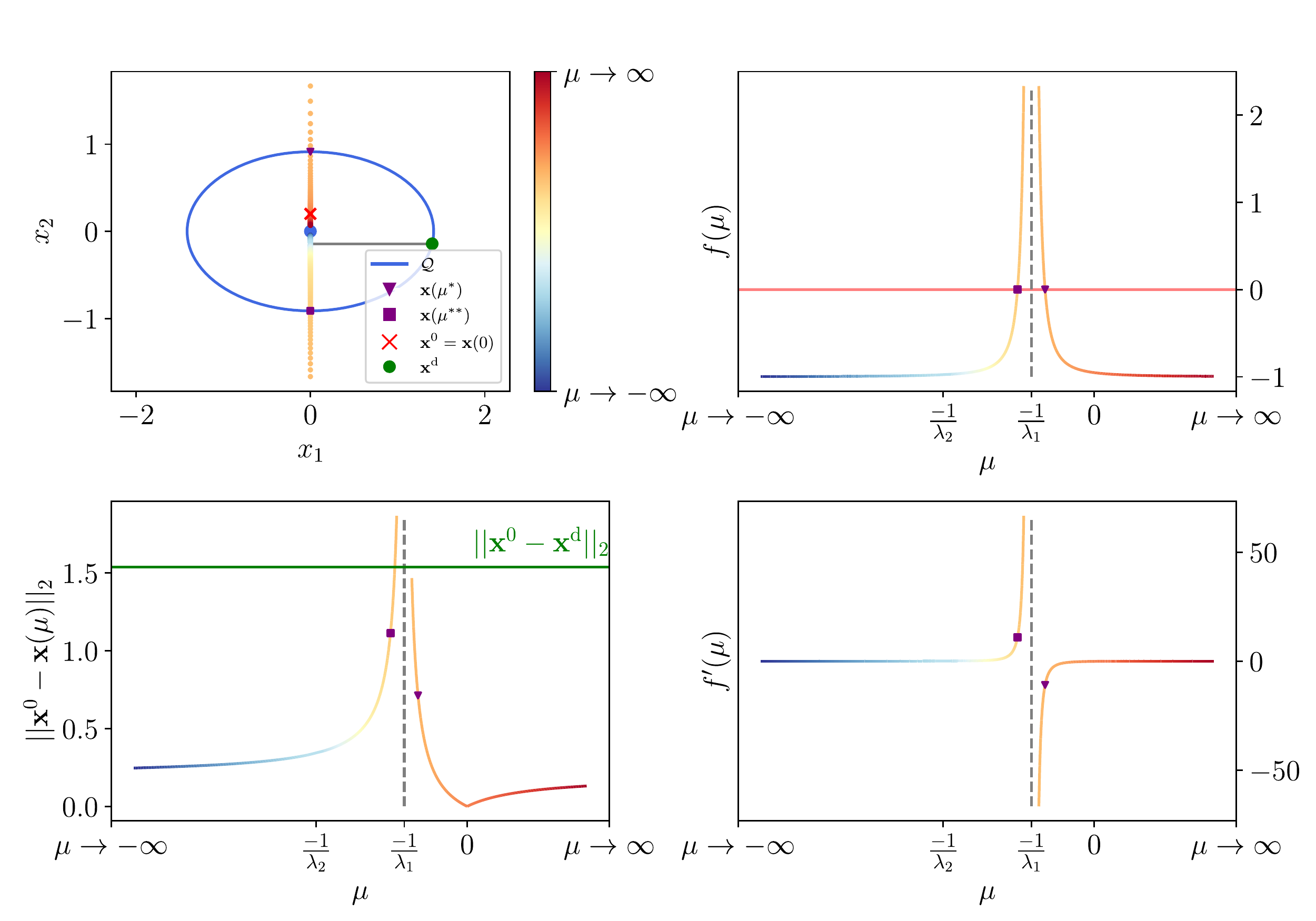}
	\end{center}
	\caption{Image of $\x(\mu)$ and graphs of $f(\mu)$, $f^\prime(\mu)$, $\norm{\x(\mu) - \x^0}_2$ for $\mu$ ranging on $]-\infty, +\infty[$ in the degenerate elliptic case. The optimal solution is the root of $f$ and not $\x^\mathrm{d}_2$: the green line showing $\norm{\x^0 - \x^\mathrm{d}_2}$ is above the purple triangle in the lower left figure.}
	\label{fig:x_trajectories_mu_degenerate}
	\end{figure}

\begin{figure}
	\begin{center}
		\includegraphics[width=\textwidth]{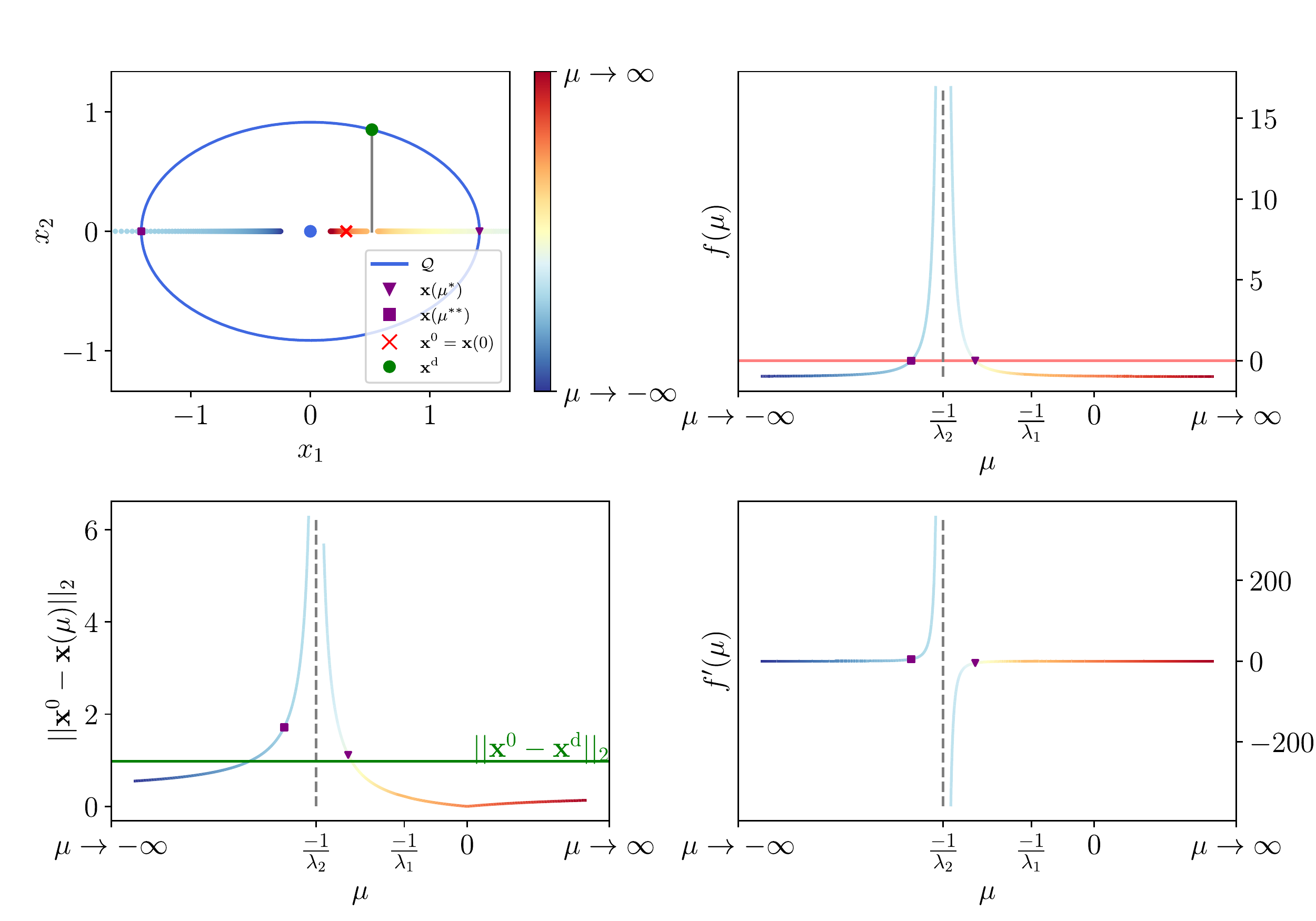}
	\end{center}
	\caption{Image of $\x(\mu)$ and graphs of $f(\mu)$, $f^\prime(\mu)$, $\norm{\x(\mu) - \x^0}_2$ for $\mu$ ranging on $]-\infty, +\infty[$ in the degenerate elliptic case. The optimal solution is not the root of $f$ but $\x^\mathrm{d}_1$: the green line showing $\norm{\x^0 - \x^\mathrm{d}_1}_2$ is below the purple triangle.}
	\label{fig:x_trajectories_mu_degenerate_2}
	\end{figure}

	\label{ssub:2D_example_of_degenerate_projections}
	\cref{fig:x_trajectories_mu_degenerate,fig:x_trajectories_mu_degenerate_2} show two examples of \emph{degenerate} projections onto an ellipse. \cref{fig:x_trajectories_mu_degenerate} depicts an example where the optimal solution is given by the KKT point corresponding to the root of $f$. \cref{fig:x_trajectories_mu_degenerate_2} depicts an example where the optimal solution is given by the KKT point corresponding to $\mu = -\frac{1}{\lambda_2}$. Notice that in these (degenerate) cases, one of the asymptotes of the $\x(\mu)$ image disappears, and the image is hence along one of the axes. Moreover, one of the discontinuities of $f$, $f^\prime$ and $\norm{\x(\mu) - \x^0}_2$ disappears as, \eg, $\lim\limits_{\mu \to -\frac{1}{\lambda_k}} f(\mu) = \sum\limits_{i=1, i\neq k}^n \lambda_i \left ( \frac{x^0_i}{1 - \frac{\lambda_i}{\lambda_k}} \right )^2 -1  \neq \infty$.

\begin{figure}
	\begin{center}
		\includegraphics[width=\textwidth]{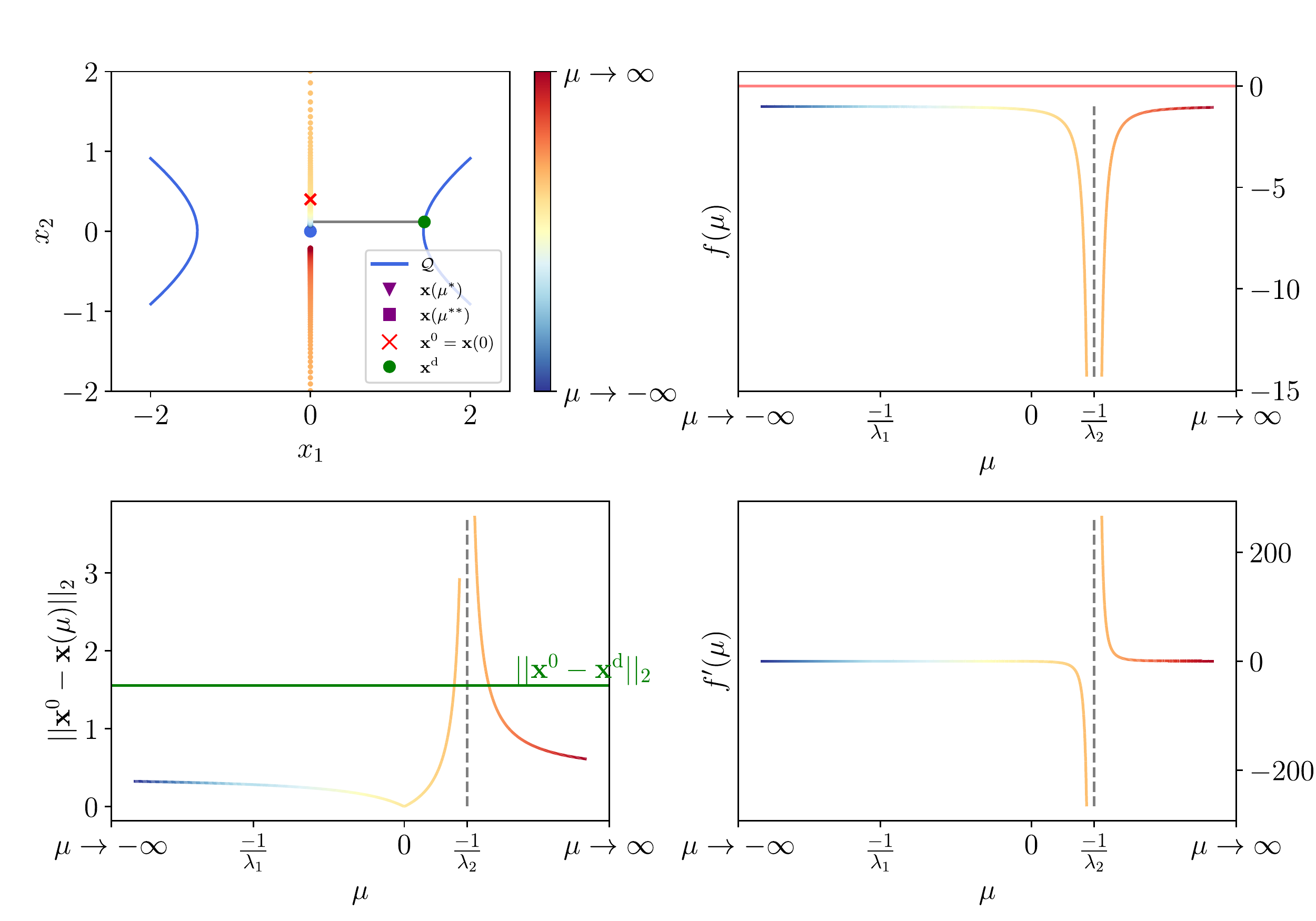}
	\end{center}
	\caption{Image of $\x(\mu)$ and graphs of $f(\mu)$, $f^\prime(\mu)$, $\norm{\x(\mu) - \x^0}_2$ for $\mu$ ranging on $]-\infty, +\infty[$ in the degenerate hyperbolic case. The optimal solution is $\x^\mathrm{d}_1$, as $f$ has no root.}
	\label{fig:x_trajectories_mu_hyperbola_degenerate_no_root}
	\end{figure}

\begin{figure}
	\begin{center}
		\includegraphics[width=\textwidth]{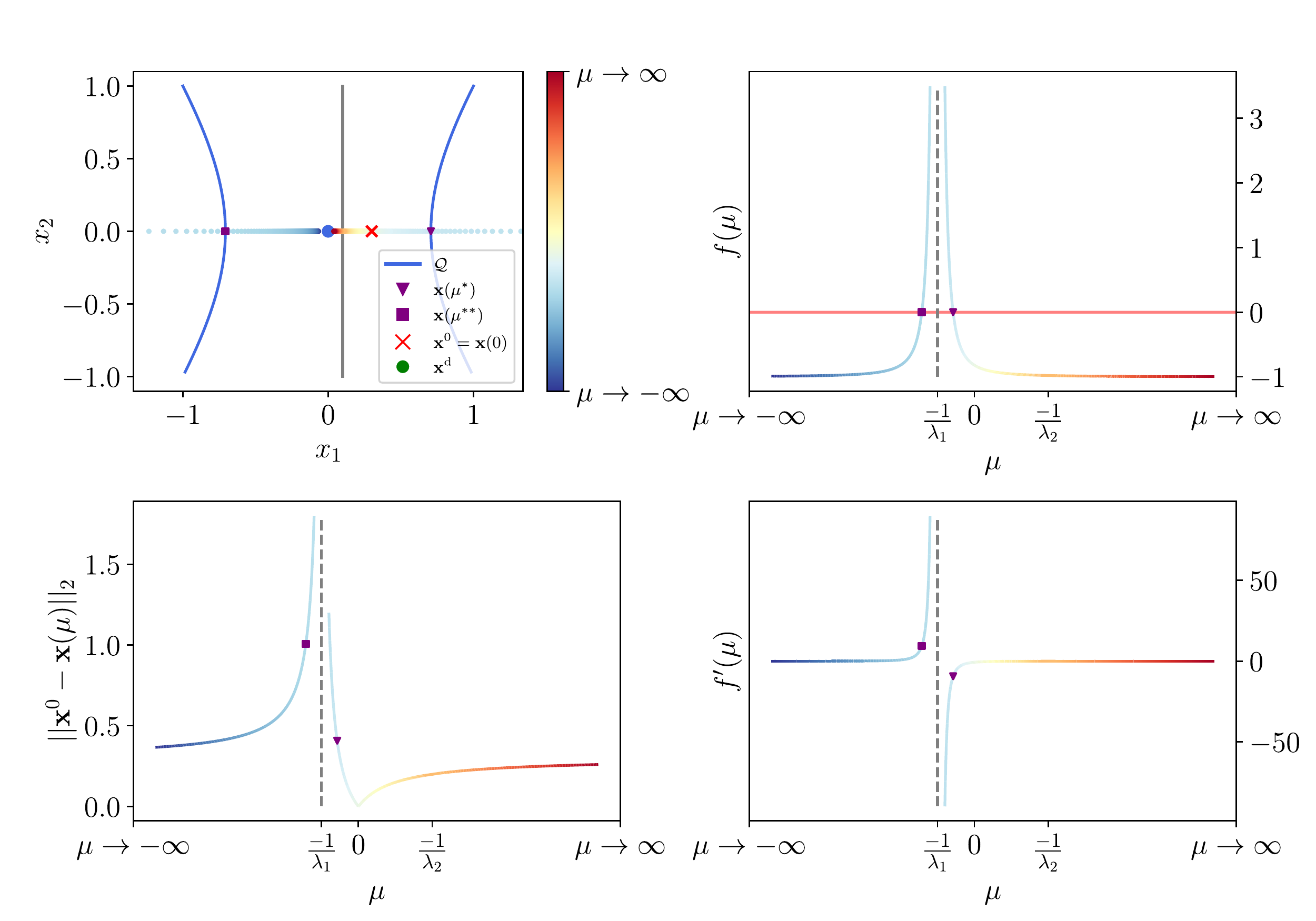}
	\end{center}
	\caption{Image of $\x(\mu)$ and graphs of $f(\mu)$, $f^\prime(\mu)$, $\norm{\x(\mu) - \x^0}_2$ for $\mu$ ranging on $]-\infty, +\infty[$ in the degenerate hyperbolic case. The optimal solution is $\x(\mu^*)$, the root of $f$. There is no additional KKT point $\x^\mathrm{d}_k$; the grey line in the upper left panel has no intersection with the quadric.}
	\label{fig:x_trajectories_mu_hyperbola_degenerate_no_green}
	\end{figure}

	\cref{fig:x_trajectories_mu_hyperbola_degenerate_no_root,fig:x_trajectories_mu_hyperbola_degenerate_no_green} show two examples of \emph{degenerate} projections onto a hyperbola. \cref{fig:x_trajectories_mu_hyperbola_degenerate_no_root} depicts an example where $f$ has no root. This is not an issue because the optimal solution is given, in this case, as one of the $\x^\mathrm{d}_k$ depicted in green which are derived in~\cref{eq:x_d_degenerate}. \cref{fig:x_trajectories_mu_hyperbola_degenerate_no_green} shows an example where there is no intersection between the grey line and the quadric, and therefore no $\x^\mathrm{d}_k$. This is not an issue because then there must exist a root $\mu^*$ of $f$ on $\mathcal{I}$, which is the optimal solution (purple triangle). Remark that, if $\sum_{i=1}^n \lambda_i (x^0_i)^2-1 > 0$, then both $\x(\mu^*)$ and $\x^\mathrm{d}_k$ are KKT points, and one of them is the optimal solution.


	\subsubsection{Some eigenvalues are repeated}%
	\label{ssub:Some_eigenvalues_are_repeated}

	Let $\overline{\bm \lambda}$ be the vector of the unique eigenvalues of $\A$, sorted in descending order, let $k \in \left \{ 1,\hdots,\size{\overline{\bm \lambda}} \right \}$ be a given component of $\overline{\bm{\lambda}}$, let $L_k$ be a subset of $I$ corresponding to the same eigenvalue, \ie, $L_k := \left \{l \in I \st{\lambda_l=\overline{\lambda}_k} \right \}$, and let $K_{k}$ be a subset of $L_k$ where the associated component of $\x^0$ is equal to zero, \ie,
	\[ 
		K_{k} := \left \{ i \in I \st{\lambda_i = \overline{\lambda}_k, x^0_i = 0} \right  \}.
	\]
	\begin{proposition}\label{prop:L_K}
		Let $L_k$ and $K_{k}$ be defined as above. There exists a solution of~\cref{eq:grad_lagrangian_not_sphere} with $\mu^* = -\frac{1}{\overline{\lambda_k}}$ only if $L_k = K_k$.
	\end{proposition}
	\begin{proof}
		Let $(\x^*, \mu^*)$ be a solution of~\cref{eq:grad_lagrangian_not_sphere} with $\mu^*=-\frac{1}{{\overline{\lambda_k}}}$. Let us assume, for the sake of contradiction, that $L_k \neq K_{k}$, or equivalently, that $\exists i \in I$ with $\lambda_i = \overline{\lambda}_k$ but $x^0_i \neq 0$. The $i-$th component of~\cref{eq:grad_lagrangian_not_sphere} reads
		\[2(x_i - x_i^0) - 2 x_i = 0 , \]
		which does not hold.
	\end{proof}
	Remark that~\cref{prop:L_K} is a left implication and it is possible that no solution of~\cref{eq:grad_lagrangian_not_sphere} exists with $\mu^* = -\frac{1}{\overline{\lambda}_k}$ and $L_k = K_k$, see, \eg,~\cref{fig:x_trajectories_mu_hyperbola_degenerate_no_root}.

	If $\size{K_{k}}=1$, the discussion is analogous to the previous paragraph: at most two KKT solutions are obtained as the intersection between a line and the quadric, but for $\size{K_{k}}>1$, we have to take the intersection between a plane $\pi := \left \{ \bm x \in \Rn \st{x_{k'} = 0 \, \forall k' \in K_{k}} \right \}$, and the quadric. Geometrically, the intersection---if there is one, \ie, if the argument of the square root below is positive---will be a $\size{K_{k}}-1$ hypersphere in the corresponding subspace of $\Rn$:

	\begin{align*}
		\pi \bigcap \Q = & \left \{ \vphantom{\frac{1}{\overline{\lambda}_k} \left ( 1- \sum\limits_{j\in I \setminus K_{k}} \lambda_j \left( \frac{x^0_j}{1-\frac{\lambda_j}{\overline{\lambda}_k}} \right)^2 \right )}  \x \in \Rn \text{ s.t } x_i = \frac{x^0_i}{1- \frac{\lambda_i}{\overline{\lambda}_k}} \text{if } i \notin K_{k}, \right .  \\ 		\label{eq:hyper_sphere} \numberthis
				 & \left . \sum_{i\in K_{k}} x_i^2 =  \frac{1}{\overline{\lambda}_k} \left ( 1- \sum\limits_{j\in I \setminus K_{k}} \lambda_j \left( \frac{x^0_j}{1-\frac{\lambda_j}{\overline{\lambda}_k}} \right)^2 \right ) \right \}
	\end{align*}
and every point belonging to this hypersphere is a KKT point. Moreover, all the points in this hypersphere achieve the same value for the objective function of~\cref{eq:proj_quadric_std}. Hence, for the purpose of finding one of the optimal solutions of~\cref{eq:proj_quadric_std}, we can keep in our list of candidates just one element of~\cref{eq:hyper_sphere}. In particular, we can arbitrarily select \emph{one} solution that lies in the first orthant by setting to zero all components of $K_{k}$ except one ($k'$):
	\begin{equation}
		\label{eq:x_d_degenerate_full}
		(\x^{\mathrm{d}}_k)_i = \begin{cases}
			\frac{x^0_i}{1- \frac{\lambda_i}{\overline{\lambda}_k}} & \text{if } i \notin K_{k} ,\\
			\sqrt{ \frac{1}{\overline{\lambda}_k} \left ( 1- \sum\limits_{j\in I \setminus K_{k}} \lambda_j \left( \frac{x^0_j}{1-\frac{\lambda_j}{\overline{\lambda}_k}} \right)^2 \right )}		      & \text{if } i=k' ,
			\\
			0																	      & \text{if } i \in K_{k}, i\neq k' . \\
	\end{cases} 
	\end{equation}
	Any $k' \in K_k$ works, let us choose without loss of generality $k' := \min\limits_{i \in K_k} i$.
	\index{subspace}
	As a matter of fact, this is equivalent to restricting the search to the subspace $\{ \x \in \Rn \st{x_i = 0 \quad \forall i\in K_{k} \setminus \{k'\}  } \}$, because all solutions of the hypersphere have the same objective.
	In this subspace, the problem is analogous to the case $\size{K_{k}}=1$, \ie, the intersection between a line and a quadric.

	\subsection{Bringing everything together}%
	\label{sub:Bringing_everything_together}
	Let us give a full characterization of an optimal solution to~\cref{eq:proj_quadric_std}.

	\begin{proposition}\label{prop:optimal_full}
		There is an optimal solution of~\cref{eq:proj_quadric_std} in the set $ \{\bm x(\mu^*)\} \bigcup \bm X^d $ where

		\begin{itemize}
			\item $\x(\mu^*)$ is defined by~\cref{eq:x_not_sphere}, where $\mu^*$ is the unique root of $f$ on $\mathcal{I} = ]e_1, e_2[$, and $e_1$, $e_2$ are given by~\cref{eq:e_1_e_2_degenerate};
			\item $\bm X^{\mathrm{d}} := \left \{ \x^{\mathrm{d}}_{k} \st{k = 1,\hdots,\size{\overline{\bm \lambda}}, \text{and } \size{K_k}>0} \right \}$ as defined in~\cref{sub:degenerate_case}.
		\end{itemize}

	\end{proposition}
	\begin{proof}
		Since the quadric is central, no point fulfils the LICQ condition. Hence, the solution of~\cref{eq:proj_quadric_std}--- which exists by~\cref{prop:existence_projection_quadric}--- must be a KKT point.
		The KKT points are the solutions of~\cref{eq:grad_lagrangian_not_sphere}, \ie,
		\[ x_i (1+\mu\lambda_i) = x^0_i \quad \text{ for all } i \in I \text{ and } \sum_{i\in I} \lambda_i x_i^2 =1 \]

		Hence $(\bm x^*, \mu^*)$ is a solution of the KKT conditions~\cref{eq:grad_lagrangian_not_sphere} if and only if one of the following holds:
		\begin{itemize}
			\item[(i)] $\mu^* \neq -\frac{1}{\lambda_i}$ for all $i$, $\mu^*$ is a root of $f$ defined in \cref{eq:f_mu}, and $\bm x^*$ satisfies~\cref{eq:x_not_sphere}.
			\item[(ii)] $\mu^* = -\frac{1}{\lambda_k}$ for some $k$, $x^0_i=0$ for all $i$ such that $\lambda_i = \lambda_k$, and, letting $K_k$ be the set of those $i$'s, $\bm x^* \in \pi \bigcap \Q$  defined in~\cref{eq:hyper_sphere}.
		\end{itemize}
		In case (i), we have seen that the smallest objective of~\cref{eq:proj_quadric_std} is given by the--- possibly nonexistent---unique root of $f$.
		In case (ii), all the points in~\cref{eq:hyper_sphere}---which may be empty---achieve the same value for the objective of~\cref{eq:proj_quadric_std}, and~\cref{eq:x_d_degenerate_full}---defined if and only if~\cref{eq:hyper_sphere} is nonempty---is one of those points. Finally, recall that~\cref{prop:existence_projection_quadric} proves the existence of one optimal solution to~\cref{eq:proj_quadric_std}: either case (i) or case (ii) will provide a solution.

	\end{proof}

	The full procedure to compute the projection of any point $\bm x^0$ onto a nonempty and non-cylindrical central quadric is given in~\cref{alg:exact_projection_quadric}.

\begin{algorithm}
	\caption{Exact projection onto a non-cylindrical central quadric in normal form}
	\label{alg:exact_projection_quadric}
	\begin{algorithmic}
		\REQUIRE $\bm \lambda$, the eigenvalues corresponding to~\cref{eq:proj_quadric_std}, and $\x^0 \in \Rn_+$

		\STATE $e_1 \gets  \max \left \{ \max_{\{i\in I \st{\lambda_i>0, x^0_i \neq 0}\}} -\frac{1}{\lambda_i}, -\infty \right \}$  
		\COMMENT{If the inner $\max$ is empty, $e_1 = -\infty$}
		\STATE $e_2 \gets \min \left \{ \min_{\{i\in I \st{\lambda_i < 0, x^0_i \neq 0}\}} -\frac{1}{\lambda_i}, +\infty \right \}$  
		\COMMENT{If the inner $\min$ is empty, $e_2 = +\infty$}
		\STATE $\bm D \gets \textrm{diag}(\bm \lambda)$
		\STATE $\x(\mu) \gets (\bm I + \mu \D)^{-1}  \x^0 $
		\STATE $f(\mu) \gets  \sum\limits_{i=1, x^0_i \neq 0}^n \lambda_i \left(\frac{x_i^0}{1+\mu \lambda_i}\right )^2 -1 $

		\IF{$e_1 \neq -\infty$} 
		\IF{$e_2 = +\infty$}
		\STATE $\mu_0 \gets \textrm{bisection}(f, e_1)$
		\COMMENT{Use bisection search (see~\cite[Chapter 2.1]{Burden01}) to find $\mu^0 > e_1$ with $f(\mu^0) > 0$}
		\STATE $\mu^* \gets \textrm{root}(f, \mu_0)$ \COMMENT{Using Newton with starting point $\mu^0$}
		\ELSE
		\STATE $\mu^* \gets \textrm{root}(f, e_1, e_2)$ \COMMENT{Using~\cref{alg:double_newton}}
		\ENDIF	
		\ENDIF	\COMMENT{See contrapositive of~\cref{prop:uniqueness_root_hyperboloid_degenerate}: $e_1 = -\infty \Rightarrow f$ has no root on $\mathcal{I}$}
		\STATE $\overline{\bm \lambda} \gets \textrm{unique}(\bm\lambda)$
		\STATE $\bm X^\mathrm{d} \gets [\,]$
		\FOR{$k = 1, \hdots, \size{\overline{\bm \lambda}}$}
		\STATE 	$K_{k} \gets \left \{ i \in I \st{\lambda_i = \overline{\lambda}_k, x^0_i = 0} \right  \} $
		\STATE 	$L_{k} \gets \left \{ i \in I \st{\lambda_i = \overline{\lambda}_k} \right  \} $
	
		\IF{$K_k = L_k \textbf{ and } 	  \frac{1}{\overline{\lambda}_k} \left ( 1- \sum\limits_{j\in I \setminus K_{k}} \lambda_j \left( \frac{x^0_j}{1-\frac{\lambda_j}{\overline{\lambda}_k}} \right)^2 \right )  > 0$ }
		\STATE $\x^{\mathrm{d}}_k \gets$ \cref{eq:x_d_degenerate_full}
		\STATE $\bm X^\mathrm{d} \textrm{.append(}\x^{\mathrm{d}}_k $)
		\ENDIF
		\ENDFOR
		\RETURN $\argmin\limits_{\x \in  \left \{ \x(\mu^*) \right \}  \bigcup \bm X^{\mathrm{d}}}   \norm{\x^0 - \x}_2   $
		\COMMENT{The $\min$ is taken over at most $n+1$ values}
	\end{algorithmic}
\end{algorithm}

\begin{figure}
	\centering
	\begin{subfigure}[t]{.475\textwidth}
		\centering
		\includegraphics[width=\linewidth]{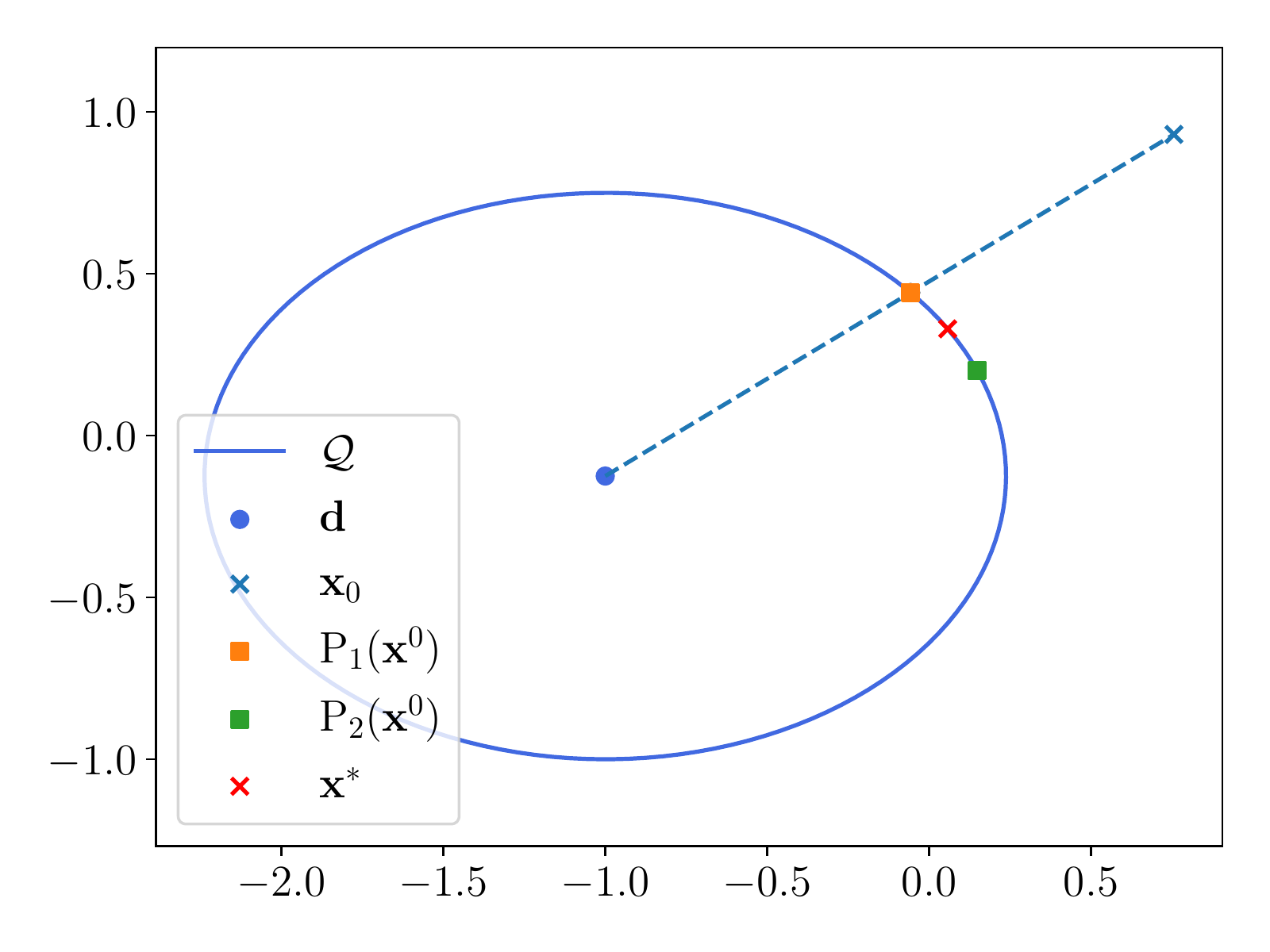}
		\caption{Illustration of the exact and quasi-projections for an ellipse.}
		\label{fig:exact_vs_quasi_proj_ellipse}
	\end{subfigure}%
	\hfill
	\begin{subfigure}[t]{.475\textwidth}
		\centering
		\includegraphics[width=\linewidth]{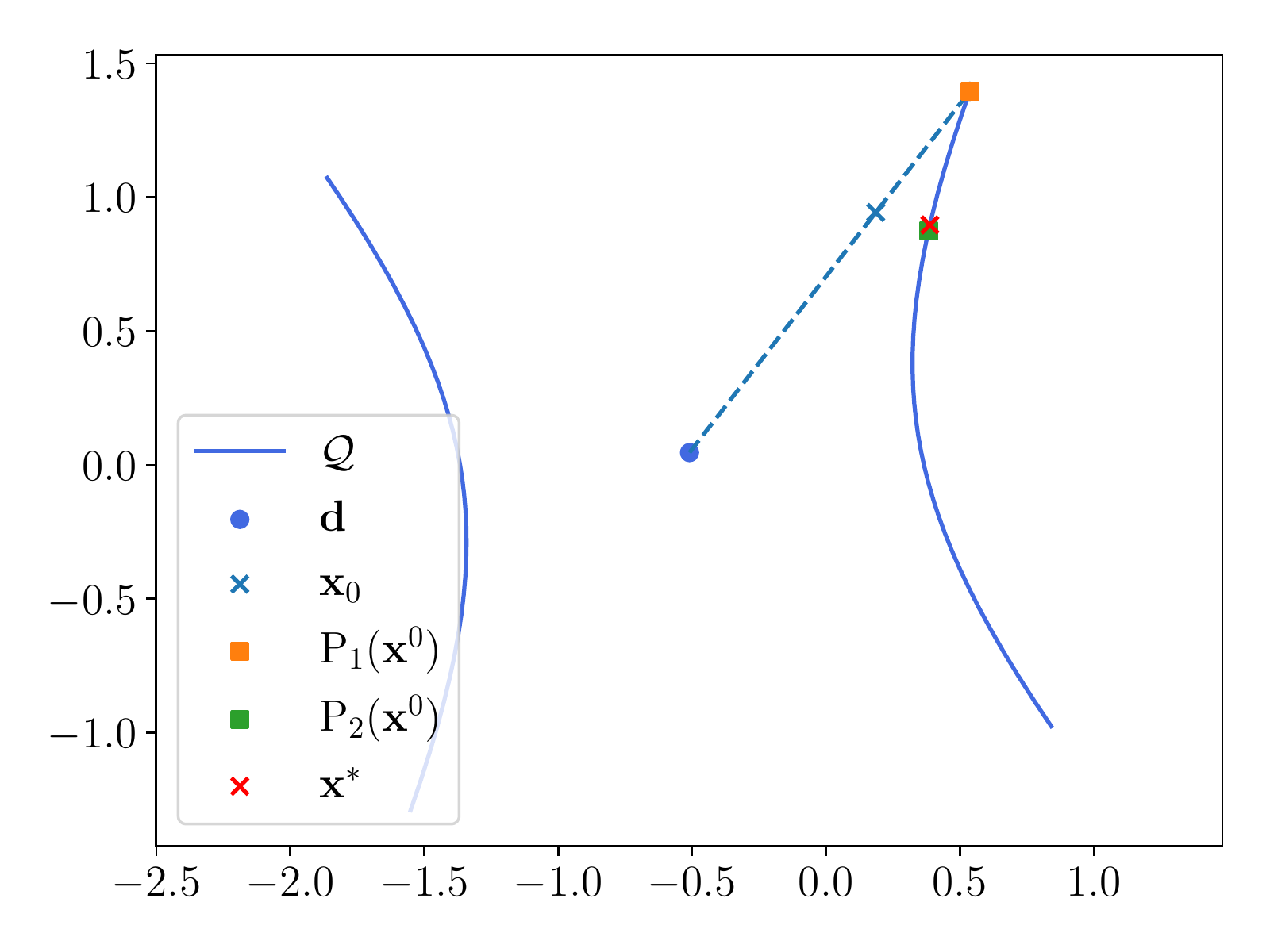}
		\caption{Illustration of the exact and quasi-projections for an hyperbola.}
		\label{fig:exact_vs_quasi_proj_hyperbole}
	\end{subfigure}
	\caption{Comparison between the exact and quasi-projections onto a 2D quadric.}
	\label{fig:comparison_quasi_exact_2D}
\end{figure}

\subsection{Quasi-projection onto the quadric}%
\label{sub:Quasi-projection_on_the_quadric}

The procedure detailed in~\cref{alg:exact_projection_quadric} is an exact projection, but it requires computing the full eigenvalue decomposition of $\A$, including the eigenvectors, which may be expensive for problems of large dimension. In this subsection, we detail a geometric procedure which allows us, under some conditions (see~\cref{ssub:Failure_of_the_quasi-projection}), to map a given point to the feasible set $\Q$ of~\cref{eq:proj_quadric}. We refer to such a mapping as a \emph{quasi-projection}.

\begin{definition}{Quasi-projection.}
	Let $\x \in \Rn$, a \emph{quasi-projection} on the quadric $\Q$ is a mapping \[ \mathrm{P} : \mathcal{D} \to \Q : \x \mapsto \mathrm{P}(\x) \, ,  \]
	where $\mathcal{D}$ is a nonempty subset of $\Rn$.
\end{definition}

Note that this definition is broad, and includes the projection operator. Ideally, $\mathcal{D}$ should be $\Rn$, but we allow the quasi-projection to fail to map some points.

This quasi-projection is inspired from the retraction in~\cite{BorSelBouAbs2014,vh22}, and from the following observation: the projection of a given point $\x$ onto a sphere that can be analytically computed by looking at the intersection between the sphere and the half line defined by the sphere centre and $\x$. 
As the quadric $\Q$ is by assumption a \emph{central quadric}, it is tempting to approximate the projection by the same mechanism described above for the sphere. This yields a first variant of the \emph{quasi-projection}. The second variant is obtained by searching for the intersection between the quadric and the line passing through $\x^0$ along the direction $\nabla \Psi(\x^0)$.

We are looking at the intersection between the quadric and the line starting from $\x^0$ along some direction $\bm \xi$. The intersections are parametrized as $\x^0 + \beta \bm \xi$, where $\beta$ satisfies $\Psi(\x^0 + \beta \bm{\xi})=0$, or equivalently, $b_1 \beta^2 + b_2 \beta + b_3 = 0$, for appropriate $b_i$'s. The $b_i$'s are given in~\cref{alg:quasi-projection-quadric}, see~\cite[\S 3.2]{BorSelBouAbs2014} for more details. 
In order to select the point that is the closest to $\x^0$ among both intersections, the $\beta$ which is the closest to zero is chosen.

We detail two variants of our \emph{quasi-projection}:
\begin{itemize}
	\item $\bm{\xi} = \d - \x^0$: this is analogous to the retraction used in~\cite{BorSelBouAbs2014}, is referred to as \emph{centre-based} quasi-projection, and is denoted by $\mathrm{P}_1$;
	\item $\bm{\xi} = \nabla \Psi(\x^0) = 2 \A \x^0 + \Ab$: the direction is given by the gradient of the level curve of $\Psi$ at $\x^0$. We refer to it as \emph{gradient-based} quasi-projection and denote it as $\mathrm{P}_2$.
\end{itemize}
The quasi-projection procedure is given in~\cref{alg:quasi-projection-quadric} and depicted in~\cref{fig:comparison_quasi_exact_2D} for both strategies. 


\begin{algorithm}
	\caption{Quasi-projection onto the quadric}
	\label{alg:quasi-projection-quadric}
	\begin{algorithmic}
		\REQUIRE{$\x^0 \in \Rn$, a non-cylindrical central quadric $\Q$ with centre $\d$, a direction $\bm \xi$}
		\STATE{$b_1 \gets \Tr{\bm \xi}  \A \bm \xi$}
		\STATE{$b_2 \gets 2 \Tr{\x^0}  \A \bm \xi + \Tr{\Ab} \bm \xi $}
		\STATE{$b_3 \gets \Tr{\x^0}  \A \x^0 + \Tr{\Ab} \x^0 + \Ac$}
		\IF{$b_2^2 - 4 b_1 b_3 < 0$}
		\RETURN \texttt{None} \COMMENT{it may happen that no intersection is available see, \eg,~\cref{fig:hyperbole_KO}}
		\ELSE
		\STATE{$\Delta \gets \sqrt{b_2^2 - 4 b_1 b_3}$}
		\STATE $\beta^+ \gets  \frac{-b_2 + \Delta}{2 b_1} $
		\STATE $\beta^- \gets  \frac{-b_2 - \Delta}{2 b_1} $
		\IF{$b_2 > 0$}
		\STATE $\beta \gets \beta^+  $
		\ELSE
		\STATE $\beta \gets \beta^-$
		\ENDIF
		\COMMENT{we select the closest to $\x^0$ of the two intersections}
		\RETURN $\x^0 + \beta \bm \xi$
		\ENDIF
	\end{algorithmic}
\end{algorithm}

\subsubsection{Failure of the quasi-projection}%
\label{ssub:Failure_of_the_quasi-projection}

\begin{figure}
\centering
\begin{subfigure}{.5\textwidth}
  \centering
  \includegraphics[width=\linewidth]{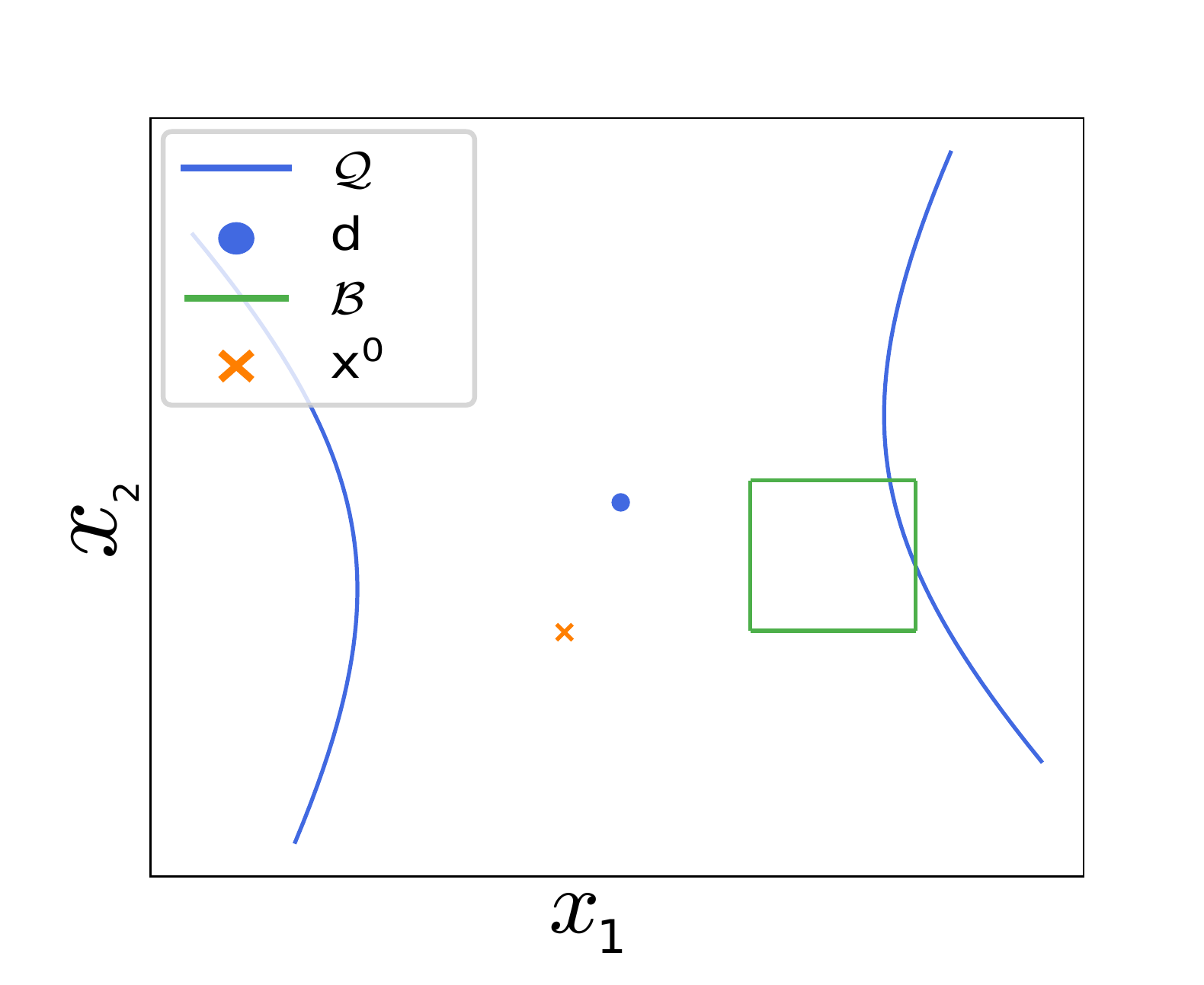}
  \caption{2D (hyperbola) case.}
  \label{subfig:projection_KO_2D}
\end{subfigure}%
\begin{subfigure}{.5\textwidth}
  \centering
  \includegraphics[width=\linewidth]{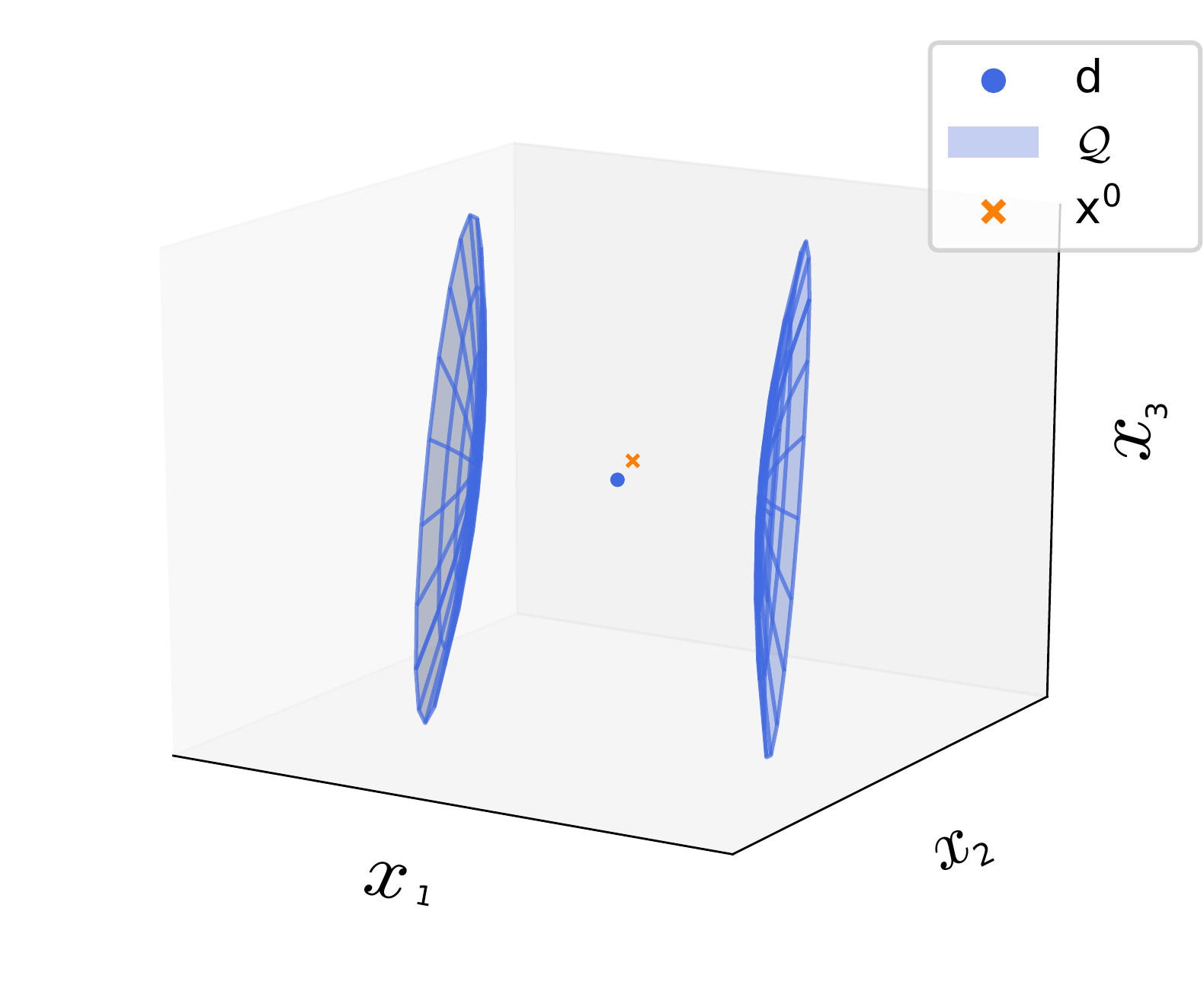}
  \includegraphics[width=\linewidth]{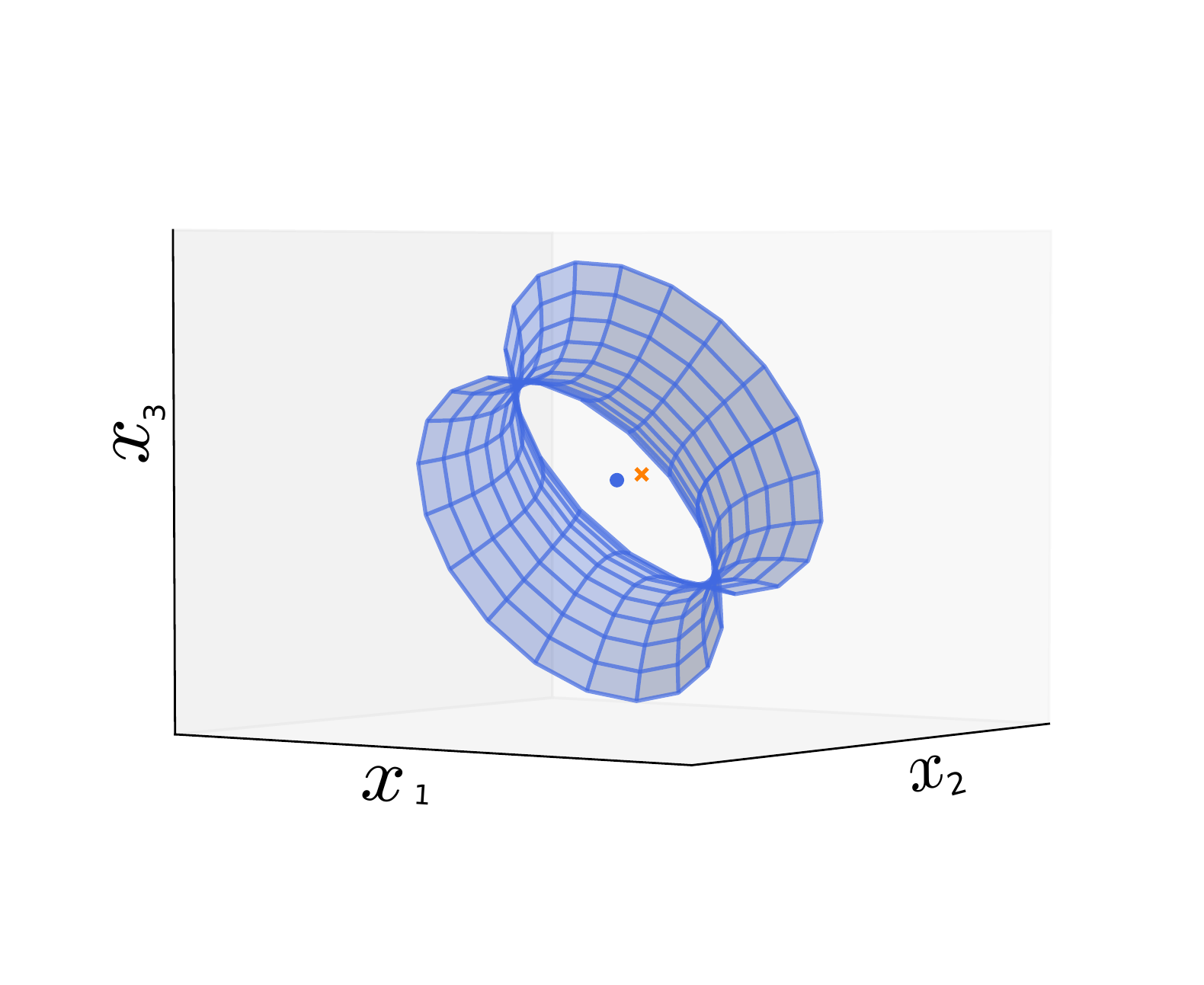}
  \caption{3D (two and one-sheet hyperboloid) cases.}
  \label{subfig:projection_KO_3D}
\end{subfigure}
\caption{Illustration of a failure of the quasi-projection ($\mathrm{P}_1$): the point $\bm x^0$ cannot be mapped to the quadric using~\cref{alg:quasi-projection-quadric}. Indeed, the line defined by $\bm x^0$ and $\d$ does not intersect with the quadric.}
\label{fig:hyperbole_KO}
\end{figure}

Remark that for $\x^0 = \bm{0}$, $\mathrm{P}_1$ is not defined, and for the hyperboloid case, the set where $\mathrm{P}_1$ is not defined ($\Rn \setminus \mathcal{D}$), is a closed set including $\bm{0}$. Examples of nontrivial points that cannot be mapped using $\mathrm{P}_1$ are provided in~\cref{fig:hyperbole_KO}. Indeed, in these cases, there is no intersection between the quadric and the line starting from $\x^0$. We tackle this issue by resorting to the exact projection from~\cref{alg:exact_projection_quadric} whenever this situation occurs.

There are points $\x^0$ for which $\mathrm{P}_2$ returns \texttt{None}, but it returns a point when $\x^0$ is close enough to $\Q$.


\subsubsection{Features of the quasi-projection}%
\label{ssub:Features_of_the_quasi-projection}

In general the quasi-projection is not exact, in the sense that the resulting point is \emph{not} the optimal solution of~\cref{eq:proj_quadric}. However, we expect the quasi-projection to be close to optimality when the point is close enough to the quadric. Such a behaviour is observed in our simulations in~\cref{sub:Alternate_projection_versus_Gurobi}. Also, in the specific case when the quadric is a sphere, then both $\mathrm{P}_1$ and $\mathrm{P}_2$ solve~\cref{eq:proj_quadric}.

\section{Splitting methods for the projection onto the intersection of a box and a quadric}%
\label{sec:Splitting_methods}

This section is devoted to the analysis of the projection onto a feasible set, $\Omega$, which is the intersection between a box and a non-cylindrical central quadric.
Let $\mathcal{B}$ be a nonempty $n$-dimensional hyper-cube or \emph{box}, aligned with the axes:
\begin{equation}
	\label{eq:box}
	\mathcal{B} = \left \{ \bm x \in \Rn \st{ \xmin_i \leq x_i \leq \xmax_i, \, \forall i=1\hdots n } \right \},
\end{equation}
for given lower and upper bounds $\bm \xmin$ and $\bm \xmax$, and $\Q$ a quadric. The optimization problem at hand reads
\begin{align*} \label{eq:P}  
	\min_{\bm x \in \Rn} 	& \norm{\bm x - \bm x^0}_2^2 \numberthis \\ 
	\text{s.t. } 	& \bm x \in \mathcal{B} , \\
			& \bm x \in \Q .
\end{align*}

Note that what is developed in this article can be easily extended to a polytope $\mathcal{P}$ and a Cartesian product of $\size{T}$ quadrics $\times_{t = 1 \hdots \size{T}} \mathcal{Q}(\Psi_t)$. This is discussed in~\cref{sub:Extensions_and_further_researches}.

In particular, we study two splitting algorithms: the Douglas-Rachford (DR) scheme, and the alternating projection method (AP). We consider three variants of the latter: one based on the exact projection and two based on the quasi-projection from~\cref{sub:Quasi-projection_on_the_quadric}, that approximates the projection via a geometric construction.
Splitting algorithms exploit the separable structure of the problem, since the projection onto each of the sets that define the intersection, $\Omega:= \mathcal{B} \bigcap \Q$, is easy to compute. They recently have been widely studied, and perform particularly well on certain classes of nonconvex problems.

A first convergence result for the (local) solution of the alternating projections in the nonconvex setting is presented in~\cite[Theorem 3.2.3]{drusvyatskiy_slope_2013} for sets that intersect transversally. This result is particularized to (nonempty and closed) semi-algebraic intersections in~\cite[Theorem 7.3]{drusvyatskiy_transversality_2015} which we use in this work. A second result that is used in this paper is~\cite[Corollary 1]{li_douglasrachford_2016}, which is a convergence result for a (modified) Douglas-Rachford splitting. These two important propositions exploit the Kurdyka-Łojasiewicz properties of the indicator function of semi-algebraic sets, and are unfortunately local results: the theorem from~\cite{drusvyatskiy_transversality_2015} assumes that the starting point is \emph{near} the intersection of the two considered sets, and the theorem from~\cite{li_douglasrachford_2016} proves the convergence to a \emph{stationary point} of the problem of minimizing the distance to one of the sets, subject to being in the second set. We present an example where both methods fail to converge to a feasible point (\cref{fig:pathological_AP,fig:pathological_DR}), and propose a restart heuristic in~\cref{sub:Implementation_details}. 

The structure is the following: \S~\ref{sub:Projection_on_the_box} briefly recalls how to project onto a box, \S~\ref{sub:Alternate_projection_method_AP} details the AP  methods, and \S~\ref{sub:Douglas-Rachford_method} covers the DR splitting. A comparison table of all methods is presented in \S~\ref{sub:comparison}, a power systems application is discussed in \S~\ref{sub:Extensions_and_further_researches}, and a restart mechanism is given in \S~\ref{sub:Implementation_details}.

\subsection{Projection onto a box}%
\label{sub:Projection_on_the_box}

This projection is straightforward. Indeed given a point $\x^0$, it suffices to check for each dimension $i$ whether this point violates the lower (respectively upper) bound and replace it accordingly. This gives~\cref{alg:projection-box}.

\begin{algorithm}
	\caption{Projection onto the box~\cref{eq:box}}
	\label{alg:projection-box}
	\begin{algorithmic}
		\REQUIRE{$\x^0 \in \Rn$}
		\STATE{$\x \gets \x^0$}
		\FOR{$i = 1 \hdots n$}
			\IF{$x_i^0 < \xmin_i$}
				\STATE{$x_i \gets \xmin_i$}
			\ELSIF{$x_i^0 > \xmax_i$}
				\STATE{$x_i \gets \xmax_i$}
			\ENDIF
		\ENDFOR
		\RETURN $\x$
	\end{algorithmic}
\end{algorithm}

Note that for a more general polytope $\mathcal{P}$, the projection cannot be computed analytically. However the projection can be efficiently computed by solving the convex QP optimization problem:
\begin{align}
	\label{eq:proj_polytope}
	\min_{\x \in \mathcal{P}} \norm{\x^0 - \x}_2^2 .
\end{align}

\subsection{Alternating projection method}%
\label{sub:Alternate_projection_method_AP}

The alternating projection method can be easily built by alternately projecting onto the quadric and onto the box. This gives~\cref{alg:AP}. Depending on whether we use the exact projection or one of the two quasi-projections detailed in \S~\ref{sub:Quasi-projection_on_the_quadric}, we refer to the methods as follows: alternating projections with exact projection onto the quadric (APE), alternating projections with the centre-based quasi-projection (APC) or with the gradient-based quasi-projection (APG). $\Pr_\mathcal{X}$ stands for (one solution of) the projection onto a (non)convex set $\mathcal{X}$ and $\mathrm{P}_{\mathcal{Y}}$ the (quasi-)projection onto a set $\mathcal{Y}$.

\begin{algorithm}
	\caption{Alternating projections}
	\label{alg:AP}
	\begin{algorithmic}
		\REQUIRE $\bm x^0 \in \Rn$
		\STATE $k \gets 0$
		\WHILE{$k < \niter$ \textrm{\textbf{and not}} $\bm x^k \in \Omega$}
		\STATE $\bm y^{k+1} \gets \Pr_{\mathcal{B}}(\x^k)$
		\COMMENT{Using~\cref{alg:projection-box}}
		\STATE $\bm x^{k+1} \gets \mathrm{P}_{\mathcal{Q}}(\bm y^{k+1})$
		\COMMENT{Using~\cref{alg:exact_projection_quadric} or \cref{alg:quasi-projection-quadric}}
			\STATE $k \gets k+1$
		\ENDWHILE
		\RETURN $\bm x^k$
	\end{algorithmic}
\end{algorithm}

Assuming that the initial iterate is close enough to the intersection,~\cite{drusvyatskiy_transversality_2015} provides a convergence result for APE, which is particularized to our case in~\cref{prop:convergence_APE}. Note that~\cref{prop:convergence_APE} guarantees convergence to a point $\x^*$ in $\Omega$, but provides no guarantee about the optimality of this point, \ie, it is not true in general that $\x^* \in \argmin\limits_{\x \in \Omega} \norm{\x - \x^0}^2_2 $.

\begin{proposition}\label{prop:convergence_APE}
	If~\cref{alg:AP} with the exact projection (APE) is initialized from $\x^0 \in \Q$ and near $\mathcal{B}$, then the distance of the iterates to the intersection $\Q \bigcap \mathcal{B}$ converges to zero, and hence every limit point, $\x^*$, lies in $Q \bigcap \mathcal{B}$.
\end{proposition}
\begin{proof}
	This follows from~\cite[Theorem 7.3]{drusvyatskiy_transversality_2015}, since $\Q$ and $\mathcal{B}$ are semi-algebraic and $\mathcal{B}$ is bounded.
\end{proof}

Remark that if $\B \succ \bm 0$, then $\Q$ is also bounded and we can as well choose $\x^0 \in \mathcal{B}$ and near $\Q$.

\begin{figure}
\centering
\begin{subfigure}[t]{.475\textwidth}
  \centering
  \includegraphics[width=\linewidth]{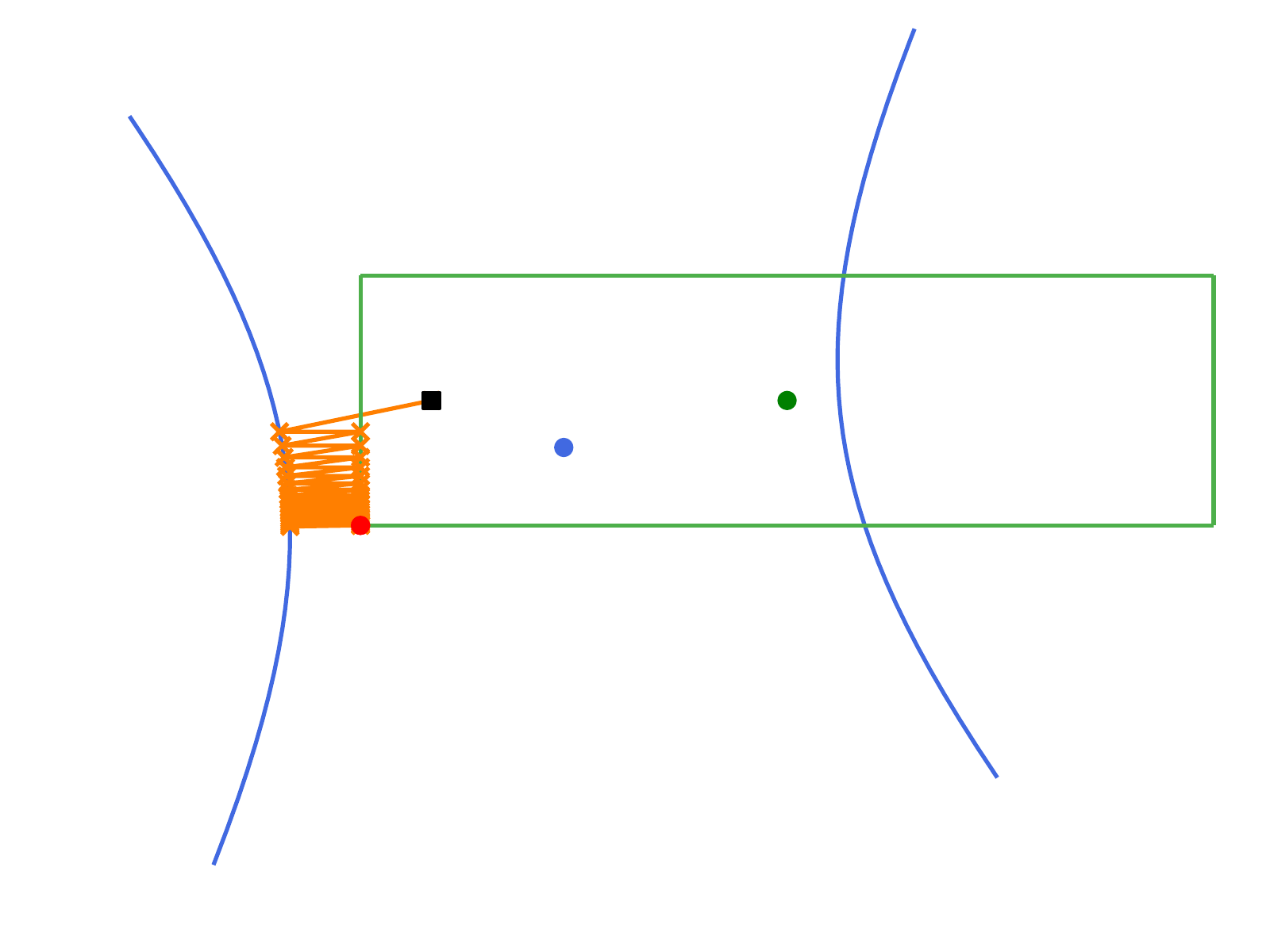}
  \caption{Alternating projections with exact projection (APE).}
  \label{subfig:pathological_APE}
\end{subfigure}%
\hfill
\begin{subfigure}[t]{.475\textwidth}
  \centering
  \includegraphics[width=\linewidth]{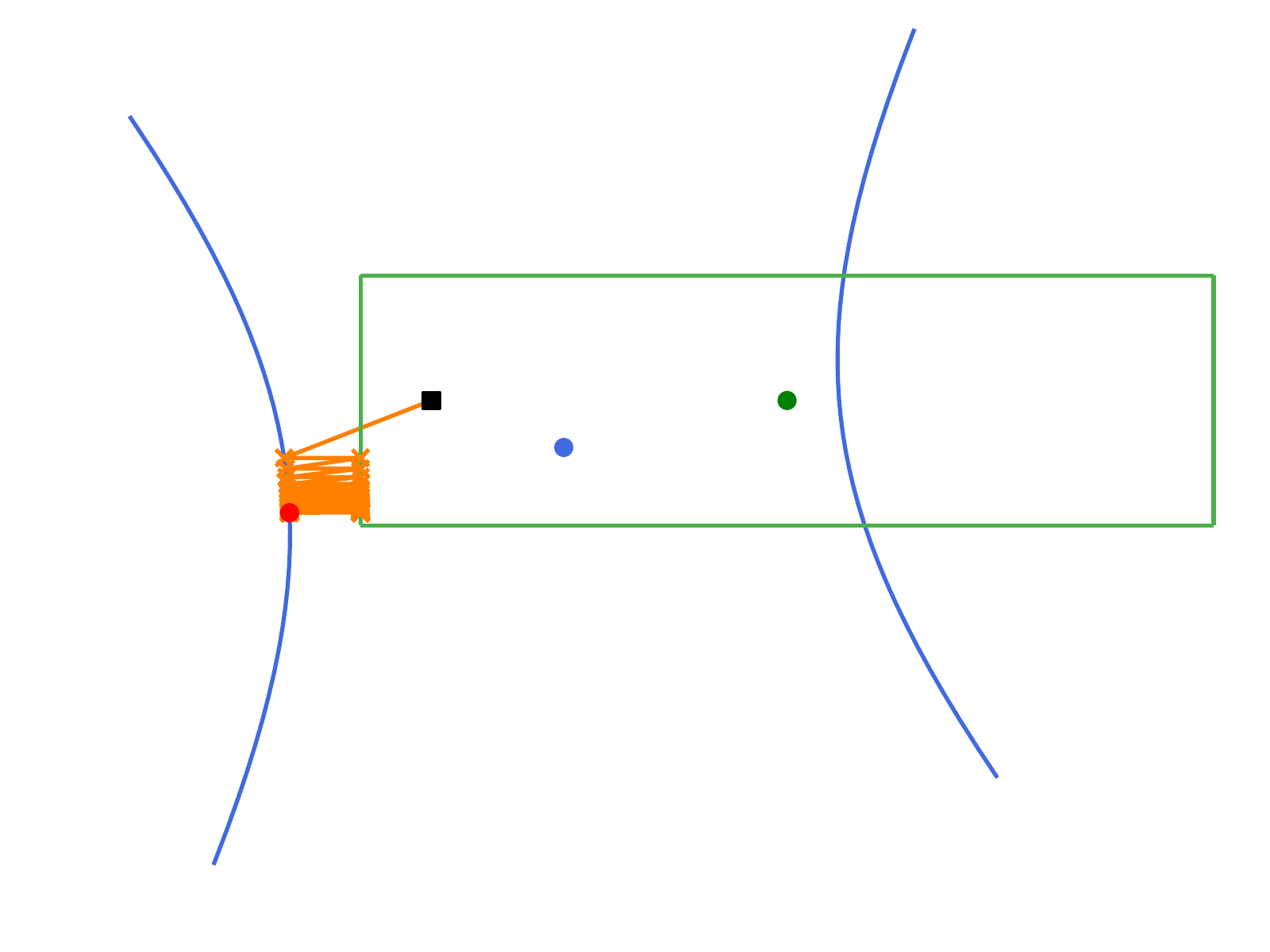}
  \caption{Gradient-based alternating projections (APG).}
  \label{subfig:pathological_APG}
\end{subfigure}
\begin{subfigure}[b]{\textwidth}
\begin{center}
  \includegraphics[width=\linewidth]{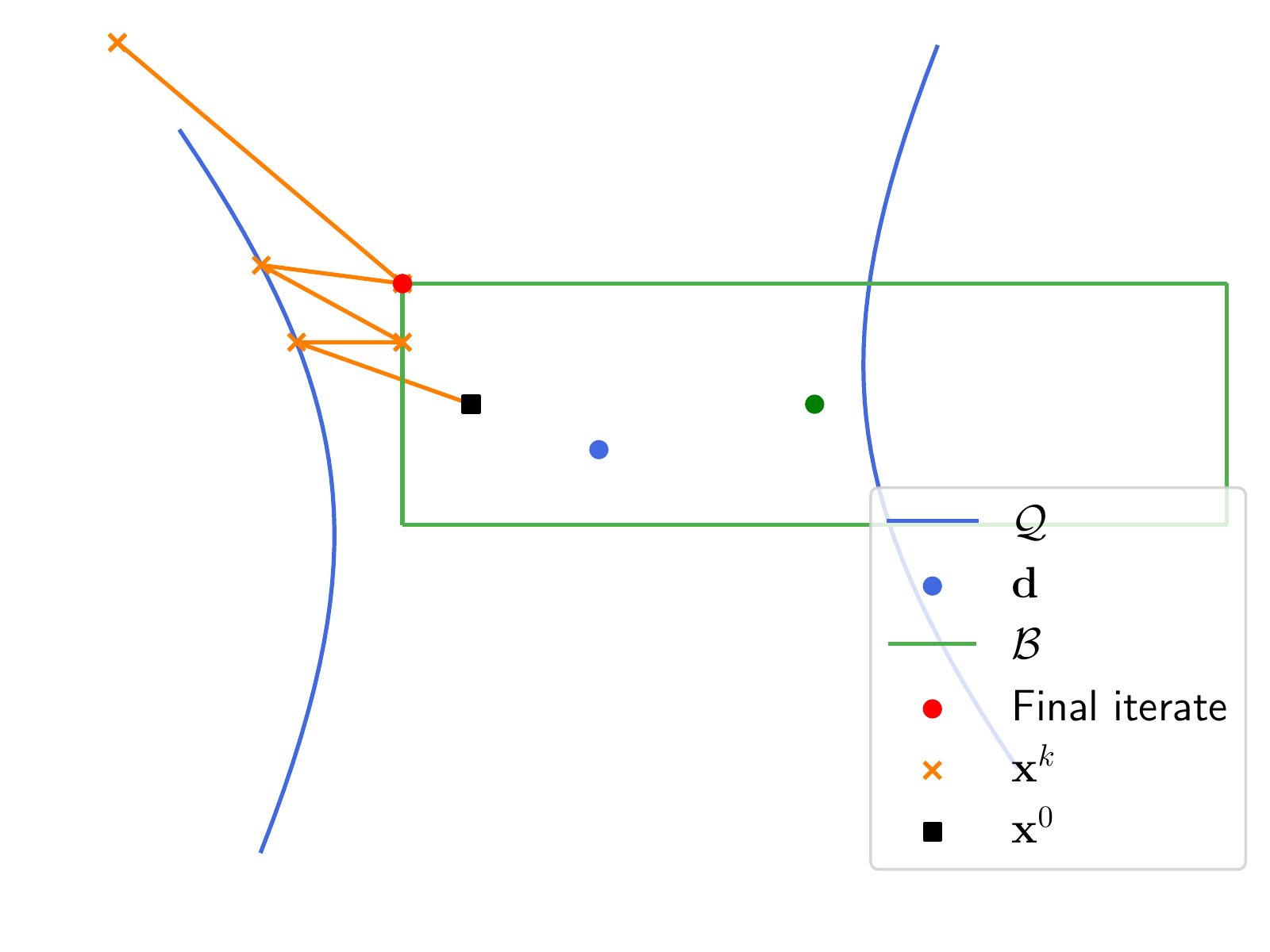}
  \caption{Centre-based alternating projections (APC).}
  \label{subfig:pathological_APC}
\end{center}
\end{subfigure}
\caption{Illustration of a (2D) pathological case where none of the proposed alternating methods converge to a feasible point of~\cref{eq:P}. We represent both $\bm x^k$ and $\bm y^k$ as orange crosses.}
\label{fig:pathological_AP}
\end{figure}

\begin{figure}
	\centering
	\begin{subfigure}[b]{.475\textwidth}
		\centering
		\includegraphics[width=\linewidth]{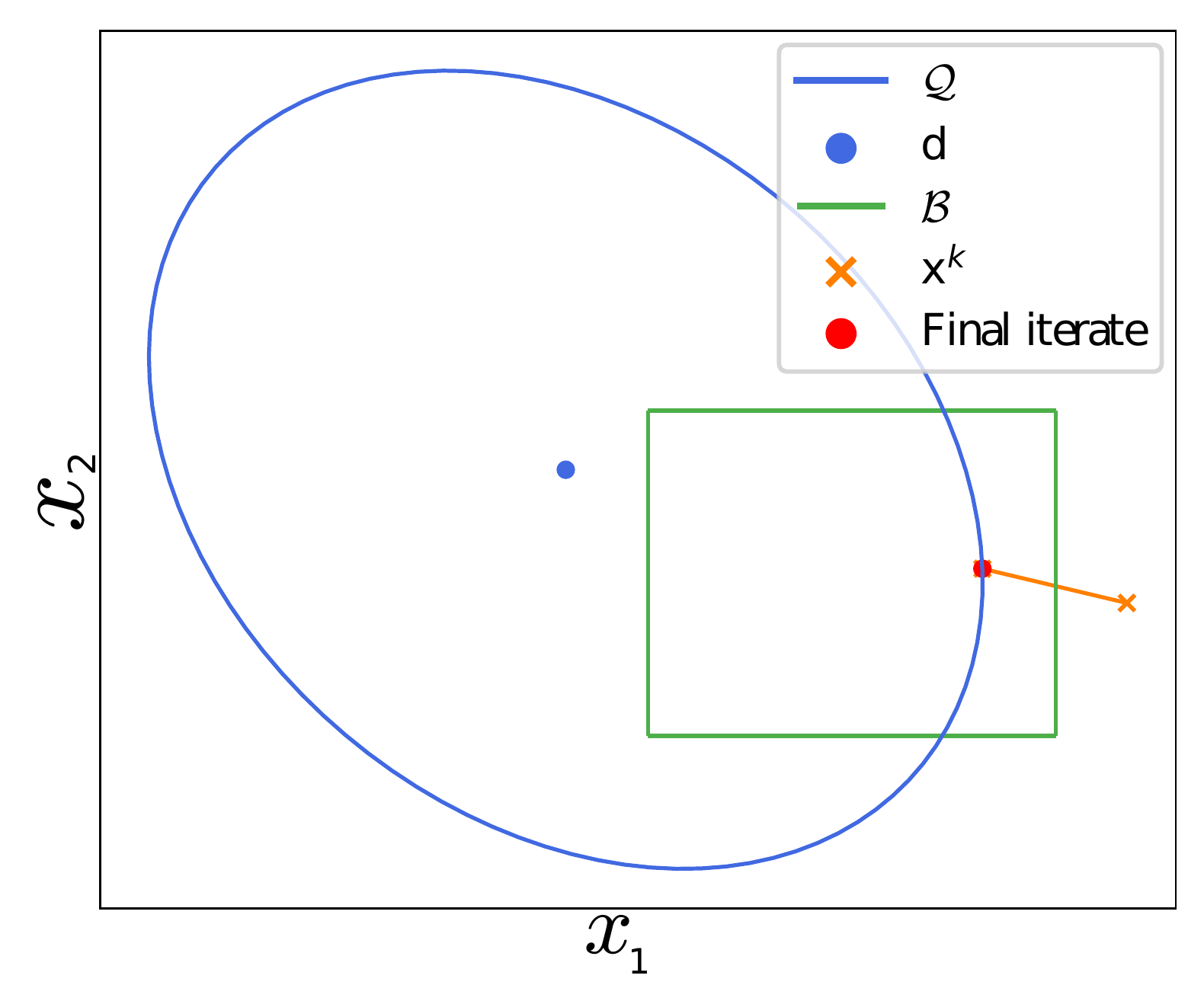}
		\caption{2D (ellipse) case.}
		\label{subfig:ellipse_oneshot}
	\end{subfigure}%
	\hfill
	\begin{subfigure}[b]{.475\textwidth}
		\centering
		\includegraphics[width=\linewidth]{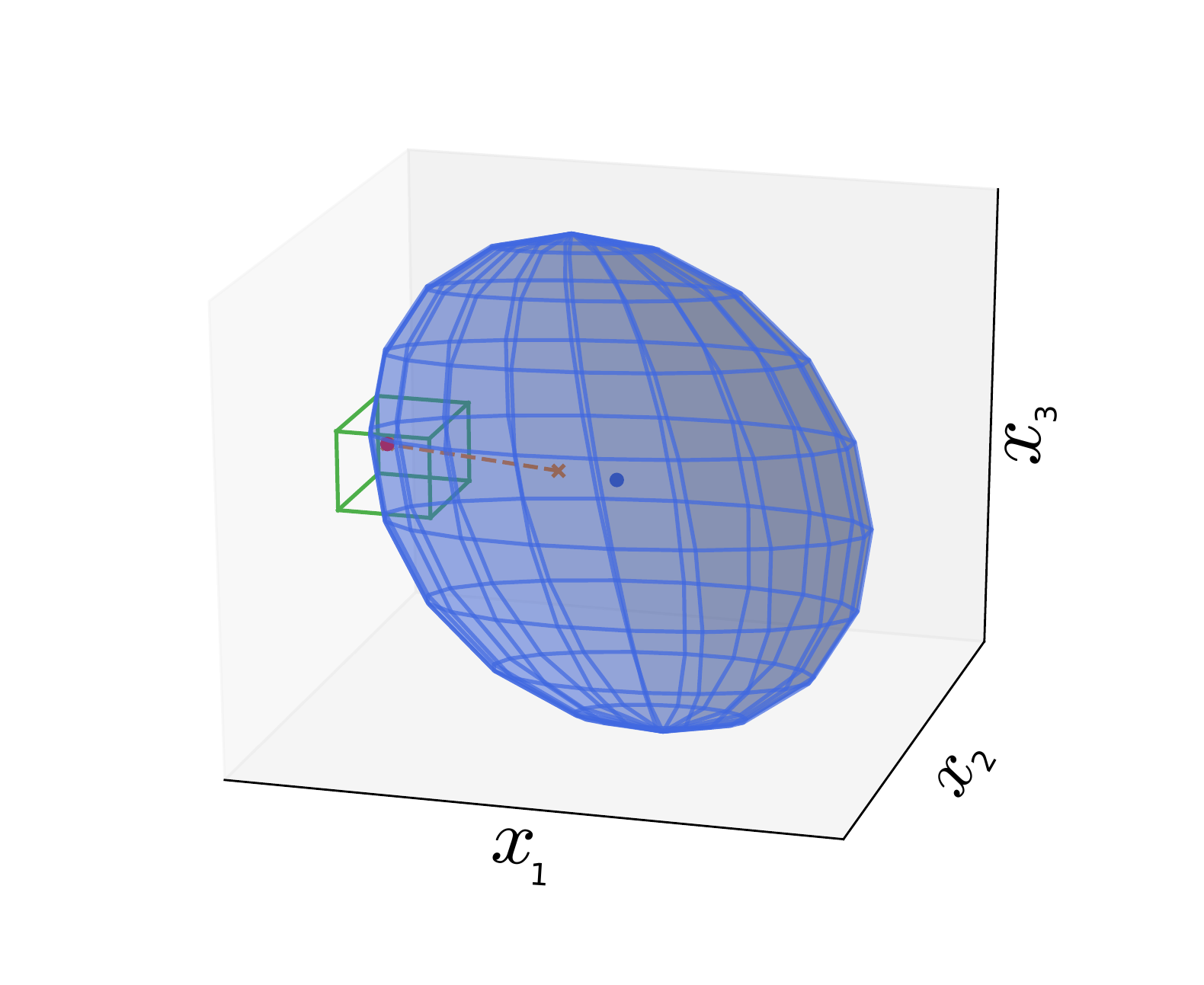}
		\label{subfig:ellipsoid_oneshot}
		\caption{3D (ellipsoid) case.}
	\end{subfigure}
	\caption{Illustration of the centre-based alternating projection method (APC). In these cases, the method converges in a single iteration as the quasi-projection from~\cref{alg:quasi-projection-quadric} yields a feasible point, \ie, a point inside the box. We represent both $\bm x^k$ and $\bm y^k$ as orange crosses.}
	\label{fig:oneshot}
\end{figure}

\begin{figure}
	\centering
	\begin{subfigure}[b]{.475\textwidth}
		\centering
		\includegraphics[width=\linewidth]{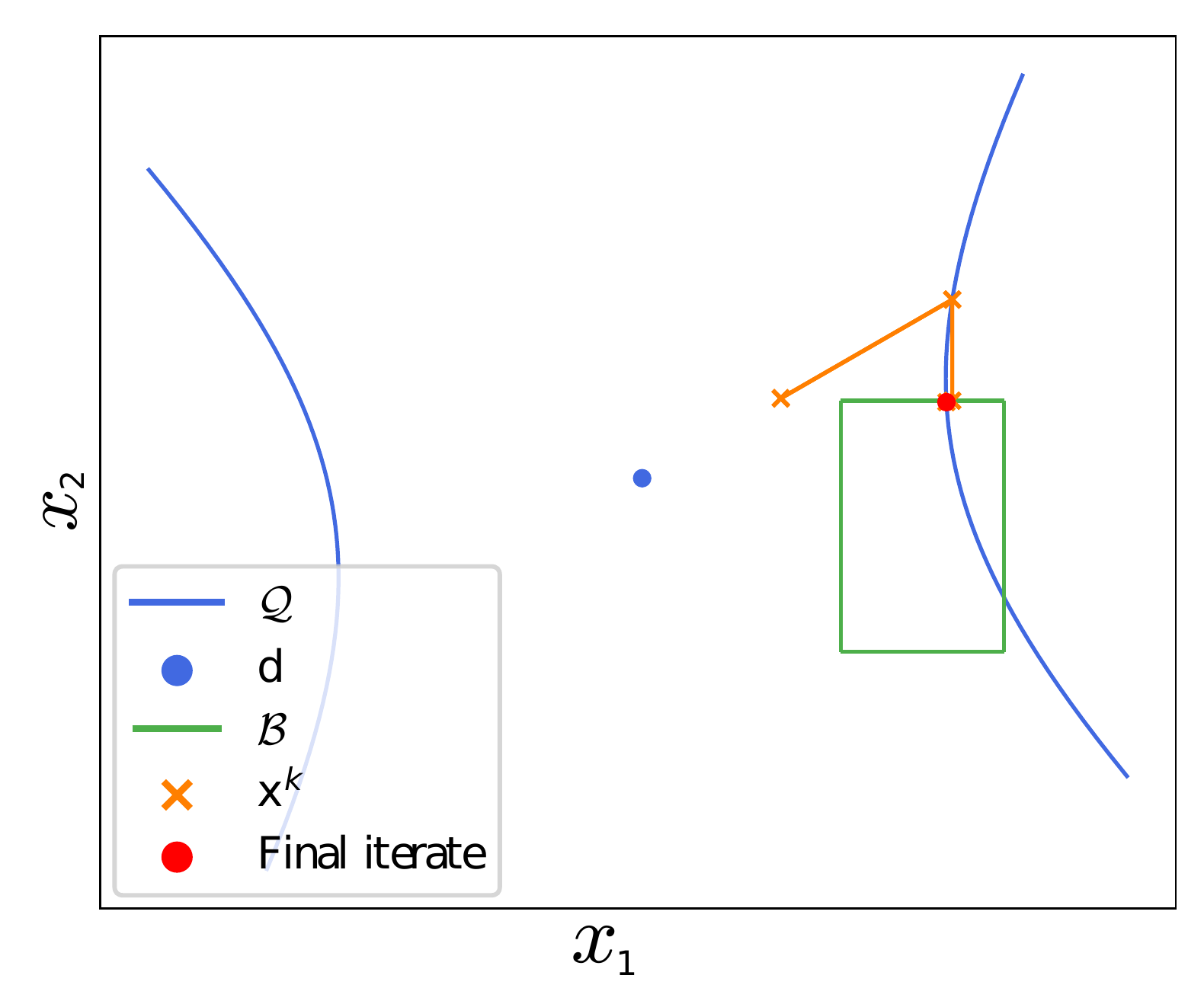}
		\caption{2D (hyperbola) case: the algorithm converges in three iterations.}
		\label{subfig:hyperbole_OK}
	\end{subfigure}%
	\hfill
	\begin{subfigure}[b]{.475\textwidth}
		\centering
		\includegraphics[width=\linewidth]{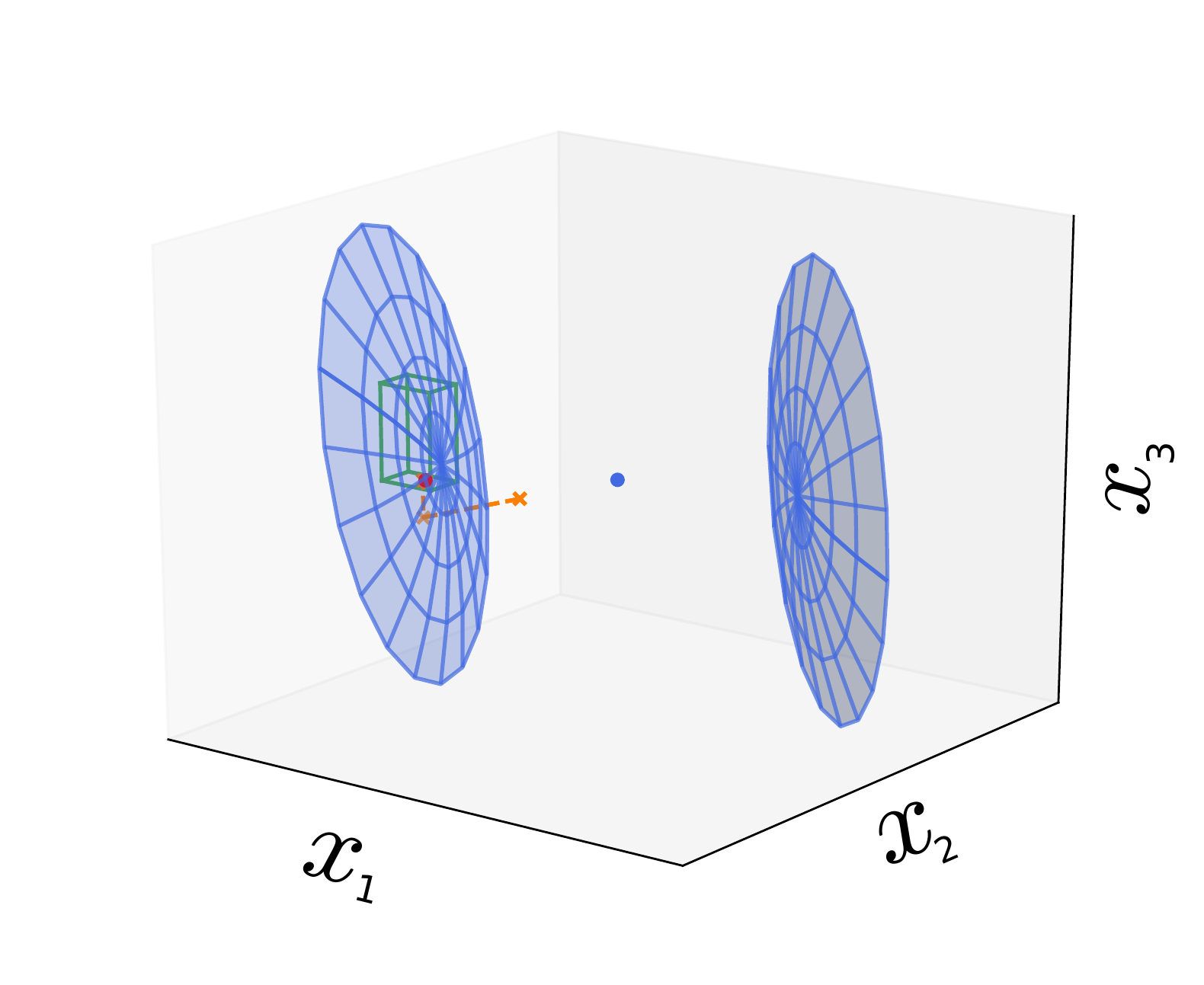}
		\caption{3D (two-sheet hyperboloid) case, the algorithm converges in three iterations.}
		\label{subfig:hyperboloid_OK}
	\end{subfigure}
	\caption{Illustration of successes of the centre-based alternating projection method (APC) on 2D and 3D hyperbolic cases. We represent both $\bm x^k$ and $\bm y^k$ as orange crosses.}
	\label{fig:hyperboloid_OK}
\end{figure}

\cref{fig:oneshot,fig:hyperboloid_OK} present examples where the alternating methods converge in a single iteration or in multiple iterations. Only APC is depicted. Notice that if APE converges in a single iteration, then the obtained solution, $\x^*$, is an optimal solution of~\cref{eq:P}, that is, $\x^* \in \argmin\limits_{\x \in \Omega} \norm{\x - \x^0}_2^2$. Figure~\ref{fig:pathological_AP} shows a pathological example where none of the alternating projection methods converge to a feasible point of~\cref{eq:proj_quadric_std}. We propose in~\cref{sub:Implementation_details} certain heuristics in order to overcome such pathological cases.


\subsection{Douglas-Rachford method}%
\label{sub:Douglas-Rachford_method}

Following~\cite{li_douglasrachford_2016}, the Douglas-Rachford splitting algorithm aims at solving 
\begin{align*}
	\min \, f(\bm x) + g(\bm x)
\end{align*}
where $f$ has a Lipschitz continuous gradient and $g$ is a proper closed function. The DR iteration starts at any $\bm y_0$ and repeats for $k = 0,1,\hdots$

\begin{align*}
	\bm x^{k+1} &= \prox_f (\bm y^k),\\ 
	\bm y^{k+1} &= \bm y^k +  \prox_g (2 \bm x^{k+1} - \bm y^k) - \bm x^{k+1},
\end{align*}
where the $\prox$ operator (with step size 1) is defined as
\begin{equation}
	\label{eq:prox}
	\prox_f(v) = \arg \min_{x\in\Rn} \left(f(x) + \frac 1 2 \|x - v\|_2^2\right) \, .
\end{equation}
Let $\mathcal{I}_{\mathcal{X}}: \mathcal{X} \mapsto \mathbb{B}$ be the indicator function of a set $\mathcal{X}$ defined as
\[ \mathcal{I}_{\mathcal{X}}(\bm x) = \begin{cases}
	0 \text{ if } x\in \mathcal{X},\\
	+\infty \text{ else.}
\end{cases} \]
If we identify $f:= \mathcal{I}_{\mathcal{B}}$ and $g := \mathcal{I}_{\mathcal{Q}}$, \ie, the indicator functions of the sets that define $\Omega$, then the DR algorithm reads

\begin{align*}
	\bm x^{k+1} &= \Pr_{\mathcal{B}}(\bm y^k)\\ 
	\bm y^{k+1} &= \bm y^k +  \Pr_{\mathcal{Q}}(2 \bm x^{k+1} - \bm y^k) - \bm x^{k+1}\\ 
\end{align*}
which can be rewritten in a compact way~\cite{bauschke_phase_2002},

\begin{equation}
	\bm x^{k+1} = \left(\Pr_{\mathcal{Q}} (2 \Pr_{\mathcal{B}} - \bm I) + (\bm I - \Pr_{\mathcal{B}}) \right ) (\bm x^k)
\end{equation}
since the proximal operator of an indicator function of a given set $X$ is the projection onto this set $\Pr_{X}$. We denote this method as DR, and explicitly state it in~\cref{alg:douglasrachford}.

\begin{algorithm}
	\caption{Douglas-Rachford splitting method (DR)}
	\label{alg:douglasrachford}
	\begin{algorithmic}
		\REQUIRE An initial point $\x^0$
		\WHILE{a termination criterion is not met}
		\STATE $\bm y^{t+1} \gets  \Pr_{\mathcal{B}}(\x^t)$
		\STATE $\bm z^{t+1} \gets \argmin_{\bm z\in \Q} \norm{2 \bm y^{t+1} - \x^t - \bm z}^2$
		\STATE $\bm x^{t+1} \gets \bm x^t + (\bm z^{t+1} - \bm y^{t+1})$
		\ENDWHILE
		\RETURN $\bm z^{t+1}$
	\end{algorithmic}
\end{algorithm}

\paragraph{Modified Douglas-Rachford}%
\label{par:Modified_Douglas-Rachford}
We now present the modification of DR splitting for the feasibility problem of~\cite{li_douglasrachford_2016}. Instead of using the indicator function for the convex set $\mathcal{B}$, the splitting is performed with the squared distance function $d^2_{\mathcal{B}}(\x) = \argmin_{\bm y \in \mathcal{B}} \norm{\x - \bm y}_2^2$, \ie,

\begin{equation}
	\label{eq:obj_dr-feasibility}
	\min_{\x \in \Q} d_{\mathcal{B}}^2(\x) \, ,
\end{equation}
which can be equivalently seen as
\begin{equation}
	\label{eq:obj_dr-feasibility_2}
	\min_{\x \in \Rn} d_{\mathcal{B}}^2(\x)  + \mathbb{1}_{\Q} (\x) \, .
\end{equation}
DR applied to~\cref{eq:obj_dr-feasibility_2} gives~\cref{alg:douglasrachford_feasibility}, denoted as DR-F.
\begin{algorithm}
	\caption{Douglas-Rachford splitting method for feasibility problems (DR-F)}
	\label{alg:douglasrachford_feasibility}
	\begin{algorithmic}
		\REQUIRE An initial point $\x^0$ and a step size parameter $\gamma >0$
		\WHILE{a termination criterion is not met}
		\STATE $\bm y^{t+1} \gets \frac{1}{\gamma +1} (\x^t + \gamma P_{\mathcal{B}}(\x^t))$
		\STATE $\bm z^{t+1} \gets \argmin_{\bm z\in \Q} \norm{2 \bm y^{t+1} - \x^t - \bm z}^2$
		\STATE $\bm x^{t+1} \gets \bm x^t + (\bm z^{t+1} - \bm y^{t+1})$
		\ENDWHILE
		\RETURN $\bm z^{t+1}$
	\end{algorithmic}
\end{algorithm}
We can use~\cite[Corollary 1]{li_douglasrachford_2016} to obtain a convergence result for the DR-F method.

\begin{proposition}\label{prop:convergence_DR-F}
	If $0<\gamma<\sqrt{\frac{3}{2}}-1$, then the sequence $\left \{ (\bm y^t, \bm z^t, \bm x^t) \right \} $ provided by~\cref{alg:douglasrachford_feasibility} converges to a point $(\bm y^*, \bm z^*, \x^*)$ which satisfies $\bm z^* = \bm y^*$, and $\bm z^*$ is a stationary point of~\cref{eq:obj_dr-feasibility}.
\end{proposition}
\begin{proof}
	Since $\Q$ and $\mathcal{B}$ are nonempty closed semi-algebraic set, with $\mathcal{B}$ being convex and compact, we satisfy the hypothesis of~\cite[Corollary 1]{li_douglasrachford_2016} for $0 < \gamma < \sqrt{\frac{3}{2}}-1$.
\end{proof}

\begin{figure}
	\centering
	\begin{subfigure}[t]{.475\textwidth}
		\centering
		\includegraphics[width=\linewidth]{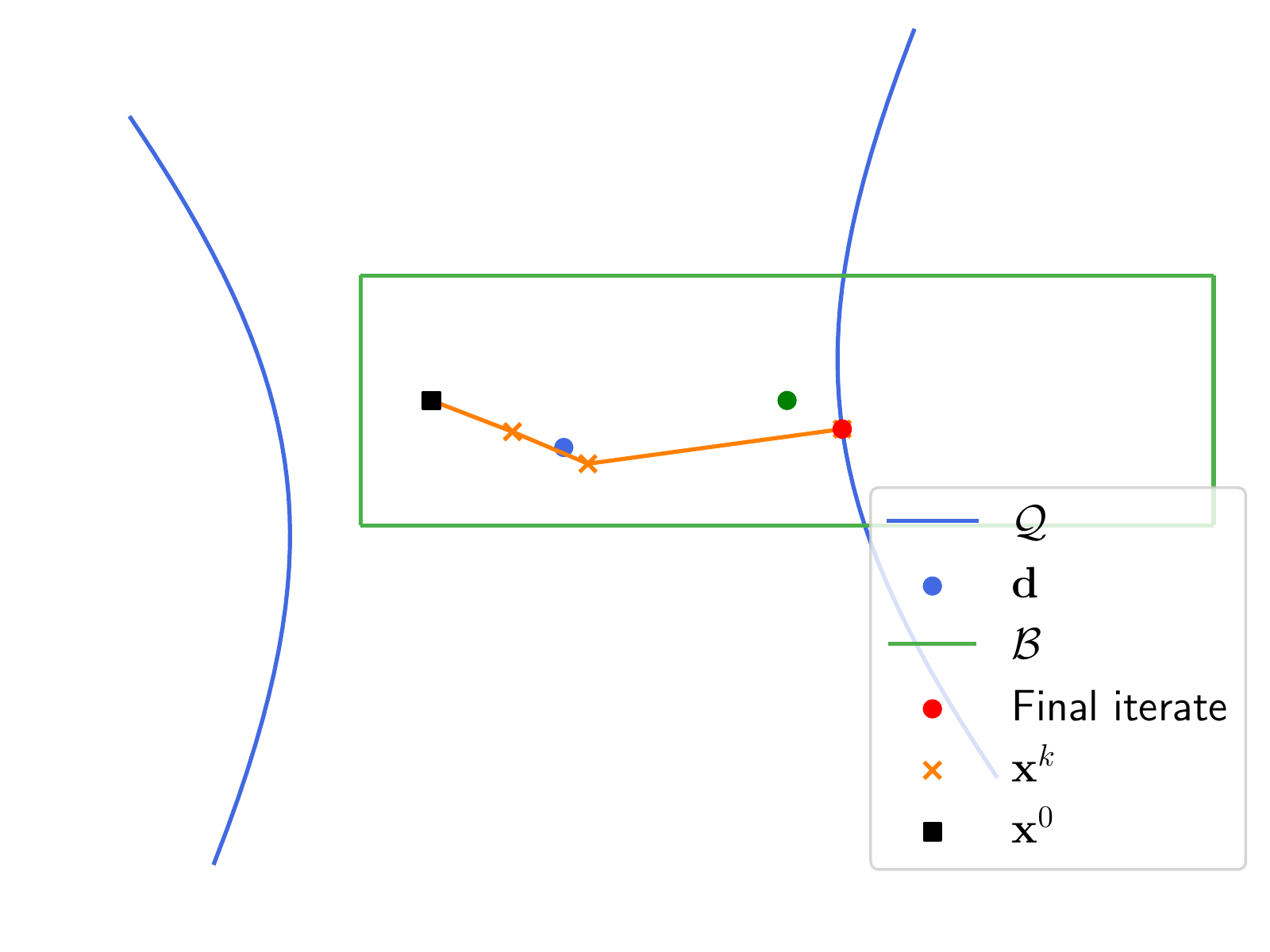}
		\caption{Douglas-Rachford (DR).}
		\label{subfig:DR_OK}
	\end{subfigure}%
	\hfill
	\begin{subfigure}[t]{.475\textwidth}
		\centering
		\includegraphics[width=\linewidth]{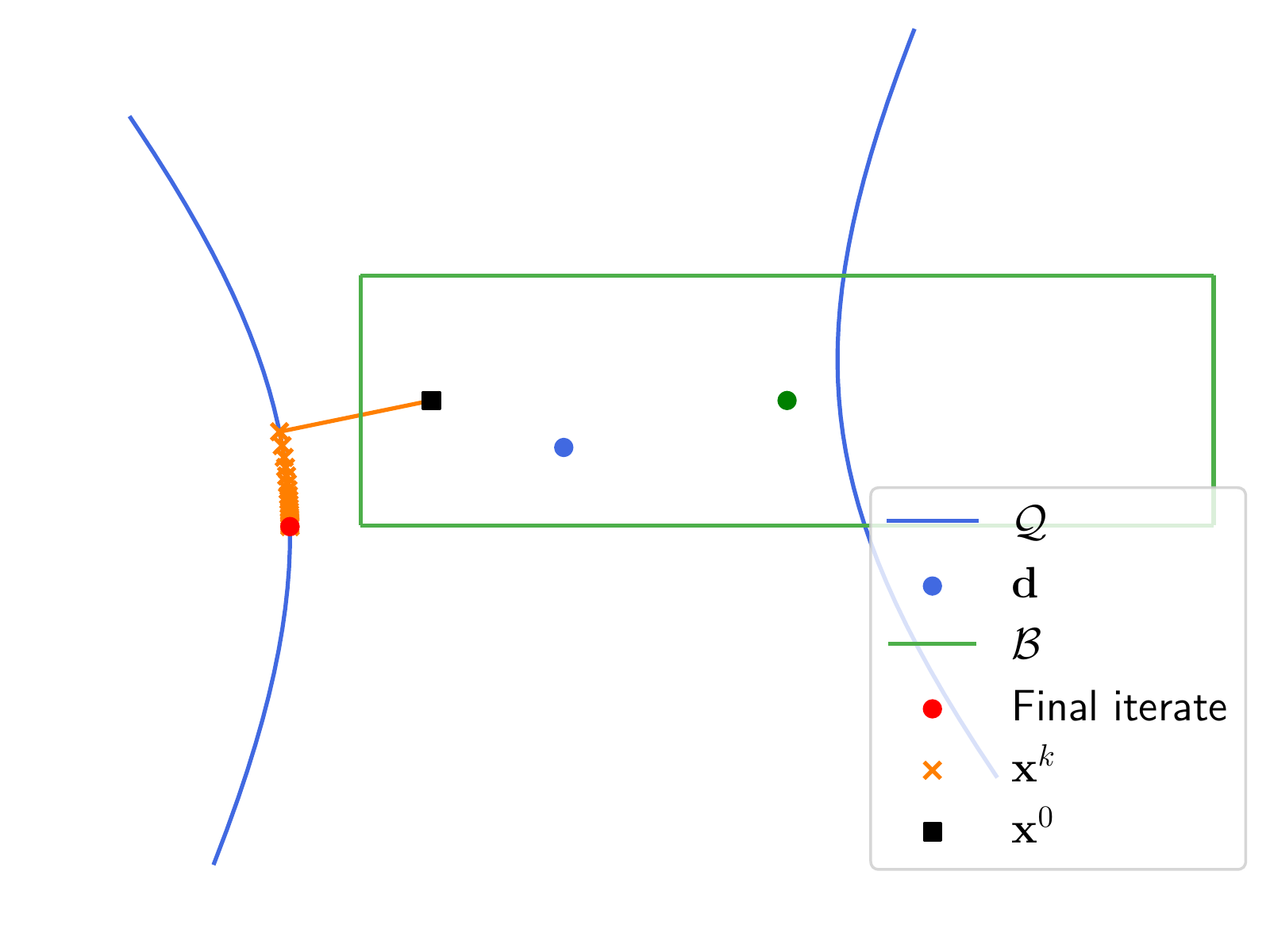}
		\caption{Modified Douglas-Rachford for the feasibility problem (DR-F).}
		\label{subfig:DR-F_OK}
	\end{subfigure}
	\caption{Illustration of the behaviour of the DR splitting algorithms on the same (pathological) case of~\cref{fig:pathological_AP}. On this problem, DR does converge to a feasible point while DR-F does not. However, DR-F converges to a stationary point (a minimum) of~\cref{eq:obj_dr-feasibility}.}
	\label{fig:pathological_DR}
\end{figure}

\subsection{Comparison}%
\label{sub:comparison}

Table~\ref{tab:comparison} compares the different complexities and convergence results of all the methods. Methods using exact projection onto the quadric require the diagonalization of $\A$ as a precomputation step, which typically costs $\mathcal{O}(n^3)$ flops.

\begin{table}
	\centering
	\caption{Comparison of the considered splitting methods.}
	\label{tab:comparison}
	{\footnotesize%
	\hspace*{-3cm}
	\begin{tabular}{lrrrrr}
		\toprule
					& APE				& APC				& APG			& DR				& DR-F 	\\ \midrule
		Complexity		& $\mathcal{O}(n^3+kn^2)$	& $\mathcal{O}(kn^2)$		& $\mathcal{O}(kn^2)$	& $\mathcal{O}(n^3+kn^2)$	& $\mathcal{O}(n^3+kn^2)$	\\
		Convergence guarantees	& Locally to a feasible  	& None 				& None			& None\footnote{There are, however, proofs of the convergence of DR in some nonconvex applications, see, \eg, the discussion in~\cite[Section~4]{aragon_douglasrachford_2020}.}				& Locally to a stationary\\
					& point of~\cref{eq:P}: \cref{prop:convergence_APE}	& 				&			&				& point of~\cref{eq:obj_dr-feasibility}: \cref{prop:convergence_DR-F}
	\end{tabular}
	} 
\end{table}

\subsection{Extensions and applications}%
\label{sub:Extensions_and_further_researches}

We can extend the splitting methods to a polytope, $\mathcal{P}$, and a Cartesian product of $m$ quadrics $\Qtot :=\bigtimes_{i=0}^m \Q_i$, and solve

\begin{align*} \label{eq:P_extension}
	\min_{\bm x \in \Rn} 	& \norm{\bm x - \bm x^0}_2^2 \numberthis \\ 
	\text{s.t. } 	& \bm x \in \mathcal{P} , \\
			& \bm x \in \Qtot .
\end{align*}

The extension of all methods described in \cref{tab:comparison} is direct: instead of computing the projection on $\mathcal{B}$---now $\mathcal{P}$---analytically, we have to resort to a QP solver. And, similarly to the retraction from~\cite{vh22}, the (quasi-)projection is obtained by working independently on each quadric $\Q_i$:

\begin{equation}
	\x^* \in \argmin_{\x \in \Qtot} \norm{\x - \x^0}_2 \Leftrightarrow \x^*_i \in \argmin_{\x_i \in \Q_i} \norm{\x_i - \x^0_i}_2 \forall i = 1 \hdots m .
\end{equation}

For example, in the practical case from~\cite{vh22}, the paper focuses on the dynamic economic dispatch problem which aims at the optimal allocation of power production among generating units at each timestep, \eg, each hour of a day. The modelling of the power losses makes the feasible set of each independent (static) economic dispatch a quadric and certain operational constraints, namely the ramping constraints, couple consecutive time steps. Hence, the full feasible set $\Omega$ is a polytope $\mathcal{P} \subseteq \mathbb{R}^{nT} $ that accounts for the power ranges (box) and ramping constraints, and of a Cartesian product of $T$ different quadrics $\Q_i \subseteq \Rn$ that model the balance constraint, \ie, that power production matches demand. The projection of a point onto $\Omega = \mathcal{P} \bigcap \Qtot$ can then be obtained using the methods described in the present paper. 

Moreover, the point that has to be projected in~\cite{vh22} is obtained as the solution of a surrogate problem defined on a relaxed set, see~\cite{vh19,vh20,vh22} for more details. And because this relaxation is close to the feasible set, the point that has to be projected is inside the box and near the quadric. This is the reason for the favorable performance of APC which is reported in~\cref{sub:Alternate_projection_versus_Gurobi}.

\subsection{Implementation details}%
\label{sub:Implementation_details}

To address the convergence issues identified in~\cref{fig:pathological_AP,fig:pathological_DR}, we add a restart mechanism  whenever this situation arises. Such situations are easily detected: the alternating method will loop between two points, and the DR or DR-F will simply converge to an infeasible point.
These problems mostly appear in the hyperboloid case, and typically occur when the method is trapped on the wrong sheet of the hyperboloid. To mitigate this issue, when detected, we use the geometric construction from the centre-based quasi-projection (\cref{alg:quasi-projection-quadric} with $\bm \xi = \x^0 - \bm d$), and select the \emph{largest} $\beta$. This is equivalent to transforming $\x^k \in \Q,  \x^k\notin \mathcal{B}$ into $-\x^k =:\x^{k+1} \in \Q$, and continuing the method from $\x^{k+1}$. Alternatively, it is also possible to consider $\x^{k+1}$ such that at least one---instead of all---of its components is the opposite of $\x^k$. If, on the other hand, $\x^k \in \mathcal{B}$, then we work analogously with respect to the centre of the box.

Such a restart mechanism is not a guarantee of convergence: the method can then be trapped into another region, or even come back to the exact same region. But in the few instances ($\approx$ once every 10000 trials) where the presented algorithms experience convergence issues, the restart results in successful convergence.

\section{Numerical experiments}%
\label{sec:Numerical_experiments}

This section is devoted to the benchmarking of the methods developed in~\cref{sec:Splitting_methods}.

\cref{sub:Douglas-Rachford_Alternate_projections_and_IPOPT} tests the five presented methods (APE, APC, APG, DR, DR-F) as well as IPOPT. IPOPT is an interesting method to benchmark against, as it is a natural candidate for solving~\cref{eq:P}. Note that IPOPT is an open-source solver that uses an embedded linear solver. The performance of IPOPT can be enhanced through the use of a dedicated commercial linear solver. In this work, we use Pardiso~\cite{pardiso-7.2a}.

We solve for small scale (\cref{fig:ellipsoid_100,fig:hyperboloid_100}) and larger scale (\cref{fig:ellipsoid_1000,fig:hyperboloid_1000}) instances of~\cref{eq:P}. For each considered dimension, $n$, we run 100 randomly independent trials in order to smooth the effect of the random selection of the problem parameters. In particular, each independent trial consists of a unique (randomly generated) set of parameters $\A, \Ab, \Ac$ and $\x^0 \in \Rn$.

The ellipsoid case is tested in~\cref{ssub:Ellipsoid_experiments} and the hyperboloid case in~\cref{ssub:Hyperboloid_experiments}.

The problems are generated as follows: we first create a quadric with $\bm A \sim \mathcal{N}(\bm \mu = \mathbb{1}, \Sigma = I)$, $\A = \frac{\bm A + \Tr{\bm A}}{2}$, $\Ab$ defined with $b_i \sim \mathcal{N}(\mu = 0, \sigma = 1)$ and $c \sim \mathcal{N}(\mu = -1, \sigma = 1)$. For the ellipsoidal case, we shift $\A$ in order to ensure that $\A \succ \bm 0$. Then, we find one feasible point and construct the box $\mathcal{B}$ around it; this allows us to ensure that the intersection of $\mathcal{B}$ and $\Q$ is nonempty. Note that this feasible point is not necessarily the centre of the box.

Then, in~\cref{sub:Alternate_projection_versus_Gurobi}, we perform two experiments for comparing APC to Gurobi for a very specific problem structure, which is a problem stemming from~\cite{vh22}, and the initial goal of the present research. The same remarks with respect to the 100 randomly generated data also apply here. 

For this second experiment, we use Gurobi as a benchmark because i) it may also be a natural method for solving~\cref{eq:P}\footnote{Since version 9.0, Gurobi supports nonconvex QCQP optimization, and Gurobi is employed widely in the power systems optimization community.} and ii) it provides lower bounds, which allow us to assess how close the returned solutions are to the global optimum. We benchmark it against APC because, even if it is the worst-performing methods among the five that are presented in~\cref{sub:Douglas-Rachford_Alternate_projections_and_IPOPT}, it still outperforms Gurobi. This behaviour is explained by the relative position of the starting point with respect to the feasible set from the problem of~\cref{sub:Alternate_projection_versus_Gurobi}: this point is inside the box and close to the quadric.

Note that, in all experiments, whenever an algorithm terminates with a timeout and returns an infeasible point, the associated objective is meaningless. In order to avoid distorting our reported results, we omit these instance in the recorded objectives; but we count the number of timeouts and record the deviation.

The deviation is computed as
\begin{equation}
	\label{eq:deviation}
	\textrm{deviation} = \abs{\Tr{\x} \A \x + \Tr{\Ab} \x + \Ac},
\end{equation}
and is an intuitive measure of how far an infeasible point is to the feasible set. The prescribed tolerance for the deviation is $10^{-6}$. This deviation does not account for the box. This is not an issue here, because none of the tests considered in the numerical experiments terminates outside of the box.

\subsection{Douglas-Rachford, Alternating Projections and IPOPT}%
\label{sub:Douglas-Rachford_Alternate_projections_and_IPOPT}

Two different settings are considered here. In~\cref{ssub:Ellipsoid_experiments}, the matrix $\A$ is chosen such that $\A \succ \bm 0$, \ie, the quadric is an \emph{ellipsoid}. This means that the quasi-projection with $\bm \xi = \x^0 - \bm d$ is well-defined: situations depicted in~\cref{fig:hyperbole_KO} cannot occur. In~\cref{ssub:Hyperboloid_experiments}, we consider the case of \emph{hyperboloids}, \ie, $\A$ is nonsingular but indefinite.

From these two experiments, it appears that both DR-F and APE are the methods that find the best solution in terms of objective. However, if the execution time is taken into account, APG reaches an objective close to the one of DR-F and APE in a significantly lower run time. APG should therefore be considered, \eg, if the eigenvector decomposition is too expensive to compute. APC works particularly well in the ellipsoidal case, but performs worse in the hyperboloidal case. IPOPT is clearly the slowest method. It achieves good solution objectives in the ellipsoidal case, but gives poorer results in the hyperboloidal case.

\subsubsection{Ellipsoid experiments}%
\label{ssub:Ellipsoid_experiments}

In these two experiments, we run small and large-scale ellipsoidal problems. The box $\mathcal{B}$ is small with respect to the quadric and the starting points $\x^0$ are  uniformly distributed inside the box.

\begin{figure}
	\begin{subfigure}[b]{.475\textwidth}
		\begin{center}
			\includegraphics[width=\textwidth]{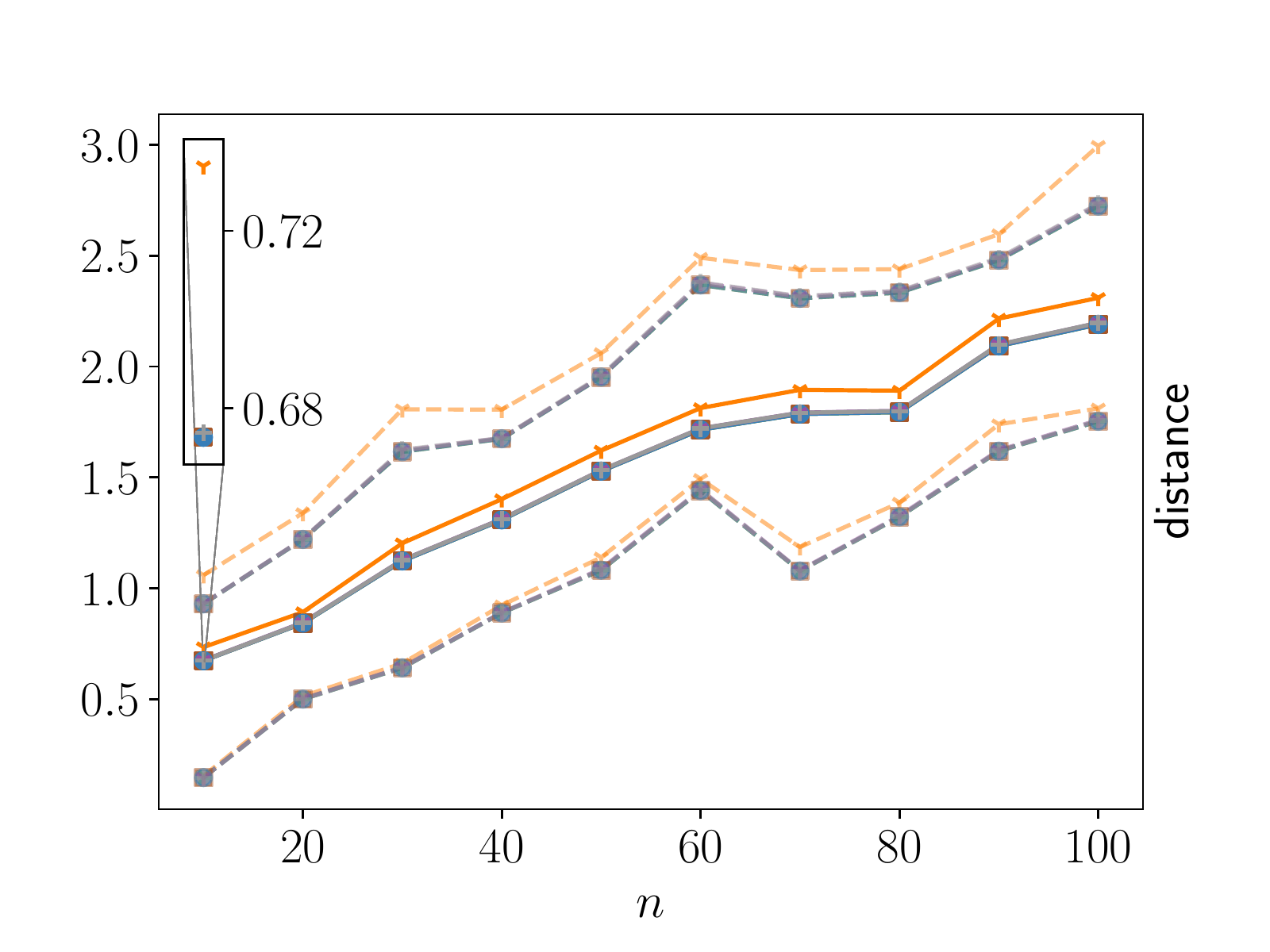}	
		\end{center}
		\caption{Distance.}
		\label{subfig:distance_ellipsoid_100}
	\end{subfigure}
	\hfill
	\begin{subfigure}[b]{.475\textwidth}
		\begin{center}
			\includegraphics[width=\textwidth]{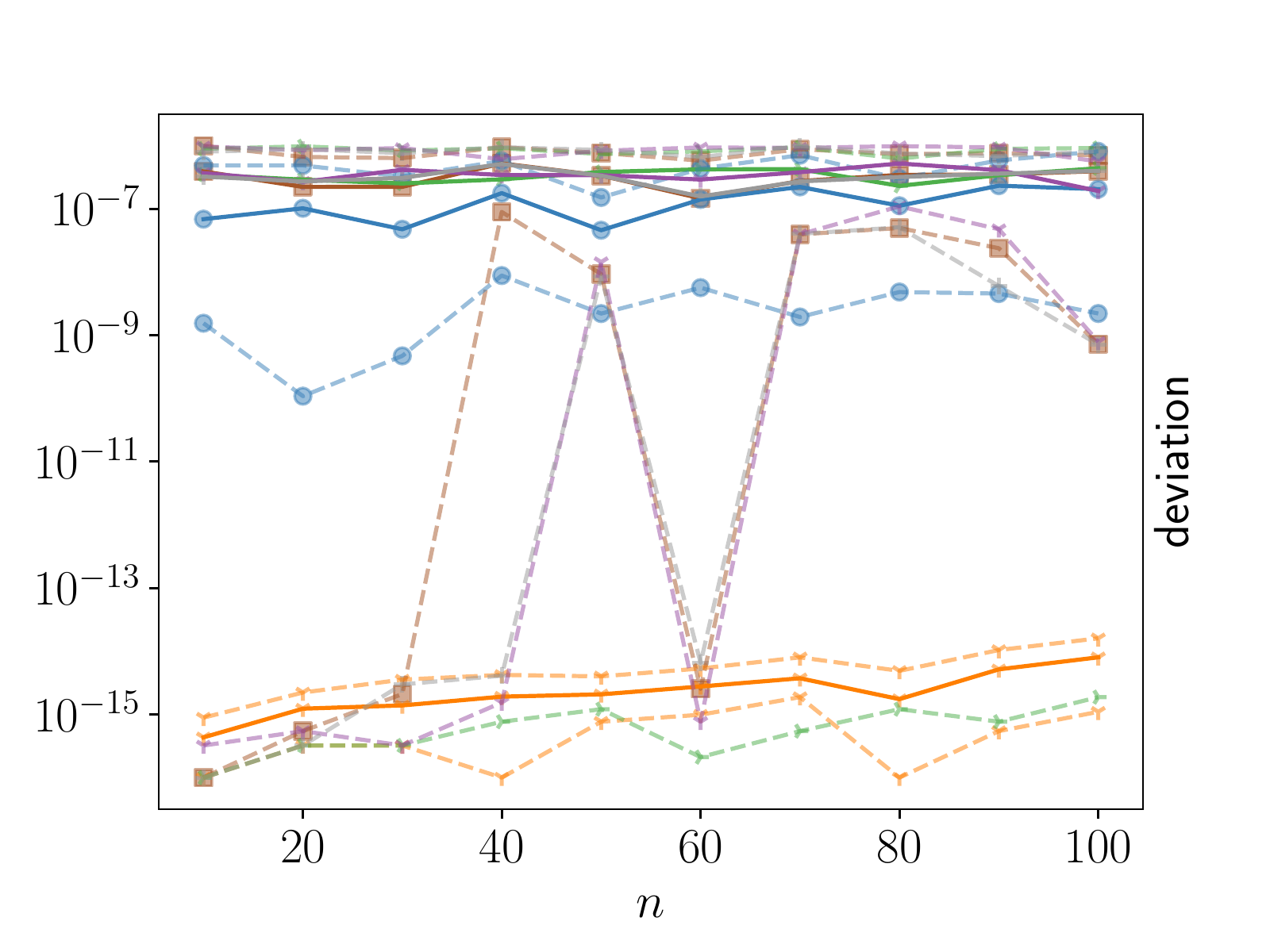}	
		\end{center}
		\caption{Deviation.}
		\label{subfig:deviation_ellipsoid_100}
	\end{subfigure}
	\begin{subfigure}[b]{.475\textwidth}
		\begin{center}
			\includegraphics[width=\textwidth]{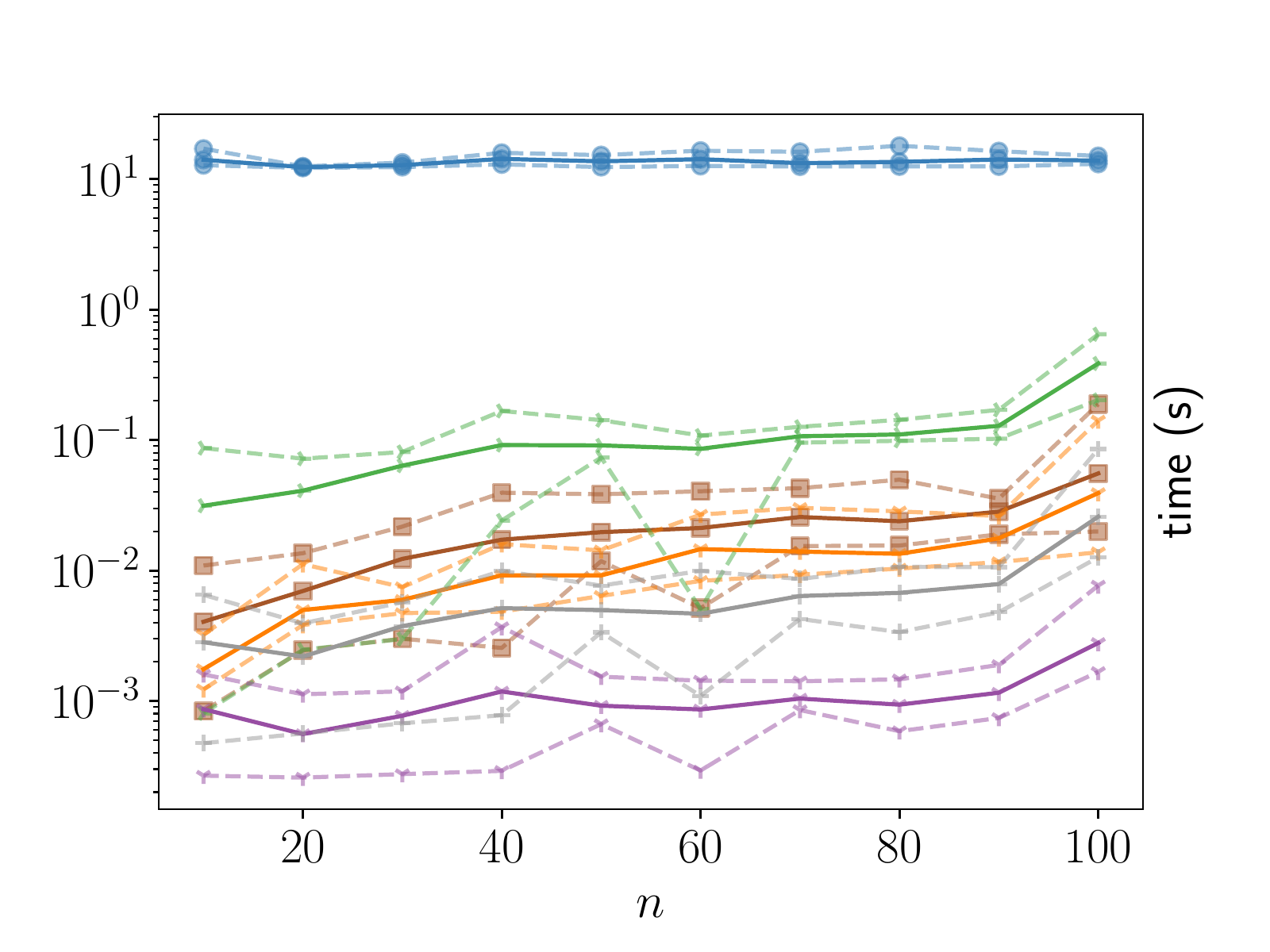}	
		\end{center}
		\caption{Execution time in seconds.}
		\label{subfig:total-time_ellipsoid_100}
	\end{subfigure}
	\hfill
	\begin{subfigure}[b]{.475\textwidth}
		\begin{center}
			\includegraphics[width=\textwidth]{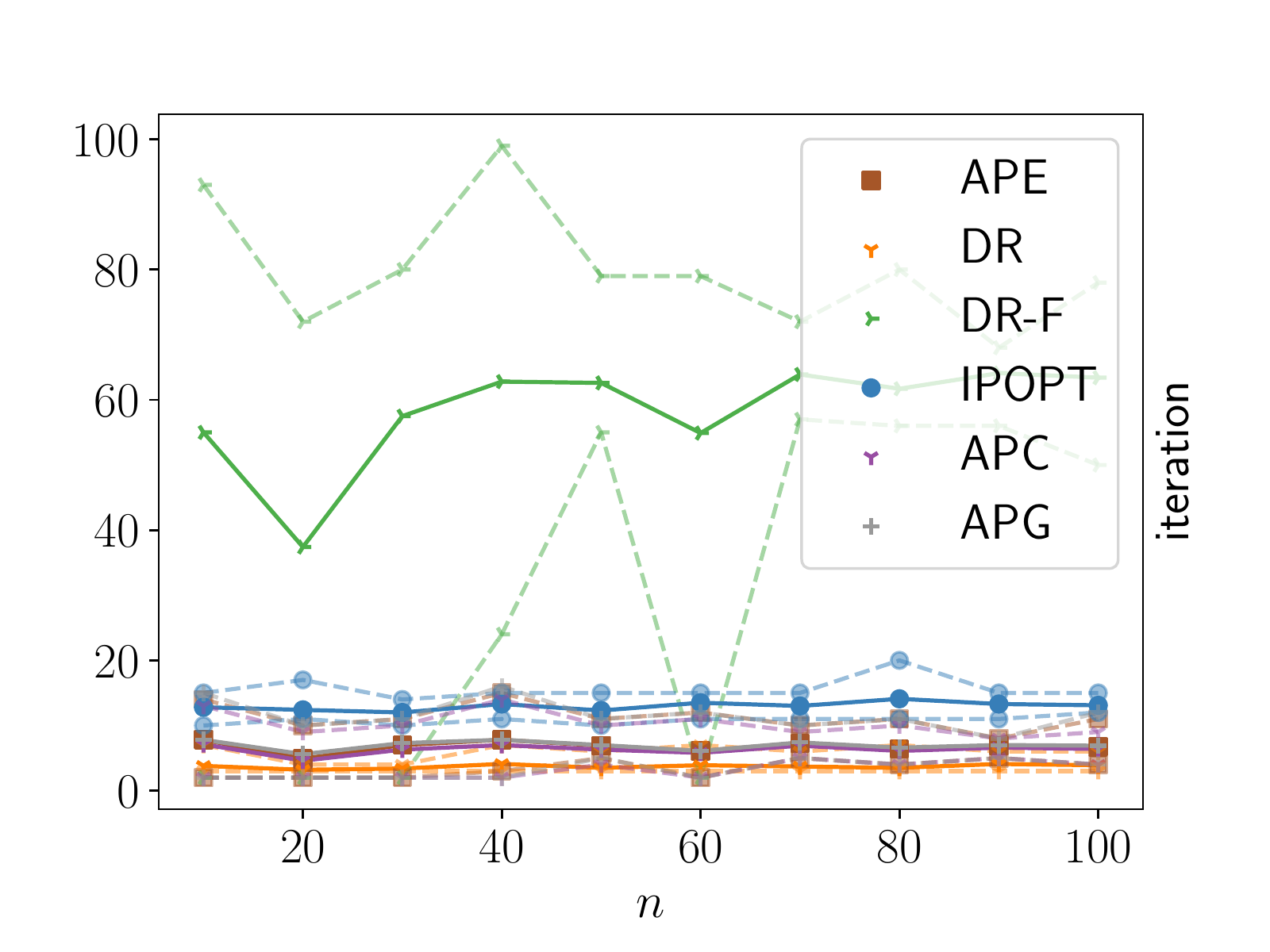}	
		\end{center}
		\caption{Number of iterations.}
		\label{subfig:iteration_ellipsoid_100}
	\end{subfigure}
	\caption{Comparison of the different methods developed in~\cref{sec:Splitting_methods}: Douglas-Rachford splitting (DR) and its modified counterpart (DR-F), alternating projections using the exact projection (APE) and the alternating projections using the quasi-projections (centre-based APC and gradient-based APG). IPOPT is used as a benchmark with standard settings and with the underlying linear solver Pardiso. Ten dimensions $n$ are considered and, for each $n$, 100 independent trials with $\A \succ \bm{0}$ are run. The top (bottom) dashed lines represent the max (min) value of the 100 trials, and the continuous line is the sample mean. The frame in the upper left of the upper left panel is a magnification around $n=10$.}
	\label{fig:ellipsoid_100}
\end{figure}

For small-scale ellipsoidal problems ($n\leq 100$,~\cref{fig:ellipsoid_100}), we observe that all methods except DR reach the same objective: APE, DR-F and IPOPT obtain the same objective, APG is within 1\% and DR within several percent. We also observe that none of the methods exceeds the prescribed deviation accuracy of $10^{-6}$, and that IPOPT provides the most feasible points. 

The number of iterations required for each method remains more or less constant when the dimension increases. Considering the running time, APC is the fastest and ten times faster than APG, which is two time faster than DR. DR-F and APE require approximately the same amount of time, which is two times slower than DR. Finally, IPOPT requires much more time than all the other methods.

\begin{figure}
	\begin{subfigure}[b]{.475\textwidth}
		\begin{center}
			\includegraphics[width=\textwidth]{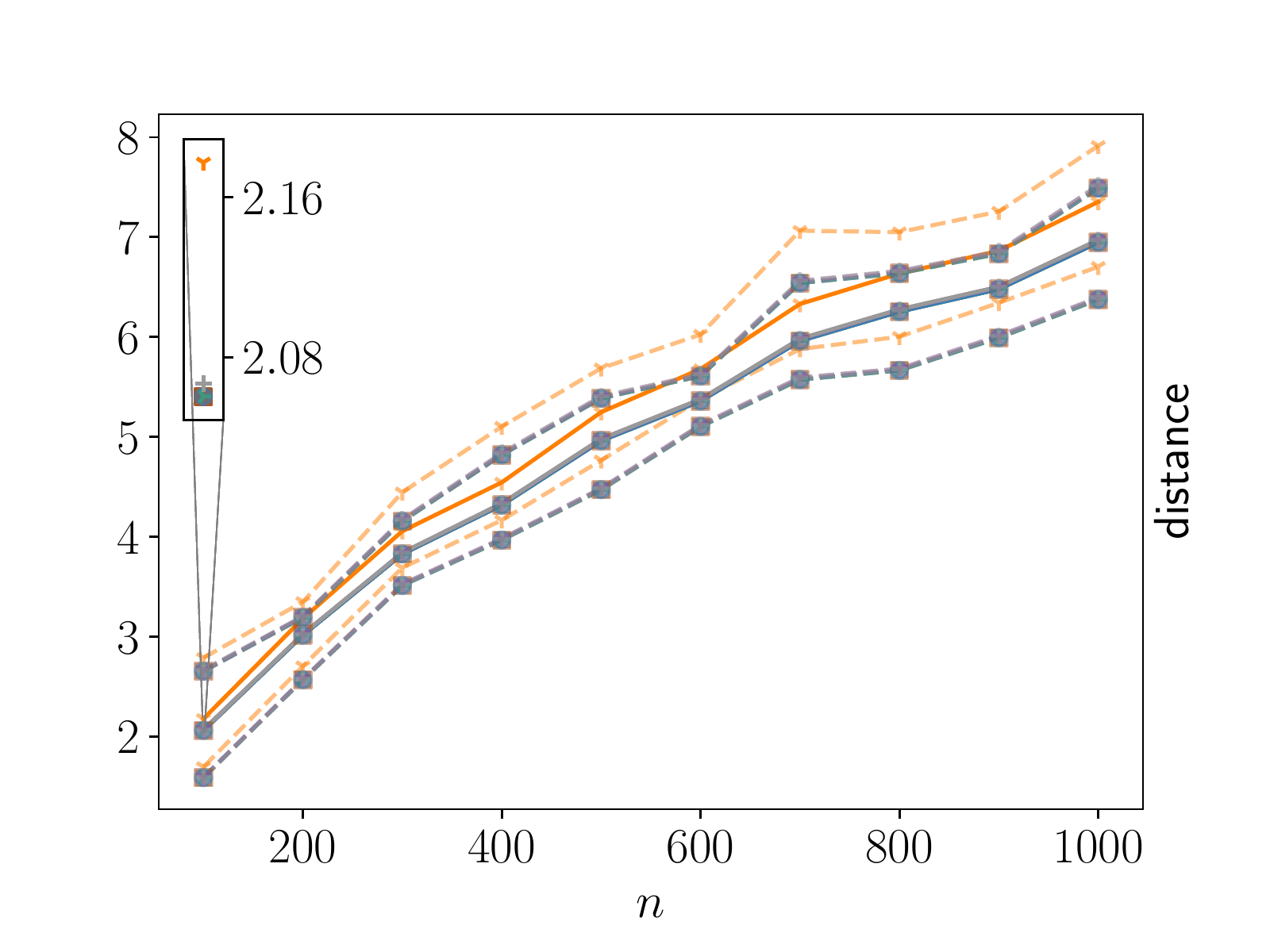}	
		\end{center}
		\caption{Distance.}
		\label{subfig:distance_ellipsoid_1000}
	\end{subfigure}
	\hfill
	\begin{subfigure}[b]{.475\textwidth}
		\begin{center}
			\includegraphics[width=\textwidth]{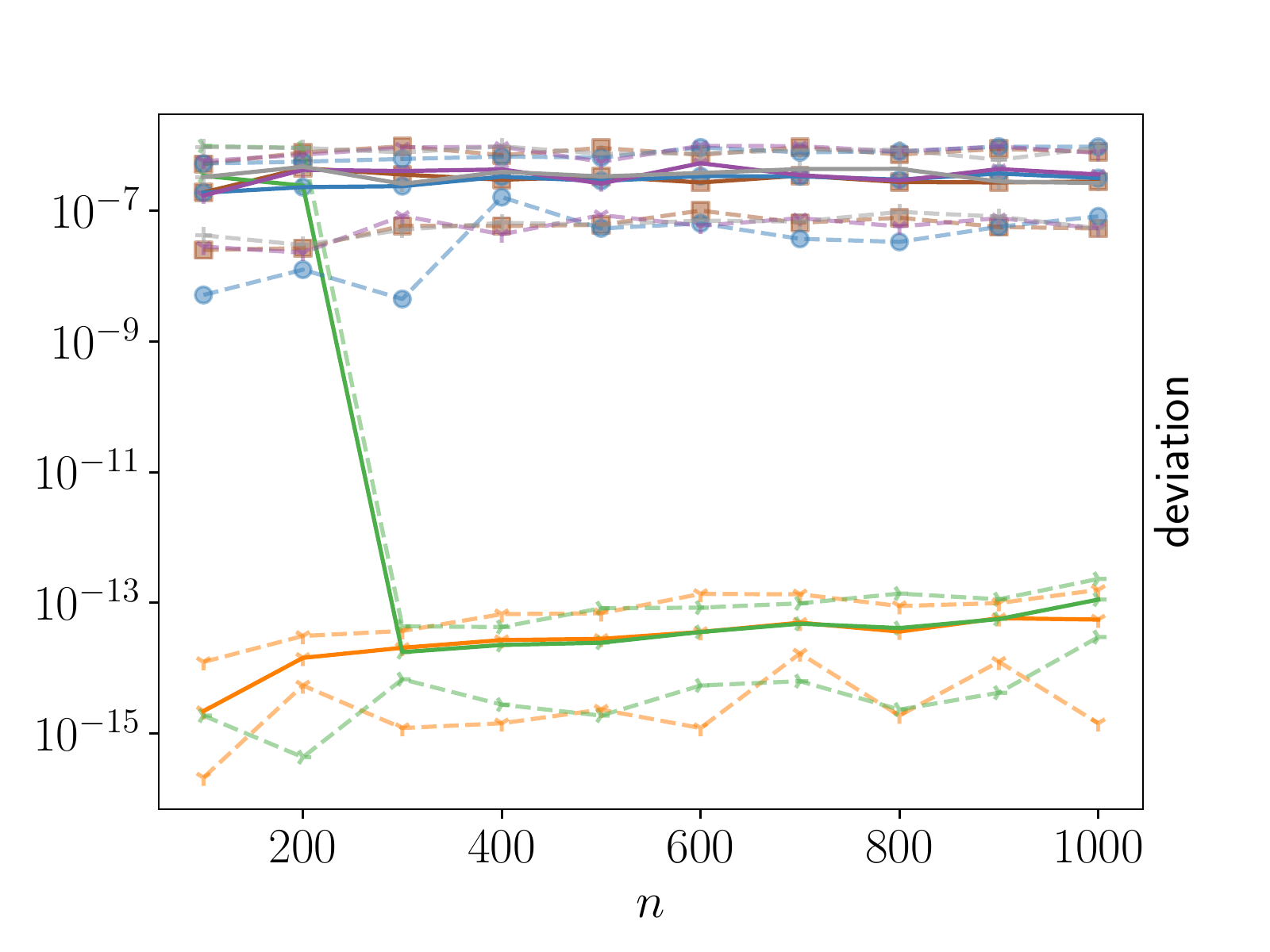}	
		\end{center}
		\caption{Deviation.}
		\label{subfig:deviation_ellipsoid_1000}
	\end{subfigure}
	\begin{subfigure}[b]{.475\textwidth}
		\begin{center}
			\includegraphics[width=\textwidth]{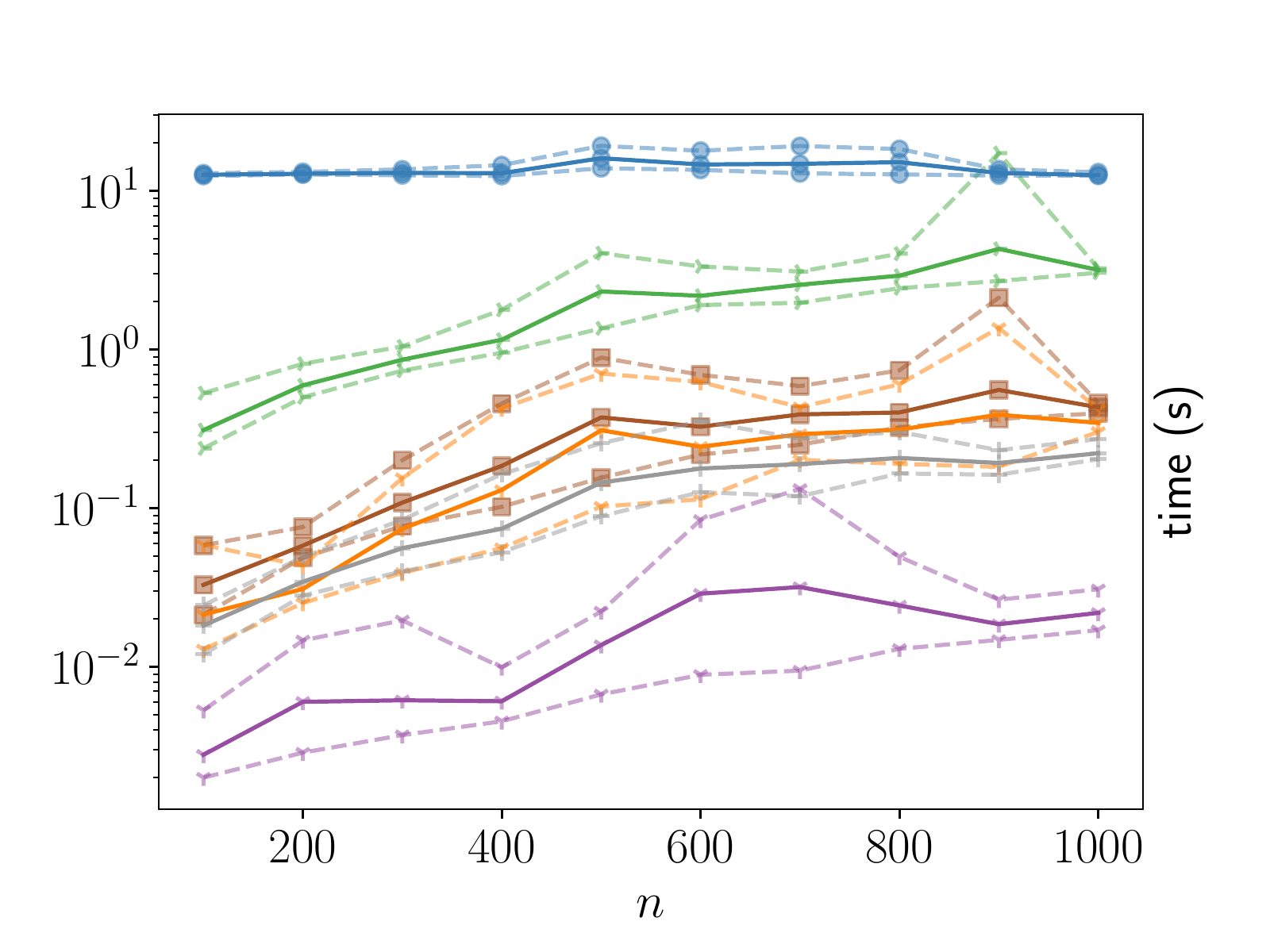}	
		\end{center}
		\caption{Execution time in seconds.}
		\label{subfig:total-time_ellipsoid_1000}
	\end{subfigure}
	\hfill
	\begin{subfigure}[b]{.475\textwidth}
		\begin{center}
			\includegraphics[width=\textwidth]{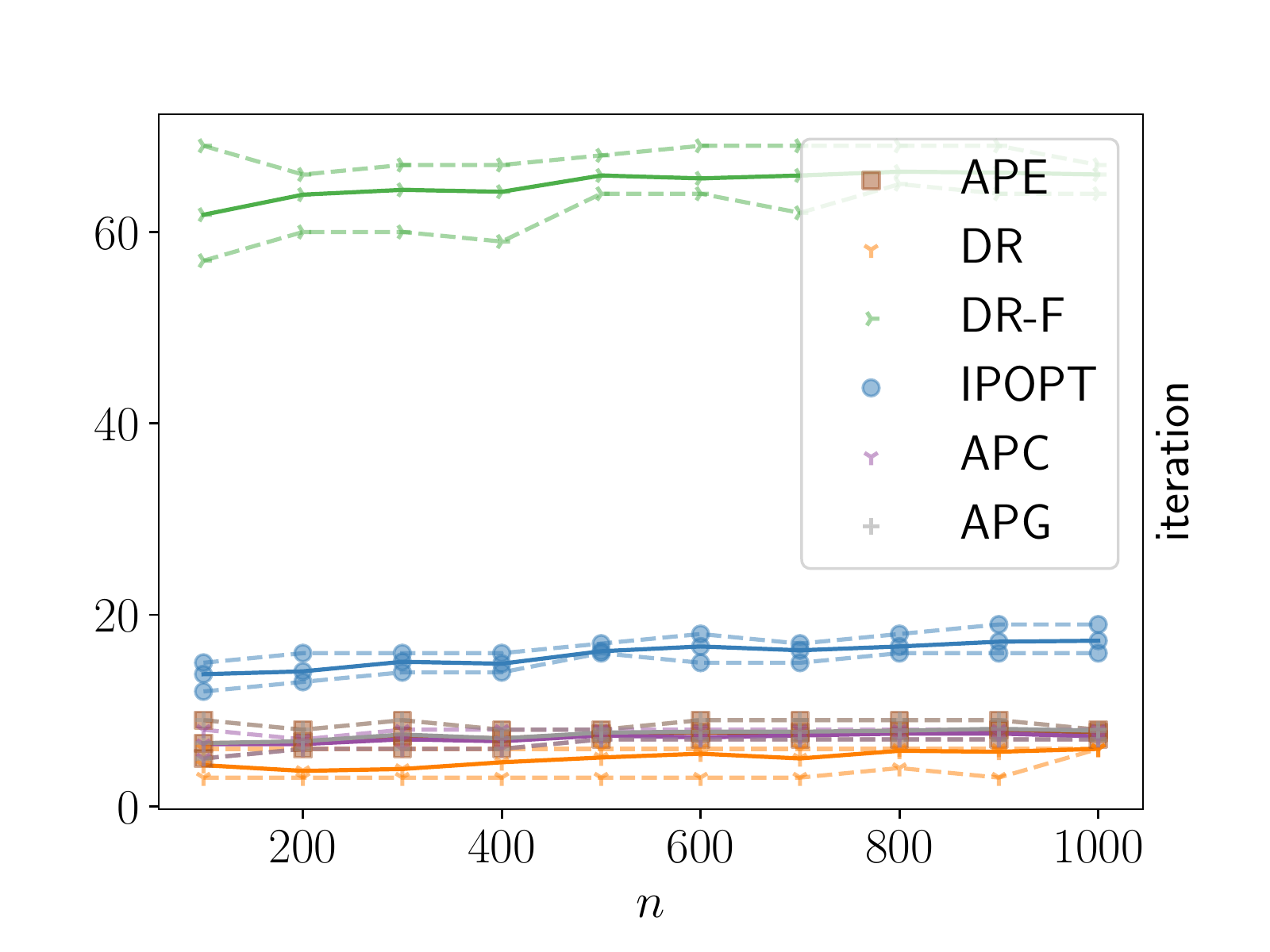}	
		\end{center}
		\caption{Number of iterations.}
		\label{subfig:iteration_ellipsoid_1000}
	\end{subfigure}
	\caption{Same as~\cref{fig:ellipsoid_100} for larger dimensions.}
	\label{fig:ellipsoid_1000}
\end{figure}

For large-scale ellipsoidal problems ($n\geq 100$,~\cref{fig:ellipsoid_1000}), the behaviour of the methods remains similar as the small-scale case. The distance increase with $n$ is simply due to the increase of $\norm{\x^* - \x^0}_2$ with $n$.

Remark that the execution time of IPOPT is remarkably stable, this is because creating the model already requires approximately 10 seconds, and this creation time does not increase much when the dimension increases. However, it should be noted that i) for much larger dimension $n>>1000$ the solving time of IPOPT increases significantly and ii) for such large dimension, it becomes crucial to use advanced linear algebra tools for, \eg, the eigenvalue decompositions and matrix products used in the methods developed here. Hence, the comparison against IPOPT when the latter relies on a dedicated linear algebra software (Pardiso) becomes less meaningful for too large $n$.


\subsubsection{Hyperboloid experiments}%
\label{ssub:Hyperboloid_experiments}

In these two experiments, we run small and large-scale hyperboloidal problems. The box $\mathcal{B}$ is large with respect to the quadric and the starting point $\x^0$ is uniformly distributed inside $\mathcal{B}$. 

\begin{figure}
	\begin{subfigure}[b]{.475\textwidth}
		\begin{center}
			\includegraphics[width=\textwidth]{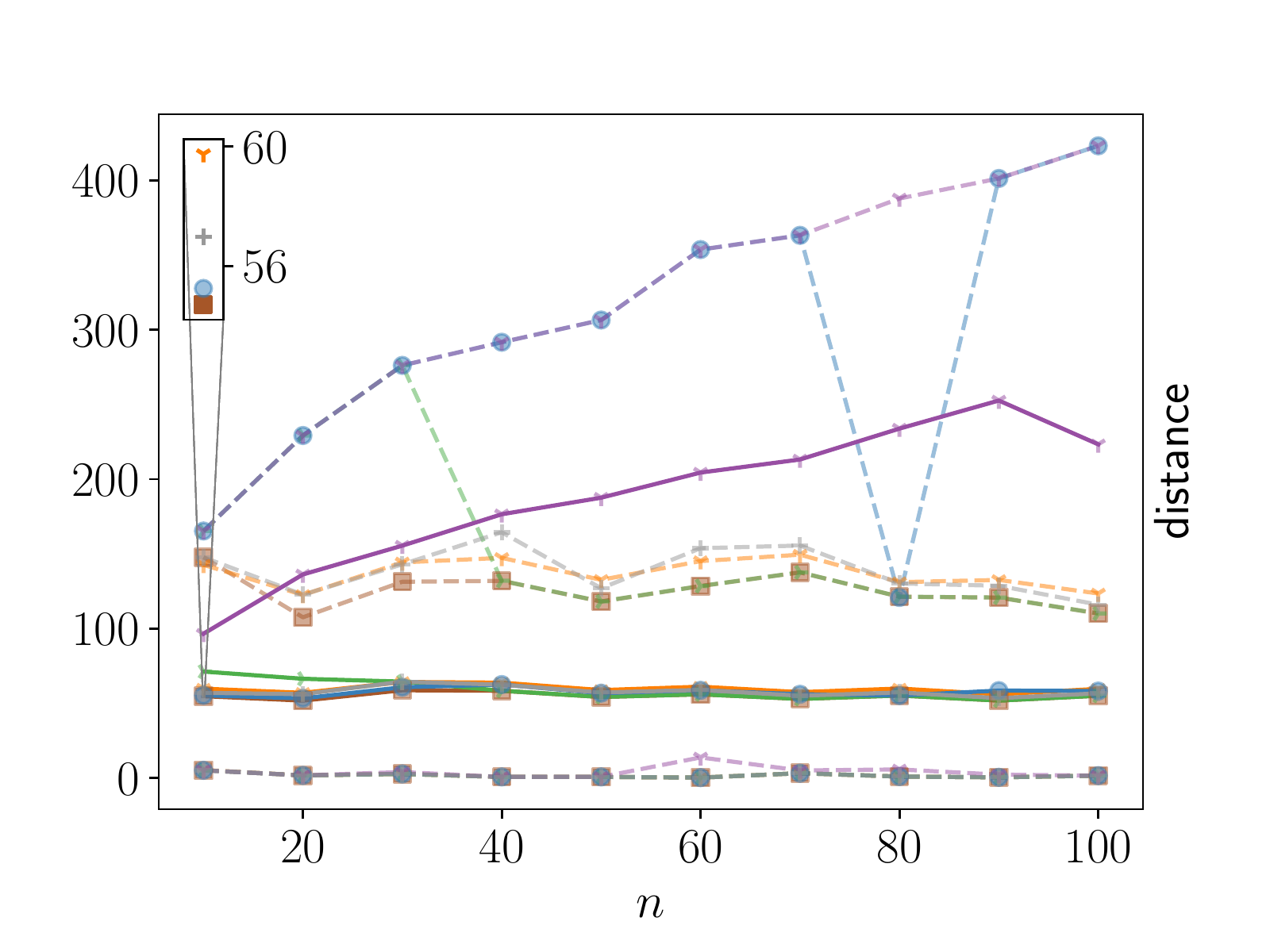}	
		\end{center}
		\caption{Distance.}
		\label{subfig:hyperboloid_100_distance}
	\end{subfigure}
	\hfill
	\begin{subfigure}[b]{.475\textwidth}
		\begin{center}
			\includegraphics[width=\textwidth]{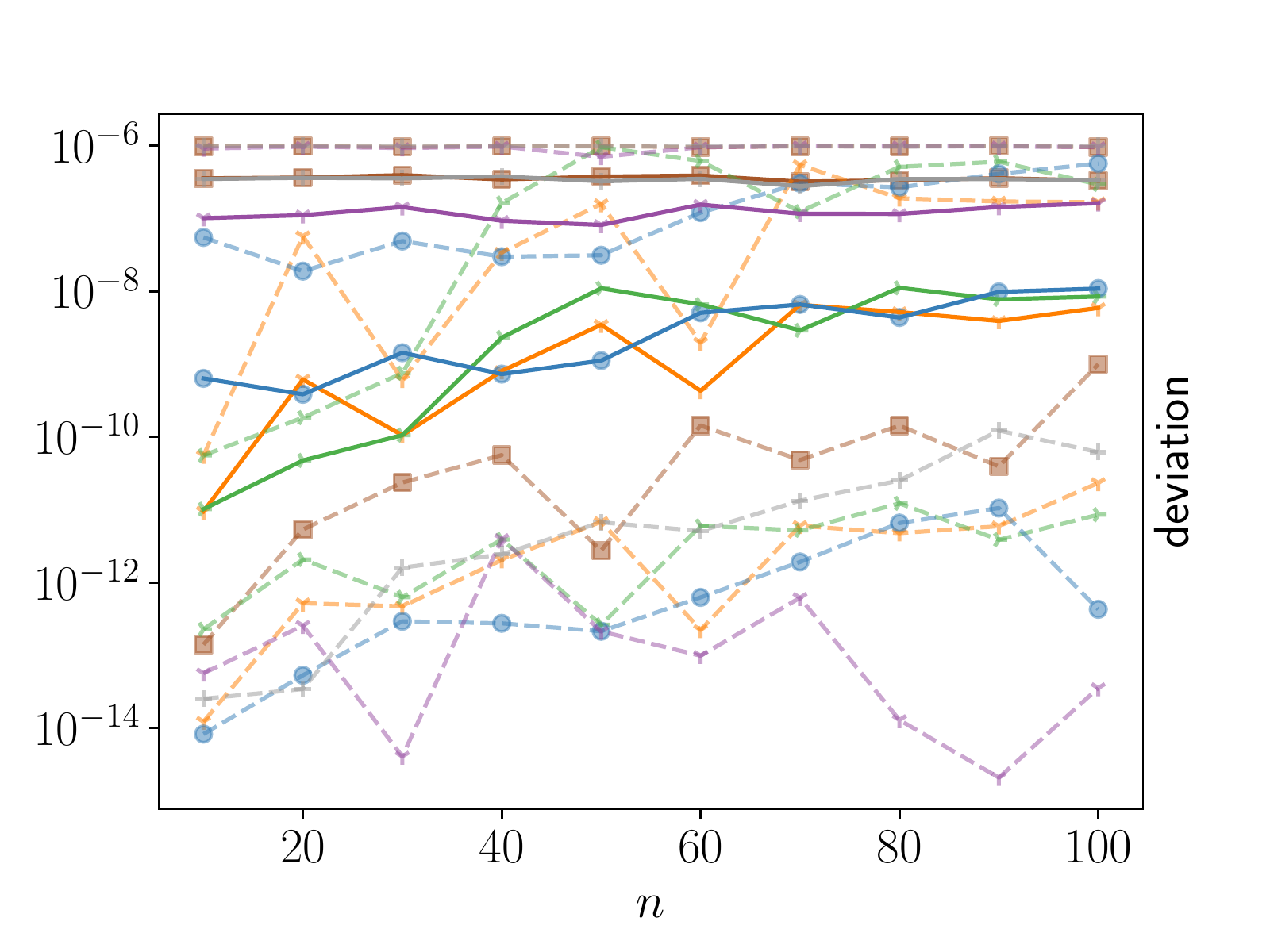}	
		\end{center}
		\caption{Deviation.}
		\label{subfig:hyperboloid_100_deviation}
	\end{subfigure}
	\begin{subfigure}[b]{.475\textwidth}
		\begin{center}
			\includegraphics[width=\textwidth]{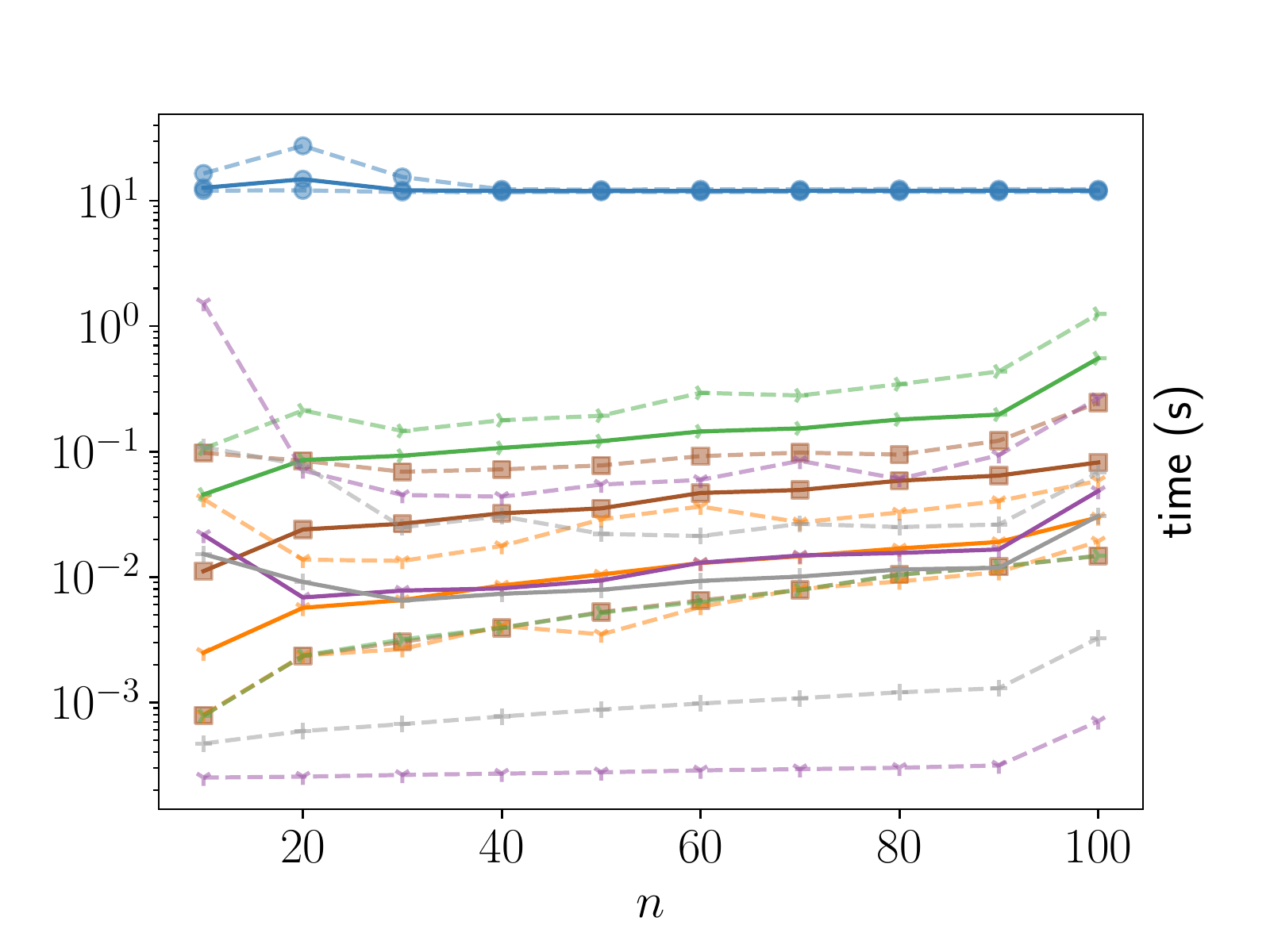}	
		\end{center}
		\caption{Execution time in seconds.}
		\label{subfig:hyperboloid_100_time}
	\end{subfigure}
	\hfill
	\begin{subfigure}[b]{.475\textwidth}
		\begin{center}
			\includegraphics[width=\textwidth]{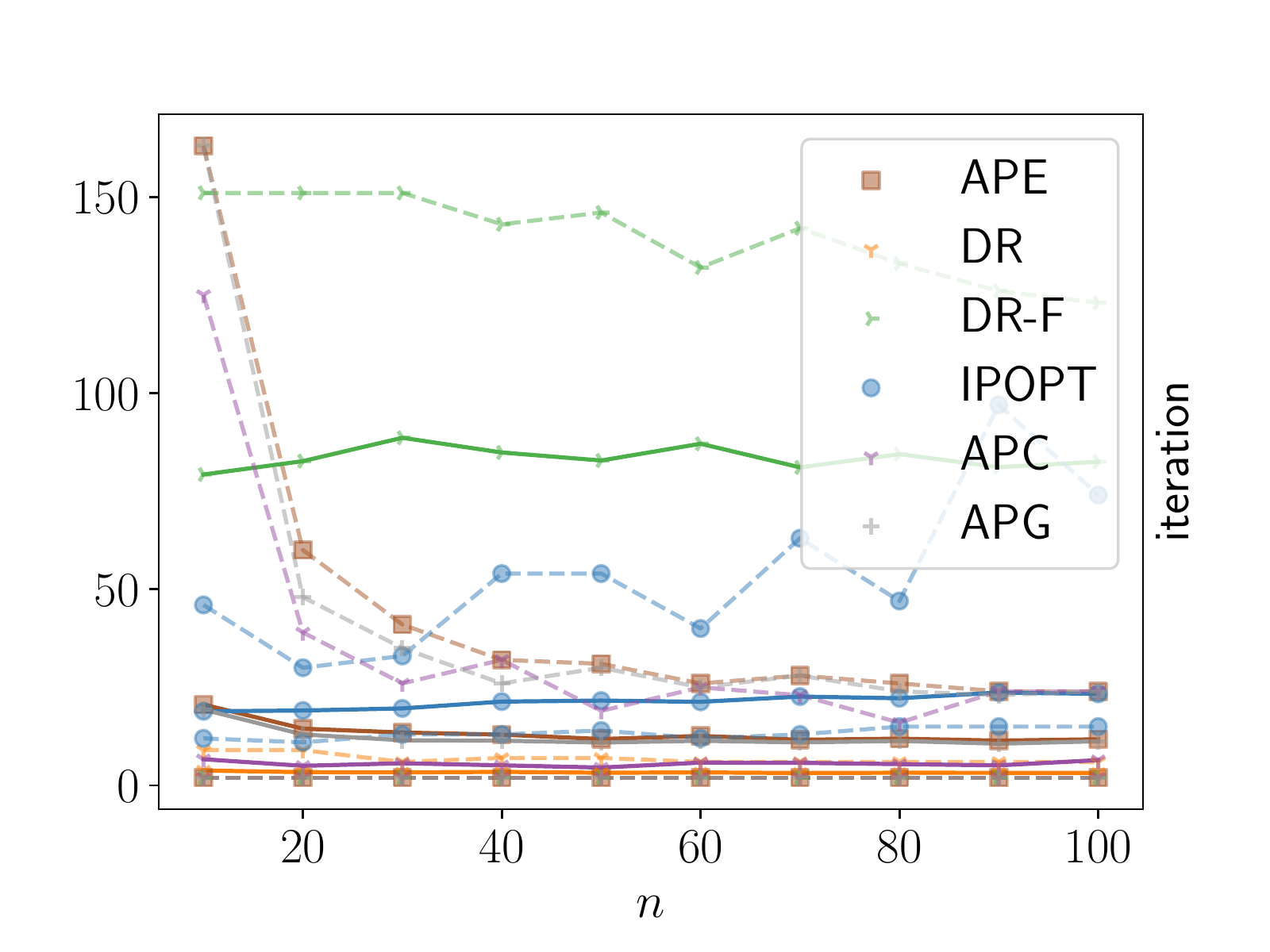}	
		\end{center}
		\caption{Number of iterations.}
		\label{subfig:hyperboloid_100_timeout}
	\end{subfigure}
	\caption{Same as~\cref{fig:ellipsoid_100} with $\bm \A \nsucc \bm 0$.}
	\label{fig:hyperboloid_100}
\end{figure}

For small-scale hyperboloidal problems ($n\leq 100$,~\cref{fig:hyperboloid_100}), we observe that the best objectives are obtained by APE. Then, within several \%, by DR-F, APG, DR and IPOPT. These four methods reach the same solution as APE most of the time, however they sometimes reach solutions that are far away from the best methods, see, \eg, the maximum curve (top dashed-lines in~\cref{subfig:hyperboloid_100_distance}) that is significantly above the maximum curves of APE. Finally, APC performs poorly in terms of objective values.
It is now APG which is the fastest method, despite its need of more iterations: the reason stems from the need of APC to resort to an exact projection whenever the situation depicted in~\cref{fig:hyperbole_KO} appears.

\begin{figure}
	\begin{subfigure}[b]{.475\textwidth}
		\begin{center}
			\includegraphics[width=\textwidth]{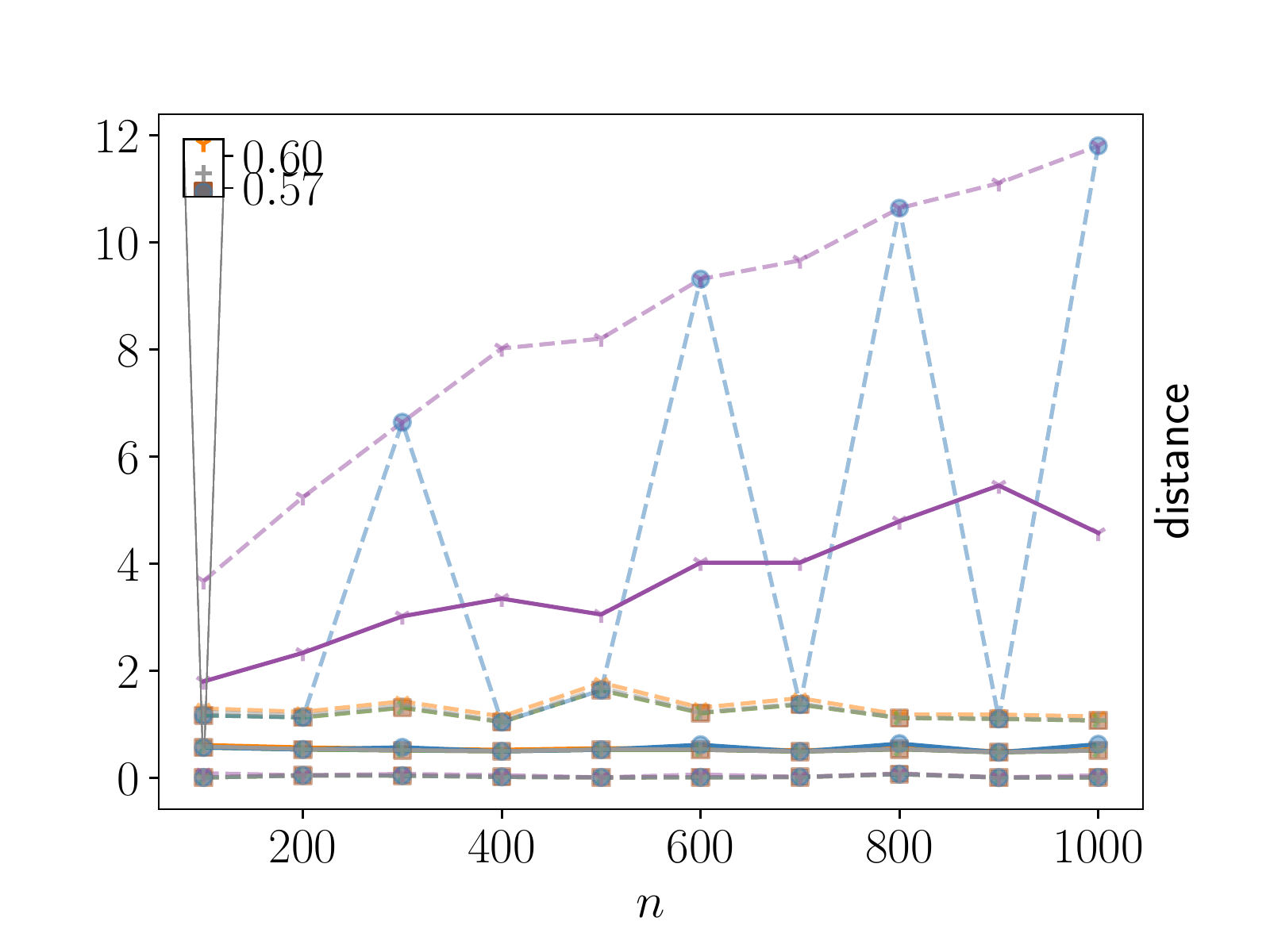}	
		\end{center}
		\caption{Distance.}
		\label{subfig:hyperboloid_1000_distance}
	\end{subfigure}
	\hfill
	\begin{subfigure}[b]{.475\textwidth}
		\begin{center}
			\includegraphics[width=\textwidth]{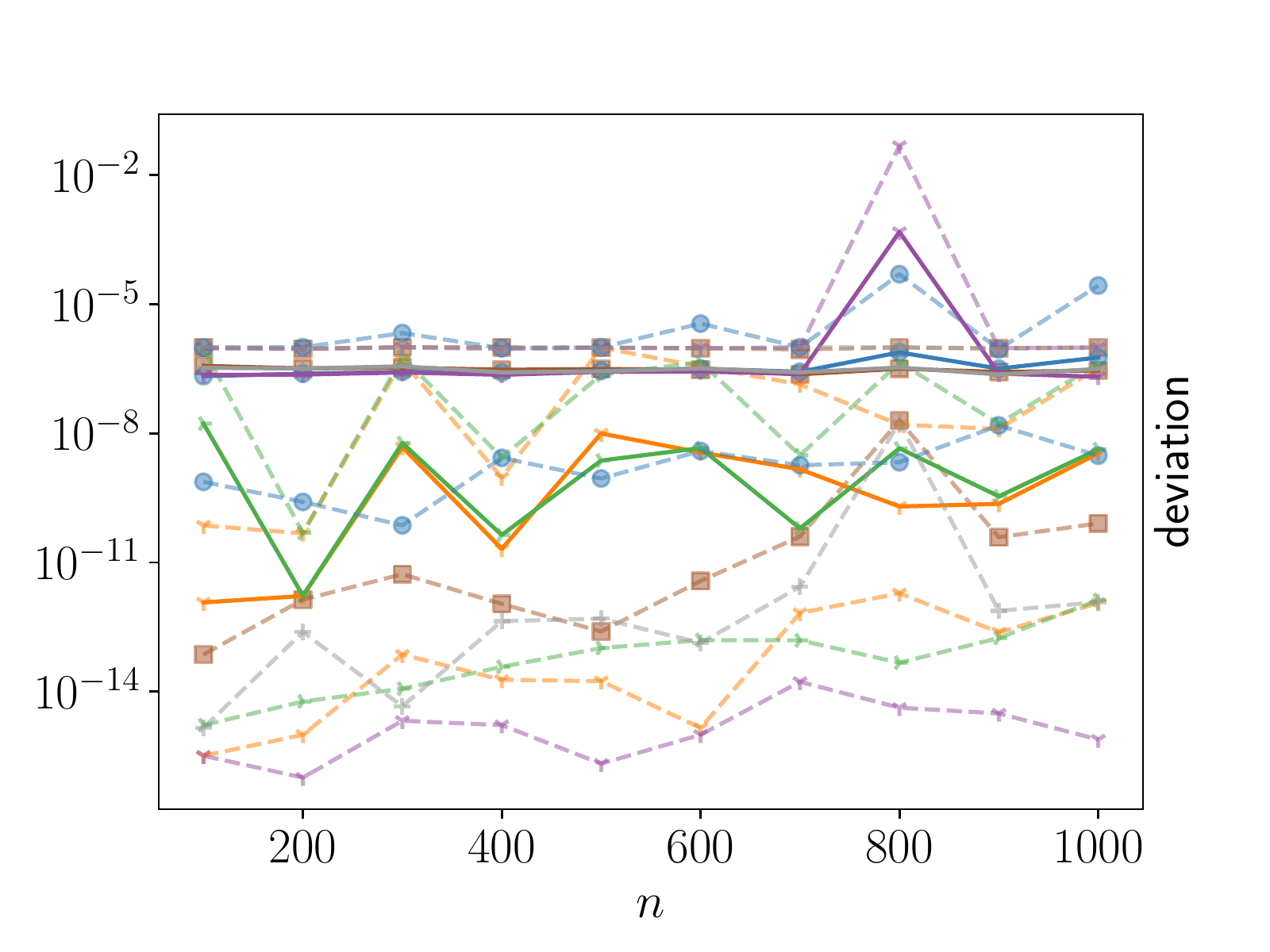}	
		\end{center}
		\caption{Deviation.}
		\label{subfig:hyperboloid_1000_deviation}
	\end{subfigure}
	\begin{subfigure}[b]{.475\textwidth}
		\begin{center}
			\includegraphics[width=\textwidth]{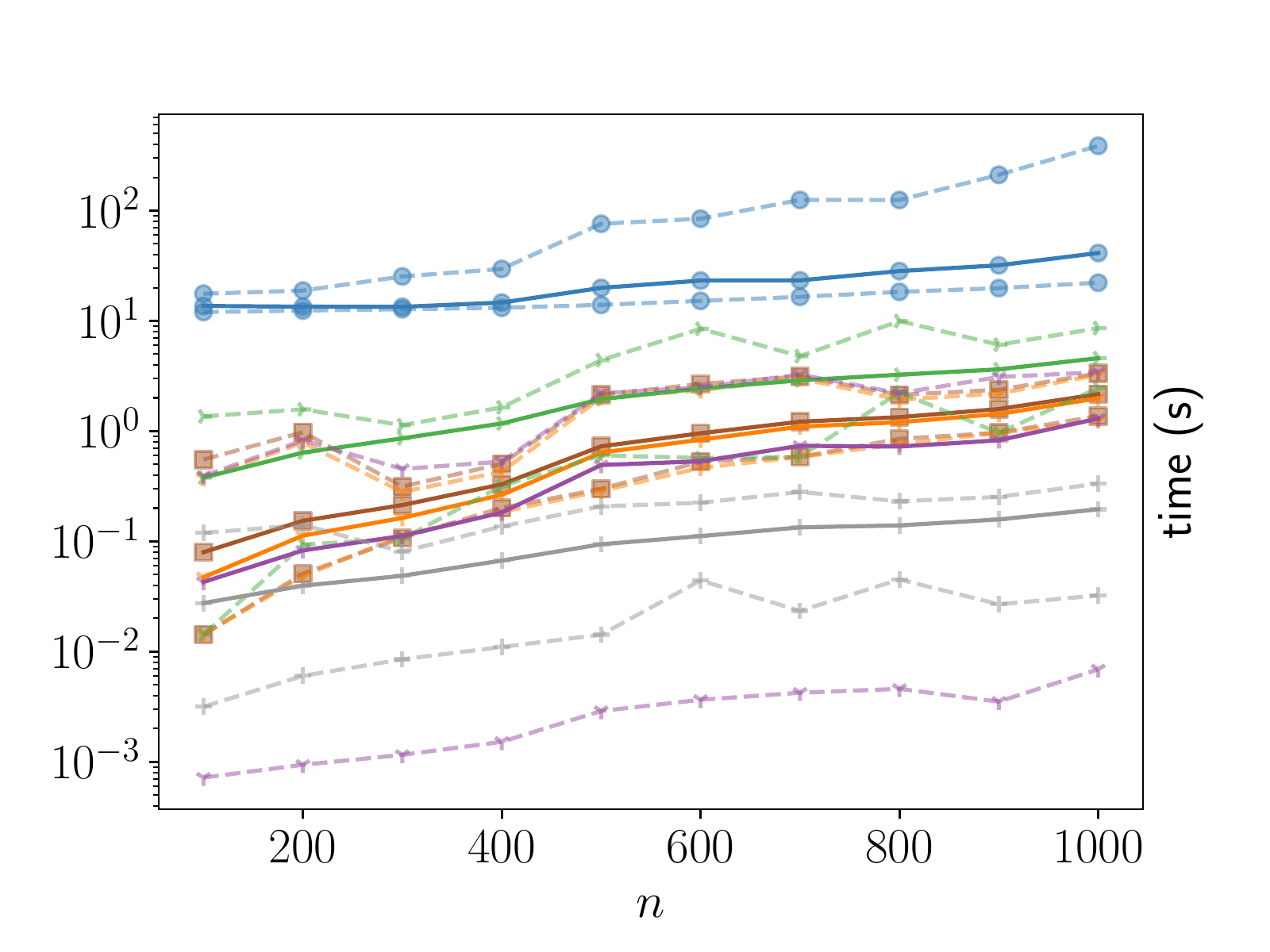}	
		\end{center}
		\caption{Execution time in seconds.}
		\label{subfig:hyperboloid_1000_time}
	\end{subfigure}
	\hfill
	\begin{subfigure}[b]{.475\textwidth}
		\begin{center}
			\includegraphics[width=\textwidth]{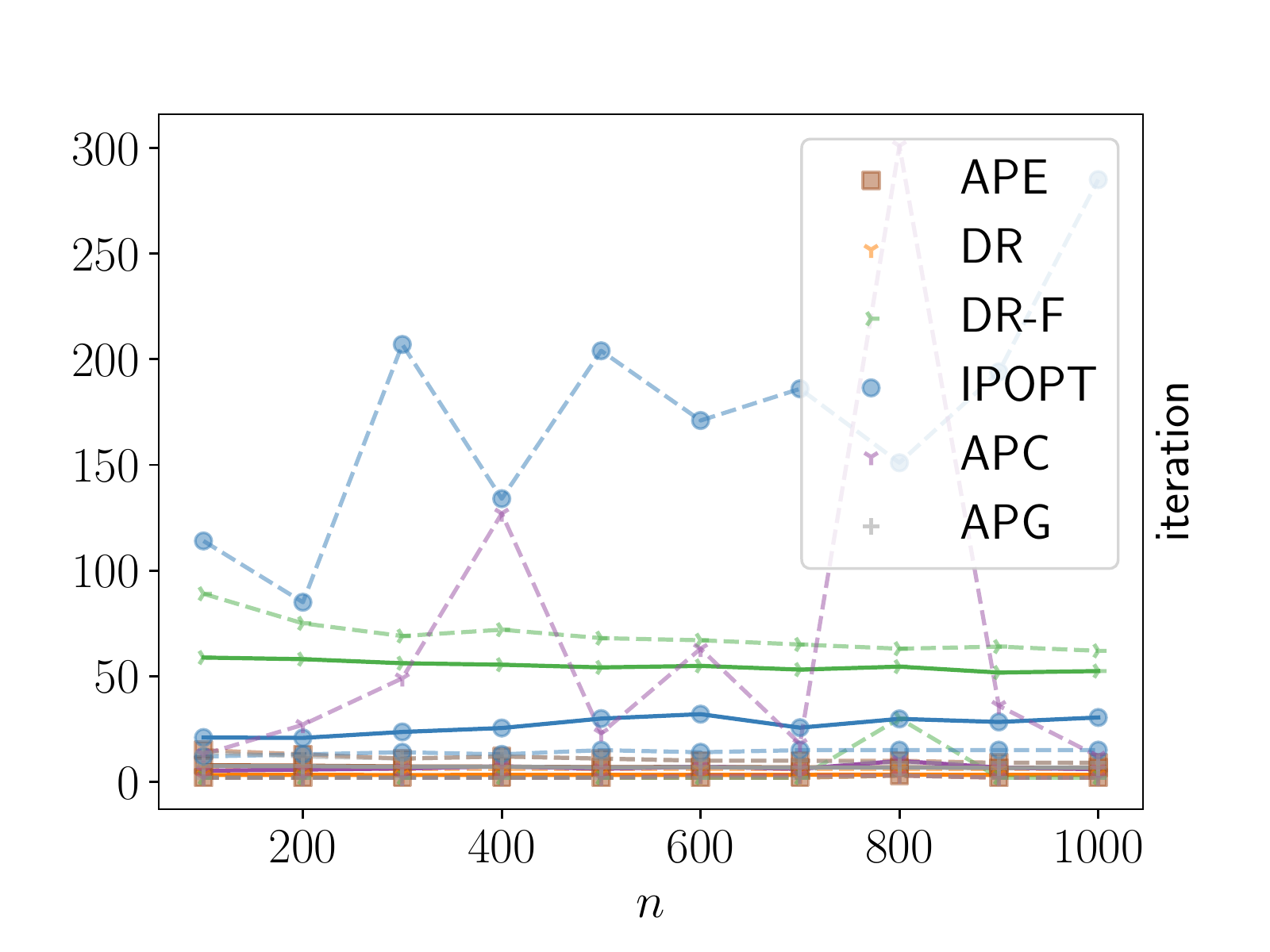}	
		\end{center}
		\caption{Number of iterations.}
		\label{subfig:hyperboloid_1000_timeout}
	\end{subfigure}
	\caption{Same as~\cref{fig:hyperboloid_100} for larger dimensions.}
	\label{fig:hyperboloid_1000}
\end{figure}

For the large-scale hyperboloidal problems ($n\geq 100$,~\cref{fig:hyperboloid_1000}), the best objectives are attained by APE and DR-F. The APG algorithm comes within one percent of their performance. The unmodified Douglas-Rachford finds objectives within several percent, and IPOPT within 10 percent, \eg, the mean objective for $n=1000$ (solid lines in~\cref{subfig:hyperboloid_1000_distance}) is around 0.51 for APE, DR-F and APG, around 0.53 for DR and around 0.63 for IPOPT. We note that the number of iterations increases with $n$, and that APG is also the fastest method. We also observe a significant increase in the execution time of IPOPT, which implies that the solving time is now larger than the 10 seconds that are required for creating the problem. Finally, we observe that both IPOPT and APC sometimes finish with a timeout, and return points above the prescribed deviation of $10^{-6}$.

\subsection{Alternating projections versus Gurobi}%
\label{sub:Alternate_projection_versus_Gurobi}

In this section, we benchmark the alternating projections with centre-based quasi-projection (APC) against an \emph{exact method}, which actually aims at finding the optimal solution of~\P. Several methods can be used to tackle \P, here we choose to use the commercial software Gurobi~\cite{gurobi}. This tool tackles the nonconvex quadratic equality via piecewise linearization, and solves the resulting mixed-integer quadratic programming (MIQP) problem. In this way, a solution along with a lower bound are obtained, and the optimal solution---up to a given tolerance---can be reached, assuming enough time is afforded to the solver.


The problem parameters are chosen so as to resemble problems from the power systems literature: the feasible set of the economic dispatch problem with power losses~\cite{vh22} typically exhibits a similar structure as the feasible set of~\cref{eq:P}. The entries of $\A$ are in the order of $10^{-5}$ except for diagonal entries ($10^{-4}$), $\Ab$ is close to 1 and $\Ac$ around -100. $\A$ represents the quadratic power losses---expected to be small---, and $\Ab$ encodes the constraint stating that the sum of the power production must be equal to the demand ($\approx-\Ac$).

The following quantities are compared:

\begin{align*}
	\textrm{Relative time} &= \log{\left (1 +\frac{t_{\alt} - t_{\exact}}{t_{\exact}}\right)} = \log \left (  \frac{\talt}{\texact}  \right ) \, \\ 
	\textrm{Relative distance} &= \log \left ( \frac{\dexact}{\dalt}  \right ) \,
\end{align*}
where $\talt$, $\texact$ are the execution times of APC and Gurobi, respectively, and $d$ is the distance between the final iterate and $\x^0$, \ie, the objective.

Note that, in order to smooth out random effects, we run $m = 100$ slightly different instances for each dimension and report the mean, such that, \eg, \[ \talt = \frac{\sum_{i=1}^m \talt^i}{m}, \]
where $\talt^i$ stands for the execution time of the $i-$th instance of the alternating projection method.

In the following experiments, it may be the case that no feasible point is found. For the alternating projection method, this can occur when the method reaches the maximum number of iterations, \eg, when the method is trapped in a cycle loop (see~\cref{fig:pathological_AP}). On the other hand, Gurobi may also fail to yield a feasible point, if the time limit criterion is attained. We do not encode such points into the relative distance.
In this way, we do not pollute the reported 
distance mean 
by a small number of instances that terminate due to a timeout. However, we also record the number of timeouts---either due to maximum iterations or time limit. 

\subsubsection{One-shot experiment}%
\label{ssub:One-shot_experiment}

In this experiment, we aim to compare the speed of both methods. We thus terminate the algorithm as soon as it finds a feasible point, hence the reference to ``one shot''. For APC, this does not affect the algorithm. On the other hand, Gurobi relies on lower and upper bounds, and terminates whenever a targeted tolerance is achieved. Here, we modify the stopping criterion such that the algorithm stops as soon as a feasible solution is obtained, no matter the objective. Hence, this is a lower bound on the execution time if the method is run with the tolerance criterion.

\begin{figure}
\begin{center}
	\includegraphics[width=\textwidth]{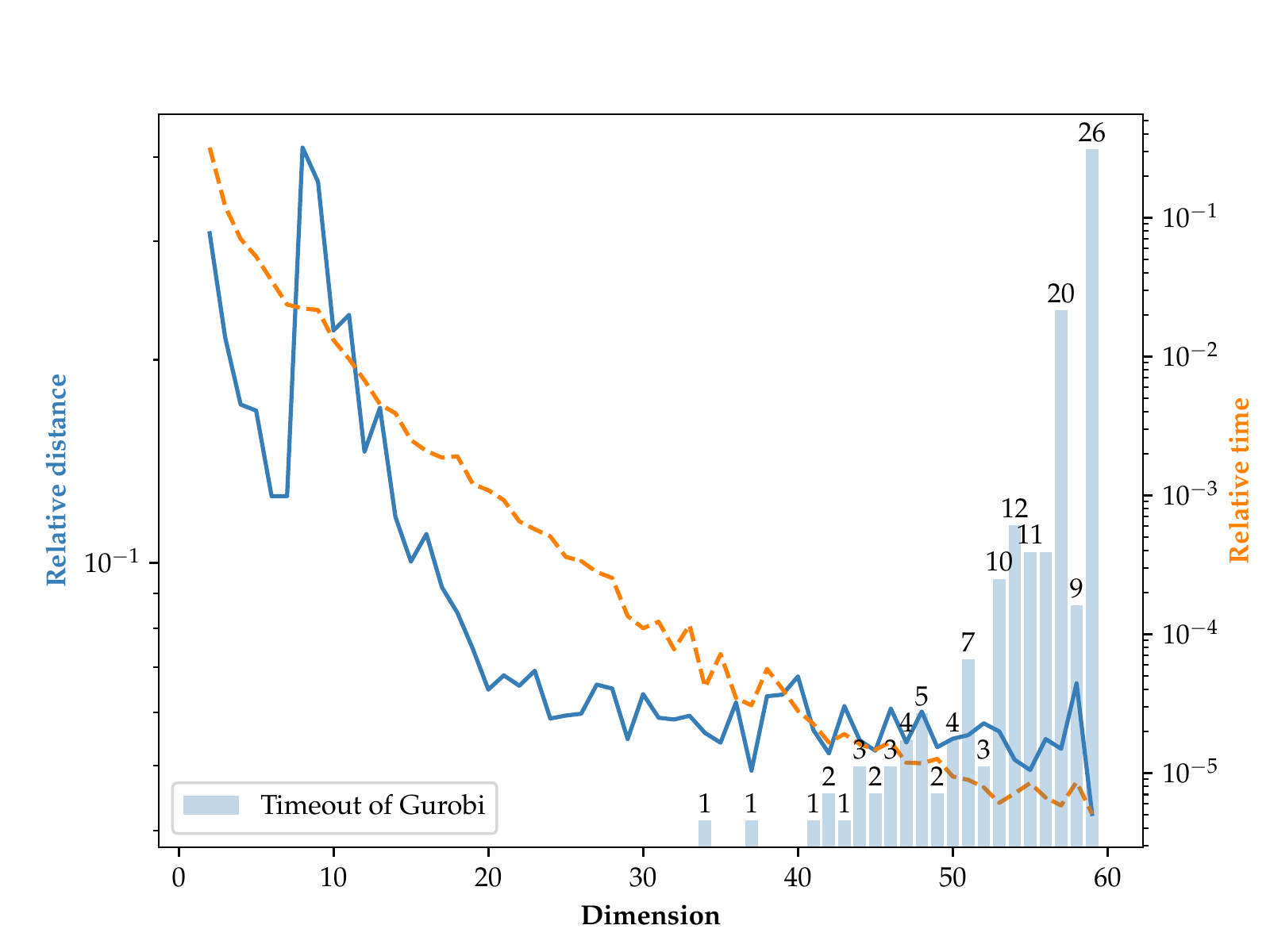}
\end{center}
\caption{Comparison between the alternating projection method with the centre-based quasi-projection (APC) and Gurobi. In this experiment, described in~\cref{ssub:One-shot_experiment}, the method terminates whenever it finds a feasible solution, no matter the objective. The timeout termination of Gurobi is set at 600 seconds, and the number of timeouts for the 100 instances is depicted as a bar plot.}
\label{fig:comparison_one_solution}
\end{figure}

\cref{fig:comparison_one_solution} presents the relative execution time and distance. We observe that APC clearly outperforms Gurobi in this experiment: for low dimension APC executes at least two times faster and reaches a better solution, for larger dimensions the difference becomes even larger, \eg, when the dimension is bigger than 40 APC accelerates by a factor of 100 000 and reaches an objective which is 10 times lower than that of Gurobi. Moreover, it should be noted that the number of timeout terminations recorded in Gurobi starts to increase for dimensions greater than 40 (see the bar plot in \cref{fig:comparison_one_solution}): hence the relative time is limited because of the time limit, this explains the saturation of the relative time for large dimensions. The relative distance does not encode the infeasibility of the points that finish with a timeout, and such points should have an infinite objective value.
The time limit criterion is set to 600 seconds. We note that a significant number of instances terminate without a solution for problems of large dimension.

\subsubsection{Multiple-shot experiment}%
\label{ssub:Multiple-shot_experiment}

\begin{figure}
\begin{center}
	\includegraphics[width=\textwidth]{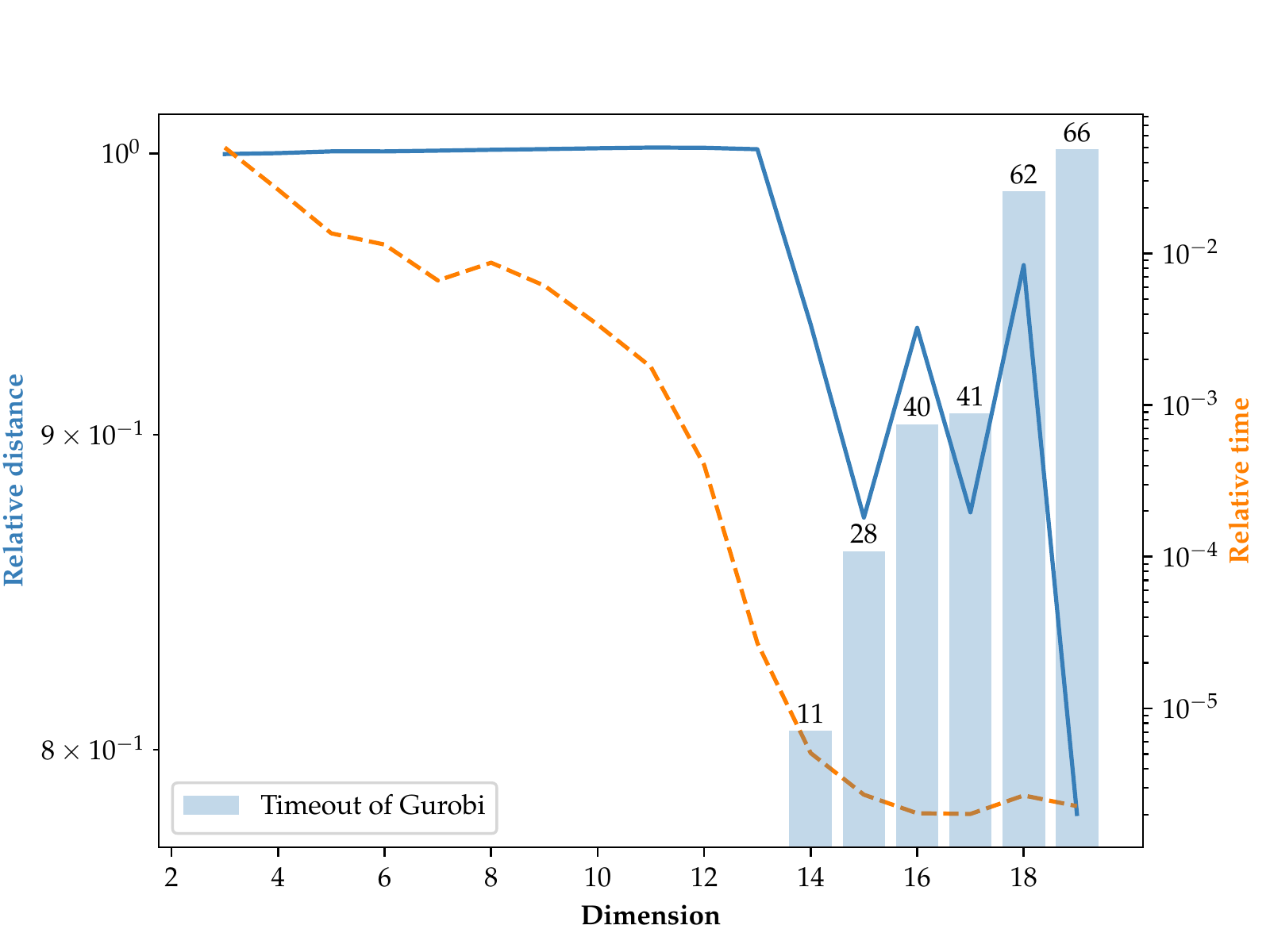}
\end{center}
\caption{Comparison between the alternating projection method with the centre-based quasi-projection (APC) and Gurobi. In this experiment, described in~\cref{ssub:Multiple-shot_experiment}, Gurobi terminates when the objective is proven to be optimal within a 1\% tolerance. The timeout termination of Gurobi is set to 600 seconds.}
	 \label{fig:Multiple-shot}
\end{figure}

In this experiment, we allow Gurobi to execute until it reaches the best solution, up to a given tolerance of one percent, or until timeout (600 seconds). \cref{fig:Multiple-shot} depicts the (mean) relative distance and execution time, for 100 runs, as well as the number of timeout terminations of Gurobi.
We observe that, for small problem instances, \ie, when the dimension is below 13, Gurobi reaches the best solution \emph{which is also very close to the one obtained via APC}. Indeed, a relative distance around one means that the solutions returned by the two algorithms are comparable. Since Gurobi does not terminate with a timeout, this implies that the solution returned by APC is, as a matter of fact, also optimal. Note that the theory does not guarantee that this should occur. 
We also note that the execution time of Gurobi is 10 to 1000 larger than that of APC.
For higher dimension, the relative distance slightly decreases and the relative time converges to $\num{1e-6}$: this is due to the increasing number of timeout terminations. In other words, Gurobi fails to find the best solution in an increasing execution time.

\section{Conclusion}%
\label{sec:Conclusion}

In this paper, the projection onto quadratic hypersurfaces, or quadrics, is investigated. We assume that the quadratic hypersurface is a non-cylindrical central quadric; however, the non-cylindrical assumption can be easily lifted by focusing on the variables that appear in the normal form. Using the method of Lagrange multipliers, we reduce this nonconvex optimization problem to the problem of finding the solutions of a system of nonlinear equations. We then show how one of the optimal solutions of the projection either lies in a finite set of computable solutions, or is a root of a scalar-valued nonlinear function. This unique root is located on a given interval, and is therefore simply computed with a Newton-Raphson scheme to which we provide suitable starting points. The cost of this projection is thus cheap, and the bottleneck is the eigenvector decomposition. This decomposition is needed for diagonalizing the matrix that is used to define the quadric.

We also propose a heuristic, referred to as quasi-projection, based on a geometric construction. This construction consists of finding the closest intersection between the quadric and a line passing by the point that we want to project. We detail two variants of the quasi-projection, depending on whether the direction of the line is computed as the level-curve gradient of the quadric, or the vector joining the centre and the point. This quasi-projection does not require eigenvector decomposition, thereby economizing in computational time.

This projection is then leveraged in the context of splitting algorithms, namely alternating projections and Douglas-Rachford splitting. This allows us to project a point onto a feasible set that is the intersection between a quadric and a box. The extension to the more general case of a Cartesian product of quadrics and a polytope is also discussed. Five methods are proposed depending on whether we use standard Douglas-Rachford splitting (DR), modified Douglas-Rachford splitting (DR-F), or one of the alternating projection methods. We detail the alternating projections with the exact projection on the quadric (APE) or one of the two quasi-projections (the centre-based APC or the gradient-based APG).

All methods are tested on problems of several dimensions, from 10 to 1000, and 100 independent trials are executed for each dimension. Using IPOPT as a benchmark, we find that APE and DF reach the best objectives and APG is within one percent, while APC and IPOPT lag behind. However, APG is much faster than the other methods and appears to achieve a good trade-off between the attained objective and execution time.

We also test APC on a case similar to the economic dispatch problem from the power systems literature, and compare it to Gurobi. We show that, in this specific case where the initial point is close to the feasible set, APC quickly reaches a solution close to or better than Gurobi, even if the execution time of Gurobi is several orders of magnitude greater. For small dimensions, Gurobi can guarantee the optimality of its solution, which shows that APC obtains the optimal solution in our examples. For higher dimension, Gurobi terminates with a timeout and with a higher objective than APC. Hence, even APC, which is the poorest of the methods that we propose in our paper in terms of performance, outperforms Gurobi in these experiments.

For the first part of the paper, namely the projection onto nonsingular quadrics in~\cref{sec:Projection_onto_a_quadric}, the extension to singular quadrics could be contemplated. In this context, the linear independence constraint qualification, LICQ, is not fulfilled any more. We should also include such points as projection candidates. A numerical comparison with the method from~\cite{sosa_algorithm_2020} is another natural research direction for further work.

For the second part of the paper,~\cref{sec:Splitting_methods}, further research may include the comparison with the alternating projection method using inexact projections, \eg, projecting on the tangent space of the previous feasible point. Such methods are discussed in~\cite{drusvyatskiy_local_2019}. 

\section*{Acknowledgement}
This work was supported by the Fonds de la Recherche Scientifique - FNRS under Grant no. PDR T.0025.18.

\bibliographystyle{plain}
\bibliography{bibliography}

\end{document}